
\input amstex

\expandafter\ifx\csname beta.def\endcsname\relax \else\endinput\fi
\expandafter\edef\csname beta.def\endcsname{%
 \catcode`\noexpand\@=\the\catcode`\@\space}

\let\atbefore @

\catcode`\@=11

\overfullrule\z@

\def\PaperA4{\hsize 6.25truein \vsize 9.63truein}

\def\PaperUS{\hsize 6.6truein \vsize 9truein} 

\def\foliorm{\ifMag\eightrm\else\ninerm\fi}

\let\@ft@\expandafter \let\@tb@f@\atbefore

\newif\ifMag
\def\Magset{\ifnum\mag>\@m\Magtrue\fi}
\Magset

\newif\ifUS

\newdimen\p@@ \p@@\p@
\def\m@ths@r{\ifnum\mathsurround=\z@\z@\else\maths@r\fi}
\def\maths@r{1.6\p@@} \def\mathsurzero{\def\maths@r{\z@}}

\mathsurround\maths@r
\font\Brm=cmr12 \font\Bbf=cmbx12 \font\Bit=cmti12 \font\ssf=cmss10
\font\Bsl=cmsl10 scaled 1200 \font\Bmmi=cmmi10 scaled 1200
\font\BBf=cmbx12 scaled 1200 \font\BMmi=cmmi10 scaled 1440

\def\atletter{\edef\atrestore{\catcode`\noexpand\@=\the\catcode`\@\space}
 \catcode`\@=11}

\newread\@ux \newwrite\@@x \newwrite\@@cd
\let\@np@@\input
\def\@np@t#1{\openin\@ux#1\relax\ifeof\@ux\else\closein\@ux\relax\@np@@ #1\fi}
\def\input#1 {\openin\@ux#1\relax\ifeof\@ux\wrs@x{No file #1}\else
 \closein\@ux\relax\@np@@ #1\fi}
\def\Input#1 {\relax} 

\def\wr@@x#1{} \def\wrs@x{\immediate\write\sixt@@n}

\def\readldf{\@np@t{\jobname.ldf}}
\def\writeldf{\def\wr@@x{\immediate\write\@@x}
 \def\cl@selbl{\wr@@x{\string\Snodef{\the\Sno}}\wr@@x{\string\endinput}%
 \immediate\closeout\@@x} \immediate\openout\@@x\jobname.ldf}
\let\cl@selbl\relax

\def\tod@y{\ifcase\month\or
 January\or February\or March\or April\or May\or June\or July\or
 August\or September\or October\or November\or December\fi\space\,
\number\day,\space\,\number\year}

\newcount\c@time
\def\h@@r{hh}\def\m@n@te{mm}
\def\wh@tt@me{\c@time\time\divide\c@time 60\edef\h@@r{\number\c@time}%
 \multiply\c@time -60\advance\c@time\time\edef
 \m@n@te{\ifnum\c@time<10 0\fi\number\c@time}}
\def\t@me{\h@@r\/{\rm:}\m@n@te}  \let\whattime\wh@tt@me
\def\today{\tod@y\wr@@x{\string\todaydef{\tod@y}}}
\def\nowtime{\t@me{\let\/\ic@\wr@@x{\string\nowtimedef{\t@me}}}}
\def\todaydef#1{} \def\nowtimedef#1{}

\def\em#1{{\it #1\/}} \def\emph#1{{\sl #1\/}}
\def\itemhalf#1{\par\hangindent1.5\parindent
 \hglue-.5\parindent\textindent{\rm#1}}
\def\fitem#1{\par\setbox\z@\hbox{#1}\hangindent\wd\z@
 \hglue-2\parindent\kern\wd\z@\indent\llap{#1}\ignore}

\def\itemflat#1{\par\setbox\z@\hbox{\rm #1\enspace}\hang\ifnum\wd\z@>\parindent
 \noindent\unhbox\z@\ignore\else\textindent{\rm#1}\fi}

\newcount\itemlet
\def\newbi{\itemlet 96} \newbi
\def\bitem{\gad\itemlet \par\hangindent1.5\parindent
 \hglue-.5\parindent\textindent{\rm\rlap{\char\the\itemlet}\hp{b})}}
\def\atem{\newbi\bitem}

\newcount\itemrm

\def\iitem{\gad\itemrm \par\hangindent1.5\parindent
 \hglue-.5\parindent\textindent{\rm\hp{v}\llap{\romannumeral\the\itemrm})}}

\def\center{\par\begingroup\leftskip\z@ plus \hsize \rightskip\leftskip
 \parindent\z@\parfillskip\z@skip \def\\{\unskip\break}}
\def\endcenter{\endgraf\endgroup}

\def\Abstract{\begingroup\narrower\nt{\bf Abstract.}\enspace\ignore}
\def\endAbs{\endgraf\endgroup}

\let\b@gr@@\begingroup \let\B@gr@@\begingroup
\def\b@gr@{\b@gr@@\let\b@gr@@\undefined}
\def\B@gr@{\B@gr@@\let\B@gr@@\undefined}

\def\@fn@xt#1#2#3{\let\@ch@r=#1\def\n@xt{\ifx\t@st@\@ch@r
 \def\n@@xt{#2}\else\def\n@@xt{#3}\fi\n@@xt}\futurelet\t@st@\n@xt}

\def\@fwd@@#1#2#3{\setbox\z@\hbox{#1}\ifdim\wd\z@>\z@#2\else#3\fi}
\def\s@twd@#1#2{\setbox\z@\hbox{#2}#1\wd\z@}

\def\r@st@re#1{\let#1\s@v@} \def\s@v@d@f{\let\s@v@}

\def\p@sk@p#1#2{\par\skip@#2\relax\ifdim\lastskip<\skip@\relax\removelastskip
 \ifnum#1=\z@\else\penalty#1\relax\fi\vskip\skip@
 \else\ifnum#1=\z@\else\penalty#1\relax\fi\fi}
\def\sk@@p#1{\par\skip@#1\relax\ifdim\lastskip<\skip@\relax\removelastskip
 \vskip\skip@\fi}

\newbox\p@b@ld
\def\poorbold#1{\setbox\p@b@ld\hbox{#1}\kern-.01em\copy\p@b@ld\kern-\wd\p@b@ld
 \kern.02em\copy\p@b@ld\kern-\wd\p@b@ld\kern-.012em\raise.02em\box\p@b@ld}

\ifx\plainfootnote\undefined \let\plainfootnote\footnote \fi

\let\s@v@\proclaim \let\proclaim\relax
\def\r@R@fs#1{\let#1\s@R@fs} \let\s@R@fs\Refs \let\Refs\relax
\def\r@endd@#1{\let#1\s@endd@} \let\s@endd@\enddocument
\let\bye\relax

\def\myR@fs{\@fn@xt[\m@R@f@\m@R@fs} \def\m@R@fs{\@fn@xt*\m@r@f@@\m@R@f@@}
\def\m@R@f@@{\m@R@f@[References]} \def\m@r@f@@*{\m@R@f@[]}

\def\Twelvepoint{\twelvepoint \let\Bbf\BBf \let\Bmmi\BMmi
\font\Brm=cmr12 scaled 1200 \font\Bit=cmti12 scaled 1200
\font\ssf=cmss10 scaled 1200 \font\Bsl=cmsl10 scaled 1440
\font\BBf=cmbx12 scaled 1440 \font\BMmi=cmmi10 scaled 1728}

\newdimen\r@f@nd \newbox\r@f@b@x \newbox\adjb@x
\newbox\p@nct@ \newbox\k@yb@x \newcount\rcount
\newbox\b@b@x \newbox\p@p@rb@x \newbox\j@@rb@x \newbox\y@@rb@x
\newbox\v@lb@x \newbox\is@b@x \newbox\p@g@b@x \newif\ifp@g@ \newif\ifp@g@s
\newbox\inb@@kb@x \newbox\b@@kb@x \newbox\p@blb@x \newbox\p@bl@db@x
\newbox\ed@b@x \newif\ifed@ \newif\ifed@s \newif\if@fl@b \newif\if@fn@m
\newbox\p@p@nf@b@x \newbox\inf@b@x \newbox\b@@nf@b@x
\newtoks\@dd@p@n \newtoks\@ddt@ks

\newdimen\b@gsize

\newif\ifp@gen@

\newif\ifamsppt

\@ft@\ifx\csname amsppt.sty\endcsname\relax

\def\p@@nt{.\kern.3em} \let\point\p@@nt

\let\proheadfont\bf \let\probodyfont\sl \let\demofont\it

\def\reffont#1{\def\r@ff@nt{#1}} \reffont\rm
\def\keyfont#1{\def\k@yf@nt{#1}} \keyfont\rm
\def\paperfont#1{\def\p@p@rf@nt{#1}} \paperfont\it
\def\bookfont#1{\def\b@@kf@nt{#1}} \bookfont\it
\def\volfont#1{\def\v@lf@nt{#1}} \volfont\bf
\def\issuefont#1{\def\iss@f@nt{#1}} \issuefont{no\p@@nt}

\headline={\hfil}
\footline={\ifp@gen@\ifnum\pageno=\z@\else\hfil\foliorm\folio\fi\else
 \ifnum\pageno=\z@\hfil\foliorm\folio\fi\fi\hfil\global\p@gen@true}
\parindent1pc

\font@\tensmc=cmcsc10
\font@\sevenex=cmex7
\font@\sevenit=cmti7
\font@\eightrm=cmr8
\font@\sixrm=cmr6
\font@\eighti=cmmi8 \skewchar\eighti='177
\font@\sixi=cmmi6 \skewchar\sixi='177
\font@\eightsy=cmsy8 \skewchar\eightsy='60
\font@\sixsy=cmsy6 \skewchar\sixsy='60
\font@\eightex=cmex8
\font@\eightbf=cmbx8
\font@\sixbf=cmbx6
\font@\eightit=cmti8
\font@\eightsl=cmsl8
\font@\eightsmc=cmcsc8
\font@\eighttt=cmtt8
\font@\ninerm=cmr9
\font@\ninei=cmmi9 \skewchar\ninei='177
\font@\ninesy=cmsy9 \skewchar\ninesy='60
\font@\nineex=cmex9
\font@\ninebf=cmbx9
\font@\nineit=cmti9
\font@\ninesl=cmsl9
\font@\ninesmc=cmcsc9
\font@\ninemsa=msam9
\font@\ninemsb=msbm9
\font@\nineeufm=eufm9
\font@\eightmsa=msam8
\font@\eightmsb=msbm8
\font@\eighteufm=eufm8
\font@\sixmsa=msam6
\font@\sixmsb=msbm6
\font@\sixeufm=eufm6

\loadmsam\loadmsbm\loadeufm
\input amssym.tex

\def\footnoterule{\kern-3\p@\hrule width5pc\kern 2.6\p@}
\def\m@k@foot#1{\insert\footins
 {\interlinepenalty\interfootnotelinepenalty
 \eightpoint\splittopskip\ht\strutbox\splitmaxdepth\dp\strutbox
 \floatingpenalty\@MM\leftskip\z@\rightskip\z@
 \spaceskip\z@\xspaceskip\z@
 \leavevmode\footstrut\ignore#1\unskip\lower\dp\strutbox
 \vbox to\dp\strutbox{}}}
\def\ftext#1{\m@k@foot{\vsk-.8>\nt #1}}
\def\pr@cl@@m#1{\p@sk@p{-100}\medskipamount
 \def\endproclaim{\endgroup\p@sk@p{55}\medskipamount}\begingroup
 \nt\ignore\proheadfont#1\unskip.\enspace\probodyfont\ignore}
\outer\def\proclaim{\pr@cl@@m} \s@v@d@f\proclaim \let\proclaim\relax
\def\demo#1{\sk@@p\medskipamount\nt{\ignore\demofont#1\unskip.}\enspace
 \ignore}
\def\enddemo{\sk@@p\medskipamount}

\def\cite#1{{\rm[#1]}} \let\nofrills\relax
 \def\Refs#1#2{\relax}

\def\big@#1#2{{\hbox{$\left#2\vcenter to#1\b@gsize{}%
 \right.\nulldelimiterspace\z@\m@th$}}}
\def\big{\big@\@ne}
\def\Big{\big@{1.5}}
\def\bigg{\big@\tw@}
\def\Bigg{\big@{2.5}}
\normallineskiplimit\p@

\def\tenpoint{\p@@\p@ \normallineskiplimit\p@@
 \mathsurround\m@ths@r \normalbaselineskip12\p@@
 \abovedisplayskip12\p@@ plus3\p@@ minus9\p@@
 \belowdisplayskip\abovedisplayskip
 \abovedisplayshortskip\z@ plus3\p@@
 \belowdisplayshortskip7\p@@ plus3\p@@ minus4\p@@
 \textonlyfont@\rm\tenrm \textonlyfont@\it\tenit
 \textonlyfont@\sl\tensl \textonlyfont@\bf\tenbf
 \textonlyfont@\smc\tensmc \textonlyfont@\tt\tentt
 \ifsyntax@ \def\big##1{{\hbox{$\left##1\right.$}}}%
  \let\Big\big \let\bigg\big \let\Bigg\big
 \else
  \textfont\z@\tenrm \scriptfont\z@\sevenrm \scriptscriptfont\z@\fiverm
  \textfont\@ne\teni \scriptfont\@ne\seveni \scriptscriptfont\@ne\fivei
  \textfont\tw@\tensy \scriptfont\tw@\sevensy \scriptscriptfont\tw@\fivesy
  \textfont\thr@@\tenex \scriptfont\thr@@\sevenex
	\scriptscriptfont\thr@@\sevenex
  \textfont\itfam\tenit \scriptfont\itfam\sevenit
	\scriptscriptfont\itfam\sevenit
  \textfont\bffam\tenbf \scriptfont\bffam\sevenbf
	\scriptscriptfont\bffam\fivebf
  \textfont\msafam\tenmsa \scriptfont\msafam\sevenmsa
	\scriptscriptfont\msafam\fivemsa
  \textfont\msbfam\tenmsb \scriptfont\msbfam\sevenmsb
	\scriptscriptfont\msbfam\fivemsb
  \textfont\eufmfam\teneufm \scriptfont\eufmfam\seveneufm
	\scriptscriptfont\eufmfam\fiveeufm
  \setbox\strutbox\hbox{\vrule height8.5\p@@ depth3.5\p@@ width\z@}%
  \setbox\strutbox@\hbox{\lower.5\normallineskiplimit\vbox{%
	\kern-\normallineskiplimit\copy\strutbox}}%
   \setbox\z@\vbox{\hbox{$($}\kern\z@}\b@gsize1.2\ht\z@
  \fi
  \normalbaselines\rm\dotsspace@1.5mu\ex@.2326ex\jot3\ex@}

\def\eightpoint{\p@@.8\p@ \normallineskiplimit\p@@
 \mathsurround\m@ths@r \normalbaselineskip10\p@
 \abovedisplayskip10\p@ plus2.4\p@ minus7.2\p@
 \belowdisplayskip\abovedisplayskip
 \abovedisplayshortskip\z@ plus3\p@@
 \belowdisplayshortskip7\p@@ plus3\p@@ minus4\p@@
 \textonlyfont@\rm\eightrm \textonlyfont@\it\eightit
 \textonlyfont@\sl\eightsl \textonlyfont@\bf\eightbf
 \textonlyfont@\smc\eightsmc \textonlyfont@\tt\eighttt
 \ifsyntax@\def\big##1{{\hbox{$\left##1\right.$}}}%
  \let\Big\big \let\bigg\big \let\Bigg\big
 \else
  \textfont\z@\eightrm \scriptfont\z@\sixrm \scriptscriptfont\z@\fiverm
  \textfont\@ne\eighti \scriptfont\@ne\sixi \scriptscriptfont\@ne\fivei
  \textfont\tw@\eightsy \scriptfont\tw@\sixsy \scriptscriptfont\tw@\fivesy
  \textfont\thr@@\eightex \scriptfont\thr@@\sevenex
	\scriptscriptfont\thr@@\sevenex
  \textfont\itfam\eightit \scriptfont\itfam\sevenit
	\scriptscriptfont\itfam\sevenit
  \textfont\bffam\eightbf \scriptfont\bffam\sixbf
	\scriptscriptfont\bffam\fivebf
  \textfont\msafam\eightmsa \scriptfont\msafam\sixmsa
	\scriptscriptfont\msafam\fivemsa
  \textfont\msbfam\eightmsb \scriptfont\msbfam\sixmsb
	\scriptscriptfont\msbfam\fivemsb
  \textfont\eufmfam\eighteufm \scriptfont\eufmfam\sixeufm
	\scriptscriptfont\eufmfam\fiveeufm
 \setbox\strutbox\hbox{\vrule height7\p@ depth3\p@ width\z@}%
 \setbox\strutbox@\hbox{\raise.5\normallineskiplimit\vbox{%
   \kern-\normallineskiplimit\copy\strutbox}}%
 \setbox\z@\vbox{\hbox{$($}\kern\z@}\b@gsize1.2\ht\z@
 \fi
 \normalbaselines\eightrm\dotsspace@1.5mu\ex@.2326ex\jot3\ex@}

\def\ninepoint{\p@@.9\p@ \normallineskiplimit\p@@
 \mathsurround\m@ths@r \normalbaselineskip11\p@
 \abovedisplayskip11\p@ plus2.7\p@ minus8.1\p@
 \belowdisplayskip\abovedisplayskip
 \abovedisplayshortskip\z@ plus3\p@@
 \belowdisplayshortskip7\p@@ plus3\p@@ minus4\p@@
 \textonlyfont@\rm\ninerm \textonlyfont@\it\nineit
 \textonlyfont@\sl\ninesl \textonlyfont@\bf\ninebf
 \textonlyfont@\smc\ninesmc \textonlyfont@\tt\ninett
 \ifsyntax@ \def\big##1{{\hbox{$\left##1\right.$}}}%
  \let\Big\big \let\bigg\big \let\Bigg\big
 \else
  \textfont\z@\ninerm \scriptfont\z@\sevenrm \scriptscriptfont\z@\fiverm
  \textfont\@ne\ninei \scriptfont\@ne\seveni \scriptscriptfont\@ne\fivei
  \textfont\tw@\ninesy \scriptfont\tw@\sevensy \scriptscriptfont\tw@\fivesy
  \textfont\thr@@\nineex \scriptfont\thr@@\sevenex
	\scriptscriptfont\thr@@\sevenex
  \textfont\itfam\nineit \scriptfont\itfam\sevenit
	\scriptscriptfont\itfam\sevenit
  \textfont\bffam\ninebf \scriptfont\bffam\sevenbf
	\scriptscriptfont\bffam\fivebf
  \textfont\msafam\ninemsa \scriptfont\msafam\sevenmsa
	\scriptscriptfont\msafam\fivemsa
  \textfont\msbfam\ninemsb \scriptfont\msbfam\sevenmsb
	\scriptscriptfont\msbfam\fivemsb
  \textfont\eufmfam\nineeufm \scriptfont\eufmfam\seveneufm
	\scriptscriptfont\eufmfam\fiveeufm
  \setbox\strutbox\hbox{\vrule height8.5\p@@ depth3.5\p@@ width\z@}%
  \setbox\strutbox@\hbox{\lower.5\normallineskiplimit\vbox{%
	\kern-\normallineskiplimit\copy\strutbox}}%
   \setbox\z@\vbox{\hbox{$($}\kern\z@}\b@gsize1.2\ht\z@
  \fi
  \normalbaselines\rm\dotsspace@1.5mu\ex@.2326ex\jot3\ex@}

\font@\twelverm=cmr10 scaled 1200
\font@\twelveit=cmti10 scaled 1200
\font@\twelvesl=cmsl10 scaled 1200
\font@\twelvebf=cmbx10 scaled 1200
\font@\twelvesmc=cmcsc10 scaled 1200
\font@\twelvett=cmtt10 scaled 1200
\font@\twelvei=cmmi10 scaled 1200 \skewchar\twelvei='177
\font@\twelvesy=cmsy10 scaled 1200 \skewchar\twelvesy='60
\font@\twelveex=cmex10 scaled 1200
\font@\twelvemsa=msam10 scaled 1200
\font@\twelvemsb=msbm10 scaled 1200
\font@\twelveeufm=eufm10 scaled 1200

\def\twelvepoint{\p@@1.2\p@ \normallineskiplimit\p@@
 \mathsurround\m@ths@r \normalbaselineskip12\p@@
 \abovedisplayskip12\p@@ plus3\p@@ minus9\p@@
 \belowdisplayskip\abovedisplayskip
 \abovedisplayshortskip\z@ plus3\p@@
 \belowdisplayshortskip7\p@@ plus3\p@@ minus4\p@@
 \textonlyfont@\rm\twelverm \textonlyfont@\it\twelveit
 \textonlyfont@\sl\twelvesl \textonlyfont@\bf\twelvebf
 \textonlyfont@\smc\twelvesmc \textonlyfont@\tt\twelvett
 \ifsyntax@ \def\big##1{{\hbox{$\left##1\right.$}}}%
  \let\Big\big \let\bigg\big \let\Bigg\big
 \else
  \textfont\z@\twelverm \scriptfont\z@\eightrm \scriptscriptfont\z@\sixrm
  \textfont\@ne\twelvei \scriptfont\@ne\eighti \scriptscriptfont\@ne\sixi
  \textfont\tw@\twelvesy \scriptfont\tw@\eightsy \scriptscriptfont\tw@\sixsy
  \textfont\thr@@\twelveex \scriptfont\thr@@\eightex
	\scriptscriptfont\thr@@\sevenex
  \textfont\itfam\twelveit \scriptfont\itfam\eightit
	\scriptscriptfont\itfam\sevenit
  \textfont\bffam\twelvebf \scriptfont\bffam\eightbf
	\scriptscriptfont\bffam\sixbf
  \textfont\msafam\twelvemsa \scriptfont\msafam\eightmsa
	\scriptscriptfont\msafam\sixmsa
  \textfont\msbfam\twelvemsb \scriptfont\msbfam\eightmsb
	\scriptscriptfont\msbfam\sixmsb
  \textfont\eufmfam\twelveeufm \scriptfont\eufmfam\eighteufm
	\scriptscriptfont\eufmfam\sixeufm
  \setbox\strutbox\hbox{\vrule height8.5\p@@ depth3.5\p@@ width\z@}%
  \setbox\strutbox@\hbox{\lower.5\normallineskiplimit\vbox{%
	\kern-\normallineskiplimit\copy\strutbox}}%
  \setbox\z@\vbox{\hbox{$($}\kern\z@}\b@gsize1.2\ht\z@
  \fi
  \normalbaselines\rm\dotsspace@1.5mu\ex@.2326ex\jot3\ex@}

\font@\twelvetrm=cmr10 at 12truept
\font@\twelvetit=cmti10 at 12truept
\font@\twelvetsl=cmsl10 at 12truept
\font@\twelvetbf=cmbx10 at 12truept
\font@\twelvetsmc=cmcsc10 at 12truept
\font@\twelvettt=cmtt10 at 12truept
\font@\twelveti=cmmi10 at 12truept \skewchar\twelveti='177
\font@\twelvetsy=cmsy10 at 12truept \skewchar\twelvetsy='60
\font@\twelvetex=cmex10 at 12truept
\font@\twelvetmsa=msam10 at 12truept
\font@\twelvetmsb=msbm10 at 12truept
\font@\twelveteufm=eufm10 at 12truept

\def\twelvetruepoint{\p@@1.2truept \normallineskiplimit\p@@
 \mathsurround\m@ths@r \normalbaselineskip12\p@@
 \abovedisplayskip12\p@@ plus3\p@@ minus9\p@@
 \belowdisplayskip\abovedisplayskip
 \abovedisplayshortskip\z@ plus3\p@@
 \belowdisplayshortskip7\p@@ plus3\p@@ minus4\p@@
 \textonlyfont@\rm\twelvetrm \textonlyfont@\it\twelvetit
 \textonlyfont@\sl\twelvetsl \textonlyfont@\bf\twelvetbf
 \textonlyfont@\smc\twelvetsmc \textonlyfont@\tt\twelvettt
 \ifsyntax@ \def\big##1{{\hbox{$\left##1\right.$}}}%
  \let\Big\big \let\bigg\big \let\Bigg\big
 \else
  \textfont\z@\twelvetrm \scriptfont\z@\eightrm \scriptscriptfont\z@\sixrm
  \textfont\@ne\twelveti \scriptfont\@ne\eighti \scriptscriptfont\@ne\sixi
  \textfont\tw@\twelvetsy \scriptfont\tw@\eightsy \scriptscriptfont\tw@\sixsy
  \textfont\thr@@\twelvetex \scriptfont\thr@@\eightex
	\scriptscriptfont\thr@@\sevenex
  \textfont\itfam\twelvetit \scriptfont\itfam\eightit
	\scriptscriptfont\itfam\sevenit
  \textfont\bffam\twelvetbf \scriptfont\bffam\eightbf
	\scriptscriptfont\bffam\sixbf
  \textfont\msafam\twelvetmsa \scriptfont\msafam\eightmsa
	\scriptscriptfont\msafam\sixmsa
  \textfont\msbfam\twelvetmsb \scriptfont\msbfam\eightmsb
	\scriptscriptfont\msbfam\sixmsb
  \textfont\eufmfam\twelveteufm \scriptfont\eufmfam\eighteufm
	\scriptscriptfont\eufmfam\sixeufm
  \setbox\strutbox\hbox{\vrule height8.5\p@@ depth3.5\p@@ width\z@}%
  \setbox\strutbox@\hbox{\lower.5\normallineskiplimit\vbox{%
	\kern-\normallineskiplimit\copy\strutbox}}%
  \setbox\z@\vbox{\hbox{$($}\kern\z@}\b@gsize1.2\ht\z@
  \fi
  \normalbaselines\rm\dotsspace@1.5mu\ex@.2326ex\jot3\ex@}

\font@\elevenrm=cmr10 scaled 1095
\font@\elevenit=cmti10 scaled 1095
\font@\elevensl=cmsl10 scaled 1095
\font@\elevenbf=cmbx10 scaled 1095
\font@\elevensmc=cmcsc10 scaled 1095
\font@\eleventt=cmtt10 scaled 1095
\font@\eleveni=cmmi10 scaled 1095 \skewchar\eleveni='177
\font@\elevensy=cmsy10 scaled 1095 \skewchar\elevensy='60
\font@\elevenex=cmex10 scaled 1095
\font@\elevenmsa=msam10 scaled 1095
\font@\elevenmsb=msbm10 scaled 1095
\font@\eleveneufm=eufm10 scaled 1095

\def\elevenpoint{\p@@1.1\p@ \normallineskiplimit\p@@
 \mathsurround\m@ths@r \normalbaselineskip12\p@@
 \abovedisplayskip12\p@@ plus3\p@@ minus9\p@@
 \belowdisplayskip\abovedisplayskip
 \abovedisplayshortskip\z@ plus3\p@@
 \belowdisplayshortskip7\p@@ plus3\p@@ minus4\p@@
 \textonlyfont@\rm\elevenrm \textonlyfont@\it\elevenit
 \textonlyfont@\sl\elevensl \textonlyfont@\bf\elevenbf
 \textonlyfont@\smc\elevensmc \textonlyfont@\tt\eleventt
 \ifsyntax@ \def\big##1{{\hbox{$\left##1\right.$}}}%
  \let\Big\big \let\bigg\big \let\Bigg\big
 \else
  \textfont\z@\elevenrm \scriptfont\z@\eightrm \scriptscriptfont\z@\sixrm
  \textfont\@ne\eleveni \scriptfont\@ne\eighti \scriptscriptfont\@ne\sixi
  \textfont\tw@\elevensy \scriptfont\tw@\eightsy \scriptscriptfont\tw@\sixsy
  \textfont\thr@@\elevenex \scriptfont\thr@@\eightex
	\scriptscriptfont\thr@@\sevenex
  \textfont\itfam\elevenit \scriptfont\itfam\eightit
	\scriptscriptfont\itfam\sevenit
  \textfont\bffam\elevenbf \scriptfont\bffam\eightbf
	\scriptscriptfont\bffam\sixbf
  \textfont\msafam\elevenmsa \scriptfont\msafam\eightmsa
	\scriptscriptfont\msafam\sixmsa
  \textfont\msbfam\elevenmsb \scriptfont\msbfam\eightmsb
	\scriptscriptfont\msbfam\sixmsb
  \textfont\eufmfam\eleveneufm \scriptfont\eufmfam\eighteufm
	\scriptscriptfont\eufmfam\sixeufm
  \setbox\strutbox\hbox{\vrule height8.5\p@@ depth3.5\p@@ width\z@}%
  \setbox\strutbox@\hbox{\lower.5\normallineskiplimit\vbox{%
	\kern-\normallineskiplimit\copy\strutbox}}%
  \setbox\z@\vbox{\hbox{$($}\kern\z@}\b@gsize1.2\ht\z@
  \fi
  \normalbaselines\rm\dotsspace@1.5mu\ex@.2326ex\jot3\ex@}

\def\m@R@f@[#1]{\mathsurzero{
 \s@ct{}{#1}}\wr@@c{\string\Refcd{#1}{\the\pageno}}\B@gr@
 \frenchspacing\rcount\z@\refkey{\k@yf@nt[##1]}\refno{\k@yf@nt[##1]}%
 \widest{AZ}\keyright\let\Key\key\let\refin\relax}
\def\widest#1{\s@twd@\r@f@nd{\r@fk@y{\k@yf@nt#1}\enspace}}
\def\widestno#1{\s@twd@\r@f@nd{\r@fn@{\k@yf@nt#1}\enspace}}
\def\widestlabel#1{\s@twd@\r@f@nd{\k@yf@nt#1\enspace}}
\def\refkey{\def\r@fk@y##1} \def\refno{\def\r@fn@##1}
\def\keyright{\def\r@fit@m{\hang\textindent}}
\def\keyflat{\def\r@fit@m##1{\setbox\z@\hbox{##1\enspace}\hang\noindent
 \ifnum\wd\z@<\parindent\indent\hglue-\wd\z@\fi\unhbox\z@}}

\def\R@fb@x{\global\setbox\r@f@b@x} \def\K@yb@x{\global\setbox\k@yb@x}
\def\ref{\par\b@gr@\r@ff@nt\R@fb@x\box\voidb@x\K@yb@x\box\voidb@x
 \@fn@mfalse\@fl@bfalse\b@g@nr@f}
\def\c@nc@t#1{\setbox\z@\lastbox
 \setbox\adjb@x\hbox{\unhbox\adjb@x\unhbox\z@\unskip\unskip\unpenalty#1}}
\def\adjust#1{\relax\ifmmode\penalty-\@M\null\hfil$\clubpenalty\z@
 \widowpenalty\z@\interlinepenalty\z@\offinterlineskip\endgraf
 \setbox\z@\lastbox\unskip\unpenalty\c@nc@t{#1}\nt$\hfil\penalty-\@M
 \else\endgraf\c@nc@t{#1}\nt\fi}
\def\adjustnext#1{\P@nct\hbox{#1}\ignore}
\def\adjustend#1{\def\@djp@{#1}\ignore}
\def\addtoks#1{\global\@ddt@ks{#1}\ignore}
\def\addnext#1{\global\@dd@p@n{#1}\ignore}

\def\cl@s@{\adjust{\@djp@}\endgraf\setbox\z@\lastbox
 \global\setbox\@ne\hbox{\unhbox\adjb@x\ifvoid\z@\else\unhbox\z@\unskip\unskip
 \unpenalty\fi}\egroup\ifnum\c@rr@nt=\k@yb@x\global\fi
 \setbox\c@rr@nt\hbox{\unhbox\@ne\box\p@nct@}\P@nct\null
 \the\@ddt@ks\global\@ddt@ks{}}
\def\@p@n#1{\def\c@rr@nt{#1}\setbox\c@rr@nt\vbox\bgroup\let\@djp@\relax
 \hsize\maxdimen\nt\the\@dd@p@n\global\@dd@p@n{}}
\def\b@g@nr@f{\bgroup\@p@n\z@}
\def\key{\cl@s@\ifvoid\k@yb@x\@p@n\k@yb@x\k@yf@nt\else\@p@n\z@\fi}
\def\label{\cl@s@\ifvoid\k@yb@x\global\@fl@btrue\@p@n\k@yb@x\k@yf@nt\else
 \@p@n\z@\fi}
\def\no{\cl@s@\ifvoid\k@yb@x\gad\rcount\global\@fn@mtrue
 \K@yb@x\hbox{\k@yf@nt\the\rcount}\fi\@p@n\z@}
\def\labelno{\cl@s@\ifvoid\k@yb@x\gad\rcount\@fl@btrue
 \@p@n\k@yb@x\k@yf@nt\the\rcount\else\@p@n\z@\fi}
\def\by{\cl@s@\@p@n\b@b@x} \def\paper{\cl@s@\@p@n\p@p@rb@x\p@p@rf@nt\ignore}
\def\jour{\cl@s@\@p@n\j@@rb@x} \def\yr{\cl@s@\@p@n\y@@rb@x}
\def\vol{\cl@s@\@p@n\v@lb@x\v@lf@nt\ignore}
\def\issue{\cl@s@\@p@n\is@b@x\iss@f@nt\ignore}
\def\page{\cl@s@\ifp@g@s\@p@n\z@\else\p@g@true\@p@n\p@g@b@x\fi}
\def\pages{\cl@s@\ifp@g@\@p@n\z@\else\p@g@strue\@p@n\p@g@b@x\fi}
\def\inbook{\cl@s@\@p@n\inb@@kb@x}
\def\book{\cl@s@\@p@n\b@@kb@x\b@@kf@nt\ignore}
\def\publ{\cl@s@\@p@n\p@blb@x} \def\publaddr{\cl@s@\@p@n\p@bl@db@x}
\def\ed{\cl@s@\ifed@s\@p@n\z@\else\ed@true\@p@n\ed@b@x\fi}
\def\eds{\cl@s@\ifed@\@p@n\z@\else\ed@strue\@p@n\ed@b@x\fi}
\def\info{\cl@s@\@p@n\inf@b@x} \def\paperinfo{\cl@s@\@p@n\p@p@nf@b@x}
\def\bookinfo{\cl@s@\@p@n\b@@nf@b@x} \let\finalinfo\info
\def\P@nct{\global\setbox\p@nct@} \def\nopunct{\P@nct\box\voidb@x}
\def\p@@@t#1#2{\ifvoid\p@nct@\else#1\unhbox\p@nct@#2\fi}
\def\sp@@{\penalty-50 \space\hskip\z@ plus.1em}
\def\c@mm@{\p@@@t,\sp@@} \def\sp@c@{\p@@@t\empty\sp@@}
\def\p@tb@x#1#2{\ifvoid#1\else#2\@nb@x#1\fi}
\def\@nb@x#1{\unhbox#1\P@nct\lastbox}
\def\endr@f@{\cl@s@\nopunct
 \R@fb@x\hbox{\unhbox\r@f@b@x \p@tb@x\b@b@x\empty
 \ifvoid\j@@rb@x\ifvoid\inb@@kb@x\ifvoid\p@p@rb@x\ifvoid\b@@kb@x
  \ifvoid\p@p@nf@b@x\ifvoid\b@@nf@b@x
  \p@tb@x\v@lb@x\c@mm@ \ifvoid\y@@rb@x\else\sp@c@(\@nb@x\y@@rb@x)\fi
  \p@tb@x\is@b@x\c@mm@ \p@tb@x\p@g@b@x\c@mm@ \p@tb@x\inf@b@x\c@mm@
  \else\p@tb@x \b@@nf@b@x\c@mm@ \p@tb@x\v@lb@x\c@mm@ \p@tb@x\is@b@x\sp@c@
  \ifvoid\ed@b@x\else\sp@c@(\@nb@x\ed@b@x,\space\ifed@ ed.\else eds.\fi)\fi
  \p@tb@x\p@blb@x\c@mm@ \p@tb@x\p@bl@db@x\c@mm@ \p@tb@x\y@@rb@x\c@mm@
  \p@tb@x\p@g@b@x{\c@mm@\ifp@g@ p\p@@nt\else pp\p@@nt\fi}%
  \p@tb@x\inf@b@x\c@mm@\fi
  \else \p@tb@x\p@p@nf@b@x\c@mm@ \p@tb@x\v@lb@x\c@mm@
  \ifvoid\y@@rb@x\else\sp@c@(\@nb@x\y@@rb@x)\fi
  \p@tb@x\is@b@x\c@mm@ \p@tb@x\p@g@b@x\c@mm@ \p@tb@x\inf@b@x\c@mm@\fi
  \else \p@tb@x\b@@kb@x\c@mm@
  \p@tb@x\b@@nf@b@x\c@mm@ \p@tb@x\p@blb@x\c@mm@
  \p@tb@x\p@bl@db@x\c@mm@ \p@tb@x\y@@rb@x\c@mm@
  \ifvoid\p@g@b@x\else\c@mm@\@nb@x\p@g@b@x p\fi \p@tb@x\inf@b@x\c@mm@ \fi
  \else \c@mm@\@nb@x\p@p@rb@x\ic@\p@tb@x\p@p@nf@b@x\c@mm@
  \p@tb@x\v@lb@x\sp@c@ \ifvoid\y@@rb@x\else\sp@c@(\@nb@x\y@@rb@x)\fi
  \p@tb@x\is@b@x\c@mm@ \p@tb@x\p@g@b@x\c@mm@\p@tb@x\inf@b@x\c@mm@\fi
  \else \p@tb@x\p@p@rb@x\c@mm@\ic@\p@tb@x\p@p@nf@b@x\c@mm@
  \c@mm@\@nb@x\inb@@kb@x \p@tb@x\b@@nf@b@x\c@mm@ \p@tb@x\v@lb@x\sp@c@
  \p@tb@x\is@b@x\sp@c@
  \ifvoid\ed@b@x\else\sp@c@(\@nb@x\ed@b@x,\space\ifed@ ed.\else eds.\fi)\fi
  \p@tb@x\p@blb@x\c@mm@ \p@tb@x\p@bl@db@x\c@mm@ \p@tb@x\y@@rb@x\c@mm@
  \p@tb@x\p@g@b@x{\c@mm@\ifp@g@ p\p@@nt\else pp\p@@nt\fi}%
  \p@tb@x\inf@b@x\c@mm@\fi
  \else\p@tb@x\p@p@rb@x\c@mm@\ic@\p@tb@x\p@p@nf@b@x\c@mm@\p@tb@x\j@@rb@x\c@mm@
  \p@tb@x\v@lb@x\sp@c@ \ifvoid\y@@rb@x\else\sp@c@(\@nb@x\y@@rb@x)\fi
  \p@tb@x\is@b@x\c@mm@ \p@tb@x\p@g@b@x\c@mm@ \p@tb@x\inf@b@x\c@mm@ \fi}}
\def\m@r@f#1#2{\endr@f@\ifvoid\p@nct@\else\R@fb@x\hbox{\unhbox\r@f@b@x
 #1\unhbox\p@nct@\penalty-200\enskip#2}\fi\egroup\b@g@nr@f}
\def\endref{\endr@f@\ifvoid\p@nct@\else\R@fb@x\hbox{\unhbox\r@f@b@x.}\fi
 \parindent\r@f@nd
 \r@fit@m{\ifvoid\k@yb@x\else\if@fn@m\r@fn@{\unhbox\k@yb@x}\else
 \if@fl@b\unhbox\k@yb@x\else\r@fk@y{\unhbox\k@yb@x}\fi\fi\fi}\unhbox\r@f@b@x
 \endgraf\egroup\endgroup}
\def\moreref{\m@r@f;\empty}
\def\transl{\m@r@f;{\unskip\space
 {\sl English translation\ic@}:\penalty-66 \space}}
\def\endRefs{\endgraf\goodbreak\endgroup}

\hyphenation{acad-e-my acad-e-mies af-ter-thought anom-aly anom-alies
an-ti-deriv-a-tive an-tin-o-my an-tin-o-mies apoth-e-o-ses
apoth-e-o-sis ap-pen-dix ar-che-typ-al as-sign-a-ble as-sist-ant-ship
as-ymp-tot-ic asyn-chro-nous at-trib-uted at-trib-ut-able bank-rupt
bank-rupt-cy bi-dif-fer-en-tial blue-print busier busiest
cat-a-stroph-ic cat-a-stroph-i-cally con-gress cross-hatched data-base
de-fin-i-tive de-riv-a-tive dis-trib-ute dri-ver dri-vers eco-nom-ics
econ-o-mist elit-ist equi-vari-ant ex-quis-ite ex-tra-or-di-nary
flow-chart for-mi-da-ble forth-right friv-o-lous ge-o-des-ic
ge-o-det-ic geo-met-ric griev-ance griev-ous griev-ous-ly
hexa-dec-i-mal ho-lo-no-my ho-mo-thetic ideals idio-syn-crasy
in-fin-ite-ly in-fin-i-tes-i-mal ir-rev-o-ca-ble key-stroke
lam-en-ta-ble light-weight mal-a-prop-ism man-u-script mar-gin-al
meta-bol-ic me-tab-o-lism meta-lan-guage me-trop-o-lis
met-ro-pol-i-tan mi-nut-est mol-e-cule mono-chrome mono-pole
mo-nop-oly mono-spline mo-not-o-nous mul-ti-fac-eted mul-ti-plic-able
non-euclid-ean non-iso-mor-phic non-smooth par-a-digm par-a-bol-ic
pa-rab-o-loid pa-ram-e-trize para-mount pen-ta-gon phe-nom-e-non
post-script pre-am-ble pro-ce-dur-al pro-hib-i-tive pro-hib-i-tive-ly
pseu-do-dif-fer-en-tial pseu-do-fi-nite pseu-do-nym qua-drat-ic
quad-ra-ture qua-si-smooth qua-si-sta-tion-ary qua-si-tri-an-gu-lar
quin-tes-sence quin-tes-sen-tial re-arrange-ment rec-tan-gle
ret-ri-bu-tion retro-fit retro-fit-ted right-eous right-eous-ness
ro-bot ro-bot-ics sched-ul-ing se-mes-ter semi-def-i-nite
semi-ho-mo-thet-ic set-up se-vere-ly side-step sov-er-eign spe-cious
spher-oid spher-oid-al star-tling star-tling-ly sta-tis-tics
sto-chas-tic straight-est strange-ness strat-a-gem strong-hold
sum-ma-ble symp-to-matic syn-chro-nous topo-graph-i-cal tra-vers-a-ble
tra-ver-sal tra-ver-sals treach-ery turn-around un-at-tached
un-err-ing-ly white-space wide-spread wing-spread wretch-ed
wretch-ed-ly Brown-ian Eng-lish Euler-ian Feb-ru-ary Gauss-ian
Grothen-dieck Hamil-ton-ian Her-mit-ian Jan-u-ary Japan-ese Kor-te-weg
Le-gendre Lip-schitz Lip-schitz-ian Mar-kov-ian Noe-ther-ian
No-vem-ber Rie-mann-ian Schwarz-schild Sep-tem-ber}

\def\leftheadtext#1{} \def\rightheadtext#1{}

\let\nopagenumber\p@gen@false \let\putpagenumber\p@gen@true
\let\pagefirst\nopagenumber \let\pagenext\putpagenumber

\else

\amsppttrue

\let\twelvepoint\relax \let\Twelvepoint\relax \let\putpagenumber\relax
\let\logo@\relax \let\pagefirst\firstpage@true \let\pagenext\firstpage@false
\def\nopagenumber{\let\f@li@ld\folio\def\folio{\global\let\folio\f@li@ld}}

\def\ftext#1{\footnotetext""{\vsk-.8>\nt #1}}
\def\Ftext#1#2{\m@k@foot{\vsk-.8>\s@twd@\hangindent{#1}\nt #2}}

\def\m@R@f@[#1]{\Refs\nofrills{}\m@th\tenpoint
 {
 \s@ct{}{#1}}\wr@@c{\string\Refcd{#1}{\the\pageno}}
 \def\k@yf@##1{\hss[##1]\enspace} \let\keyformat\k@yf@
 \def\widest##1{\s@twd@\refindentwd{\tenpoint\k@yf@{##1}}}
 \let\Key\key \def\refin{\kern\refindentwd}}
\let\info\finalinfo \r@R@fs\Refs
\def\adjust#1{#1} \let\adjustend\relax
\let\adjustnext\adjust 

\fi

\outer\def\myRefs{\myR@fs} \r@st@re\proclaim
\def\bye{\par\vfill\supereject\cl@selbl\cl@secd\b@e} \r@endd@\b@e
\let\Cite\cite \let\Key\key \def\endpro{\par\endproclaim}
\let\d@c@\document \def\document{\d@c@\tenpoint}
\hyphenation{ortho-gon-al}

\newtoks\@@tp@t \@@tp@t\output
\output=\@ft@{\let\{\noexpand\the\@@tp@t}
\let\{\relax

\newif\ifVersion \Versiontrue
\def\p@n@l#1{\ifnum#1=\z@\else\penalty#1\relax\fi}

\def\s@ct#1#2{\ifVersion
 \skip@\lastskip\ifdim\skip@<1.5\bls\vskip-\skip@\p@n@l{-200}\vsk.5>%
 \p@n@l{-200}\vsk.5>\p@n@l{-200}\vsk.5>\p@n@l{-200}\vsk-1.5>\else
 \p@n@l{-200}\fi\ifdim\skip@<.9\bls\vsk.9>\else
 \ifdim\skip@<1.5\bls\vskip\skip@\fi\fi
 \vtop{\twelvepoint\raggedright\s@cf@nt\vp1\vsk->\vskip.16ex
 \s@twd@\parindent{#1}%
 \ifdim\parindent>\z@\adv\parindent.5em\fi\hang\textindent{#1}#2\strut}
 \else
 \p@sk@p{-200}{.8\bls}\vtop{\s@cf@nt\s@twd@\parindent{#1}%
 \ifdim\parindent>\z@\adv\parindent.5em\fi\hang\textindent{#1}#2\strut}\fi
 \nointerlineskip\nobreak\vtop{\strut}\nobreak\vskip-.6\bls\nobreak}

\def\s@bs@ct#1#2{\ifVersion
 \skip@\lastskip\ifdim\skip@<1.5\bls\vskip-\skip@\p@n@l{-200}\vsk.5>%
 \p@n@l{-200}\vsk.5>\p@n@l{-200}\vsk.5>\p@n@l{-200}\vsk-1.5>\else
 \p@n@l{-200}\fi\ifdim\skip@<.9\bls\vsk.9>\else
 \ifdim\skip@<1.5\bls\vskip\skip@\fi\fi
 \vtop{\elevenpoint\raggedright\s@bf@nt\vp1\vsk->\vskip.16ex%
 \s@twd@\parindent{#1}\ifdim\parindent>\z@\adv\parindent.5em\fi
 \hang\textindent{#1}#2\strut}
 \else
 \p@sk@p{-200}{.6\bls}\vtop{\s@bf@nt\s@twd@\parindent{#1}%
 \ifdim\parindent>\z@\adv\parindent.5em\fi\hang\textindent{#1}#2\strut}\fi
 \nointerlineskip\nobreak\vtop{\strut}\nobreak\vskip-.8\bls\nobreak}

\def\gadv{\global\adv} \def\gad#1{\gadv#1\@ne} \def\gadneg#1{\gadv#1-\@ne}

\newcount\t@@n \t@@n=10 \newbox\testbox

\newcount\Sno \newcount\Lno \newcount\Fno

\def\pr@cl#1{\r@st@re\pr@c@\pr@c@{#1}\global\let\pr@c@\relax}

\def\tagg#1{\tag"\rlap{\rm(#1)}\kern.01\p@"}
\def\l@L#1{\l@bel{#1}L} \def\l@F#1{\l@bel{#1}F} \def\<#1>{\l@b@l{#1}F}
\def\Tag#1{\tag{\l@F{#1}}} \def\Tagg#1{\tagg{\l@F{#1}}}
\def\Rem{\demo{\sl Remark}} \def\Ex{\demo{\bf Example}}
\def\Pf#1.{\demo{Proof #1}} \def\epf{\qed\enddemo}
\def\Ap@x{Appendix}
\def\Appendix{\Sno=64 \t@@n\@ne \wr@@c{\string\Appencd}
 \def\sf@rm{\char\the\Sno} \def\sf@rm@{\Ap@x\space\sf@rm} \def\sf@rm@@{\Ap@x}
 \def\s@ct@n##1##2{\s@ct\empty{\setbox\z@\hbox{##1}\ifdim\wd\z@=\z@
 \if##2*\sf@rm@@\else\if##2.\sf@rm@@.\else##2\fi\fi\else
 \if##2*\sf@rm@\else\if##2.\sf@rm@.\else\sf@rm@.\enspace##2\fi\fi\fi}}}
\def\Appcd#1#2#3{\def\Ap@@{\hglue-\l@ftcd\Ap@x}\ifx\@ppl@ne\empty
 \def\l@@b{\@fwd@@{#1}{\space#1}{}}\if*#2\entcd{}{\Ap@@\l@@b}{#3}\else
 \if.#2\entcd{}{\Ap@@\l@@b.}{#3}\else\entcd{}{\Ap@@\l@@b.\enspace#2}{#3}\fi\fi
 \else\def\l@@b{\@fwd@@{#1}{\c@l@b{#1}}{}}\if*#2\entcd{\l@@b}{\Ap@x}{#3}\else
 \if.#2\entcd{\l@@b}{\Ap@x.}{#3}\else\entcd{\l@@b}{#2}{#3}\fi\fi\fi}

\let\s@ct@n\s@ct
\def\s@ct@@[#1]#2{\@ft@\xdef\csname @#1@S@\endcsname{\sf@rm}\wr@@x{}%
 \wr@@x{\string\labeldef{S}\space{\?#1@S?}\space{#1}}%
 {
 \s@ct@n{\sf@rm@}{#2}}\wr@@c{\string\Entcd{\?#1@S?}{#2}{\the\pageno}}}
\def\s@ct@#1{\wr@@x{}{
 \s@ct@n{\sf@rm@}{#1}}\wr@@c{\string\Entcd{\sf@rm}{#1}{\the\pageno}}}
\def\s@ct@e[#1]#2{\@ft@\xdef\csname @#1@S@\endcsname{\sf@rm}\wr@@x{}%
 \wr@@x{\string\labeldef{S}\space{\?#1@S?}\space{#1}}%
 {
 \s@ct@n\empty{#2}}\wr@@c{\string\Entcd{}{#2}{\the\pageno}}}
\def\s@cte#1{\wr@@x{}{
 \s@ct@n\empty{#1}}\wr@@c{\string\Entcd{}{#1}{\the\pageno}}}
\def\theSno#1#2{\dff\?#1@S?{#2}%
 \wr@@x{\string\labeldef{S}\space{#2}\space{#1}}\fi}

\newif\ifd@bn@\d@bn@true
\def\Section{\gad\Sno\ifd@bn@\Fno\z@\Lno\z@\fi\@fn@xt[\s@ct@@\s@ct@}
\def\section{\gad\Sno\ifd@bn@\Fno\z@\Lno\z@\fi\@fn@xt[\s@ct@e\s@cte}
\let\Sect\Section \let\sect\section
\def\subsection{\@fn@xt*\subs@ct@\subs@ct}
\def\subs@ct#1{{\s@bs@ct\empty{#1}}\wr@@c{\string\subcd{#1}{\the\pageno}}}
\def\subs@ct@*#1{\vsk->\nobreak
 {\s@bs@ct\empty{#1}}\wr@@c{\string\subcd{#1}{\the\pageno}}}
 \def\Snodef#1{\Sno #1}

\def\l@b@l#1#2{\def\n@@{\csname #2no\endcsname}%
 \if*#1\gad\n@@ \@ft@\xdef\csname @#1@#2@\endcsname{\l@f@rm}\else\def\t@st{#1}%
 \ifx\t@st\empty\gad\n@@ \@ft@\xdef\csname @#1@#2@\endcsname{\l@f@rm}%
 \else\@ft@\ifx\csname @#1@#2@mark\endcsname\relax\gad\n@@
 \@ft@\xdef\csname @#1@#2@\endcsname{\l@f@rm}%
 \@ft@\gdef\csname @#1@#2@mark\endcsname{}%
 \wr@@x{\string\labeldef{#2}\space{\?#1@#2?}\space\ifnum\n@@<10 \space\fi{#1}}%
 \fi\fi\fi}
\def\labeldef#1#2#3{\dff\?#3@#1?{#2}}
\def\Labeldef#1#2#3{\dff\?#3@#1?{#2}\@ft@\gdef\csname @#3@#1@mark\endcsname{}}

\def\l@bel#1#2{\l@b@l{#1}{#2}\?#1@#2?}

\newcount\c@cite
\def\?#1?{\csname @#1@\endcsname}
\def\[{\@fn@xt:\c@t@sect\c@t@}
\def\c@t@#1]{{\c@cite\z@\@fwd@@{\?#1@L?}{\adv\c@cite1}{}%
 \@fwd@@{\?#1@F?}{\adv\c@cite1}{}\@fwd@@{\?#1?}{\adv\c@cite1}{}%
 \relax\ifnum\c@cite=\z@{\bf ???}\wrs@x{No label [#1]}\else
 \ifnum\c@cite=1\let\@@PS\relax\let\@@@\relax\else\let\@@PS\underbar
 \def\@@@{{\rm<}}\fi\@@PS{\?#1?\@@@\?#1@L?\@@@\?#1@F?}\fi}}
\def\(#1){{\rm(\c@t@#1])}}
\def\c@t@s@ct#1{\@fwd@@{\?#1@S?}{\?#1@S?\relax}%
 {{\bf ???}\wrs@x{No section label {#1}}}}
\def\c@t@sect:#1]{\c@t@s@ct{#1}} \let\SNo\c@t@s@ct

\newdimen\l@ftcd \newdimen\r@ghtcd \let\nlc\relax

\def\d@tt@d{\leaders\hbox to 1em{\kern.1em.\hfil}\hfill}
\def\entcd#1#2#3{\item{\l@bcdf@nt#1}{\entcdf@nt#2}\alb\kern.9em\hbox{}%
 \kern-.9em\d@tt@d\kern-.36em{\p@g@cdf@nt#3}\kern-\r@ghtcd\hbox{}\par}
\def\Entcd#1#2#3{\def\l@@b{\@fwd@@{#1}{\c@l@b{#1}}{}}\vsk.2>%
 \entcd{\l@@b}{#2}{#3}}
\def\subcd#1#2{{\adv\leftskip.333em\entcd{}{\s@bcdf@nt#1}{#2}}}
\def\Refcd#1#2{\def\t@@st{#1}\ifx\t@@st\empty\ifx\r@fl@ne\empty\relax\else
 \R@fcd{\r@fl@ne}{#2}\fi\else\R@fcd{#1}{#2}\fi}
\def\R@fcd#1#2{\sk@@p{.6\bls}\entcd{}{\hglue-\l@ftcd\R@fcdf@nt #1}{#2}}
\def\Refline{\def\r@fl@ne} \def\Refempty{\let\r@fl@ne\empty}
\def\Appencd{\par\adv\leftskip-\l@ftcd\adv\rightskip-\r@ghtcd\@ppl@ne
 \adv\leftskip\l@ftcd\adv\rightskip\r@ghtcd\let\Entcd\Appcd}
\def\appline{\def\@ppl@ne} \def\Appempty{\let\@ppl@ne\empty}
\def\Appline#1{\def\@ppl@ne{\s@bs@ct{}{#1}}}
\def\leftcd#1{\adv\leftskip-\l@ftcd\s@twd@\l@ftcd{\c@l@b{#1}\enspace}
 \adv\leftskip\l@ftcd}
\def\rightcd#1{\adv\rightskip-\r@ghtcd\s@twd@\r@ghtcd{#1\enspace}
 \adv\rightskip\r@ghtcd}
\def\C@nt{Contents} \def\Ap@s{Appendices} \def\R@fcs{References}
\def\contents{\@fn@xt*\cont@@\cont@}
\def\cont@{\@fn@xt[\cnt@{\cnt@[\C@nt]}}
\def\cont@@*{\@fn@xt[\cnt@@{\cnt@@[\C@nt]}}
\def\cnt@[#1]{\c@nt@{M}{#1}{44}{\s@bs@ct{}{\@ppl@f@nt\Ap@s}}}
\def\cnt@@[#1]{\c@nt@{M}{#1}{44}{}}
\def\endco{\par\penalty-500\vsk>\vskip\z@\endgroup}
\def\readcd{\@np@t{\jobname.cd}}
\def\Cde{\@fn@xt*\Cde@@\Cde@}
\def\Cde@{\@fn@xt[\Cd@{\Cd@[\C@nt]}}
\def\Cde@@*{\@fn@xt[\Cd@@{\Cd@@[\C@nt]}}
\def\Cd@[#1]{\cnt@[#1]\readcd\endco}
\def\Cd@@[#1]{\cnt@@[#1]\readcd\endco}
\def\contlabeldef{\def\c@l@b}

\long\def\c@nt@#1#2#3#4{\s@twd@\l@ftcd{\c@l@b{#1}\enspace}
 \s@twd@\r@ghtcd{#3\enspace}\adv\r@ghtcd1.333em
 \def\@ppl@ne{#4}\def\r@fl@ne{\R@fcs}\s@ct{}{#2}\B@gr@\parindent\z@\let\nlc\nl
 \let\nl\relax\parskip.2\bls\adv\leftskip\l@ftcd\adv\rightskip\r@ghtcd}

\def\writecd{\immediate\openout\@@cd\jobname.cd \def\wr@@c{\write\@@cd}
 \def\cl@secd{\immediate\write\@@cd{\string\endinput}\immediate\closeout\@@cd}
 \def\closecd{\cl@secd\global\let\cl@secd\relax}}
\let\cl@secd\relax \def\wr@@c#1{} \let\closecd\relax

\def\dff{\@ft@\d@f} \def\d@f{\@ft@\def}
\def\edff{\@ft@\ed@f} \def\ed@f{\@ft@\edef}
\def\defi#1#2{\def#1{#2}\wr@@x{\string\def\string#1{#2}}}

\def\qed{\hbox{}\nobreak\hfill\nobreak{\m@th$\,\square$}}
\def\back#1 {\strut\kern-.33em #1\enspace\ignore} 

\def\hcor#1{\advance\hoffset by #1}
\def\vcor#1{\advance\voffset by #1}
\let\bls\baselineskip \let\ignore\ignorespaces
\ifx\ic@\undefined \let\ic@\/\fi
\def\vsk#1>{\vskip#1\bls} \let\adv\advance
\def\vv#1>{\vadjust{\vsk#1>}\ignore}
\def\vvn#1>{\vadjust{\nobreak\vsk#1>\nobreak}\ignore}
\def\vvv#1>{\vskip\z@\vsk#1>\nt\ignore}
\def\vvgood{\vadjust{\penalty-500}}
\def\nngood{\noalign{\penalty-500}}

\def\Goodbreak{\par\penalty-\@m}
\def\Par{\vsk.5>} \def\setparindent{\edef\Parindent{\the\parindent}}
\def\Type{\vsk.5>\bgroup\parindent\z@\tt\rightskip\z@ plus1em minus1em%
 \spaceskip.3333em \xspaceskip.5em\relax}
\def\endType{\vsk.5>\egroup\nt} 

\let\Hat\widehat \let\Tilde\widetilde \let\dollar\$ \let\ampersand\&
\let\sss\scriptscriptstyle  
\let\vp\vphantom \let\hp\hphantom \let\nt\noindent
\let\cline\centerline \let\lline\leftline \let\rline\rightline
\def\nn#1>{\noalign{\vskip#1\p@@}} \def\NN#1>{\openup#1\p@@}
\def\cnn#1>{\noalign{\vsk#1>}}
\def\Cup{\bigcup\limits} \def\Cap{\bigcap\limits}
\let\Lim\lim \def\lim{\Lim\limits} \let\Sum\sum \def\sum{\Sum\limits}
\def\Plus{\bigoplus\limits} 
\let\Prod\prod \def\prod{\Prod\limits} \let\Int\int \def\int{\Int\limits}

\def\tsum{\mathop{\tsize\Sum}\limits} 
\def\tprod{\mathop{\tsize\Prod}\limits} \def\&{.\kern.1em}
\def\nl{\leavevmode\hfill\break} \def\~{\leavevmode\@fn@xt~\m@n@s\@md@@sh}
\def\@md@@sh{\@fn@xt-\d@@sh\@md@sh} \def\@md@sh{\raise.13ex\hbox{--}}
\def\m@n@s~{\raise.15ex\mbox{-}} \def\d@@sh-{\raise.15ex\hbox{-}}

\let\procent\% \def\%#1{\ifmmode\mathop{#1}\limits\else\procent#1\fi}
\let\@ml@t\" \def\"#1{\ifmmode ^{(#1)}\else\@ml@t#1\fi}
\let\@c@t@\' \def\'#1{\ifmmode _{(#1)}\else\@c@t@#1\fi}
\let\colon\: \def\:{^{\vp{\topsmash|}}} \def\vpa{{\vp|}^{\]*}}

\let\texspace\ \def\ {\ifmmode\alb\fi\texspace} \def\.{\d@t\ignore}

\def\Oldskips{\def\>{{\!\;}} \def\]{{\!\!\;}} \def\){\>\]} \def\}{\]\]}
 \def\d@t{.\ }}
\def\Newskips{\def\d@t{.\alb\hskip.3em}
\def\>{\RIfM@\mskip.666667\thinmuskip\relax\else\kern.111111em\fi}
\def\}{\RIfM@\mskip-.666667\thinmuskip\relax\else\kern-.111111em\fi}
\def\){\RIfM@\mskip.333333\thinmuskip\relax\else\kern.0555556em\fi}
\def\]{\RIfM@\mskip-.333333\thinmuskip\relax\else\kern-.0555556em\fi}}
\Newskips

\let\n@wp@ge\newpage \def\newpage{\endgraf\n@wp@ge}
\let\=\m@th \def\mbox#1{\hbox{\m@th$#1$}}
\def\mtext#1{\text{\m@th$#1$}} \def\^#1{\text{\m@th#1}}
\def\Line#1{\kern-.5\hsize\line{\m@th$\dsize#1$}\kern-.5\hsize}
\def\Lline#1{\kern-.5\hsize\lline{\m@th$\dsize#1$}\kern-.5\hsize}
\def\Cline#1{\kern-.5\hsize\cline{\m@th$\dsize#1$}\kern-.5\hsize}
\def\Rline#1{\kern-.5\hsize\rline{\m@th$\dsize#1$}\kern-.5\hsize}

\def\Ll@p#1{\llap{\m@th$#1$}} \def\Rl@p#1{\rlap{\m@th$#1$}}
 \def\Cl@p#1{\llap{\m@th$#1$\hss}}
\def\Llap#1{\mathchoice{\Ll@p{\dsize#1}}{\Ll@p{\tsize#1}}{\Ll@p{\ssize#1}}%
 {\Ll@p{\sss#1}}}
\def\Clap#1{\mathchoice{\Cl@p{\dsize#1}}{\Cl@p{\tsize#1}}{\Cl@p{\ssize#1}}%
 {\Cl@p{\sss#1}}}
\def\Rlap#1{\mathchoice{\Rl@p{\dsize#1}}{\Rl@p{\tsize#1}}{\Rl@p{\ssize#1}}%
 {\Rl@p{\sss#1}}}
 
\def\LRtph#1#2{\setbox\z@\hbox{#1}\dimen\z@\wd\z@\hbox{\hbox to\dimen\z@{#2}}}
\def\LRph#1#2{\LRtph{\m@th$#1$}{\m@th$#2$}}

\def\Lto#1{\setbox\z@\mbox{\tsize{#1}}%
 \mathrel{\mathop{\hbox to\wd\z@{\rightarrowfill}}\limits#1}}
\def\Lgets#1{\setbox\z@\mbox{\tsize{#1}}%
 \mathrel{\mathop{\hbox to\wd\z@{\leftarrowfill}}\limits#1}}
\def\vpb#1{{\vp{\big(}}^{\]#1}}

\let\alb\allowbreak 
\def\ald{\noalign{\alb}} \let\alds\allowdisplaybreaks

\let\o\circ \let\x\times \let\ox\otimes 
\let\sub\subset  \let\tabs\+
\let\le\leqslant \let\ge\geqslant
\let\der\partial \let\8\infty \let\*\star
\let\bra\langle \let\ket\rangle
 
\let\map\mapsto  
 
\let\Vert\parallel \def\vert{\ |\ } \def\nin{\not\in}

\let\lb\lbrace \let\rb\rbrace

   \def\Rb{\bigr\rb}

\def\lsym#1{#1\alb\ldots\relax#1\alb}
\def\lc{\lsym,}   \def\lox{\lsym\ox}

 \def\Im{\mathop{\roman{Im}\>}}
\def\End{\mathop{\roman{End}\>}\nolimits}
\def\Hom{\mathop{\roman{Hom}\>}\nolimits}

\def\Ker{\mathop{\roman{Ker}\>}}

\def\sign{\mathop{\roman{sign}\)}\limits}

\def\id{\roman{id}}  
\def\1{^{-1}} \let\underscore\_ \def\_#1{_{\Rlap{#1}}}
\def\vst#1{{\lower1.9\p@@\mbox{\bigr|_{\raise.5\p@@\mbox{\ssize#1}}}}}
\def\vrp#1:#2>{{\vrule height#1 depth#2 width\z@}}
\def\vru#1>{\vrp#1:\z@>} \def\vrd#1>{\vrp\z@:#1>}
\def\qqq{\qquad\quad} 
\def\sscr#1{\raise.3ex\mbox{\sss#1}} \def\@@PS{\bold{OOPS!!!}}

\def\intcl{\mathop
 {\Rlap{\raise.3ex\mbox{\kern.12em\curvearrowleft}}\int}\limits}
\def\intcr{\mathop
 {\Rlap{\raise.3ex\mbox{\kern.24em\curvearrowright}}\int}\limits}

\def\pms{\raise.25ex\mbox{\ssize\pm}\>}
\def\mps{\raise.25ex\mbox{\ssize\mp}\>}
\def\pss{{\sscr+}} \def\mss{{\sscr-}}
\def\pmss{{\sscr\pm}} 

\let\al\alpha
\let\bt\beta
\let\gm\gamma  
\let\dl\delta \let\Dl\Delta 
 \let\eps\varepsilon \let\epsilon\eps

\let\zt\zeta
\let\tht\theta \let\Tht\Theta

\let\la\lambda \let\La\Lambda

 \let\phi\varphi

\let\om\omega \let\Om\Omega 

\def\C{\Bbb C}

\def\Z{\Bbb Z}

\def\II{\Bbb I}

\def\Zp{\Z_{\ge 0}} \def\Zpp{\Z_{>0}}

\def\difl/{differential} \def\dif/{difference}
\def\cf.{cf.\ \ignore} \def\Cf.{Cf.\ \ignore}
\def\egv/{eigenvector} \def\eva/{eigenvalue} \def\eq/{equation}
\def\lhs/{the left hand side} \def\rhs/{the right hand side}
\def\Lhs/{The left hand side} \def\Rhs/{The right hand side}
\def\gby/{generated by} \def\wrt/{with respect to} \def\st/{such that}
\def\resp/{respectively} \def\off/{offdiagonal} \def\wt/{weight}
\def\pol/{polynomial} \def\rat/{rational} \def\tri/{trigonometric}
\def\fn/{function} \def\var/{variable} \def\raf/{\rat/ \fn/}
\def\inv/{invariant} \def\hol/{holomorphic} \def\hof/{\hol/ \fn/}
\def\mer/{meromorphic} \def\mef/{\mer/ \fn/} \def\mult/{multiplicity}
\def\sym/{symmetric} \def\perm/{permutation} \def\fd/{finite-dimensional}
\def\rep/{representation} \def\irr/{irreducible} \def\irrep/{\irr/ \rep/}
\def\hom/{homomorphism} \def\aut/{automorphism} \def\iso/{isomorphism}
\def\lex/{lexicographical} \def\as/{asymptotic} \def\asex/{\as/ expansion}
\def\ndeg/{nondegenerate} \def\neib/{neighbourhood} \def\deq/{\dif/ \eq/}
\def\hw/{highest \wt/} \def\gv/{generating vector} \def\eqv/{equivalent}
\def\msd/{method of steepest descend} \def\pd/{pairwise distinct}
\def\wlg/{without loss of generality} \def\Wlg/{Without loss of generality}
\def\onedim/{one-dimensional} \def\qcl/{quasiclassical} \def\hwv/{\hw/ vector}
\def\hgeom/{hyper\-geometric} \def\hint/{\hgeom/ integral}
\def\hwm/{\hw/ module} \def\emod/{evaluation module} \def\Vmod/{Verma module}
\def\symg/{\sym/ group} \def\sol/{solution} \def\eval/{evaluation}
\def\anf/{analytic \fn/} \def\anco/{analytic continuation}
\def\qg/{quantum group} \def\qaff/{quantum affine algebra}

\def\Rm/{\^{$R$-}matrix} \def\Rms/{\^{$R$-}matrices} \def\YB/{Yang-Baxter \eq/}
\def\Ba/{Bethe ansatz} \def\Bv/{Bethe vector} \def\Bae/{\Ba/ \eq/}
\def\KZv/{Knizh\-nik-Zamo\-lod\-chi\-kov} \def\KZvB/{\KZv/-Bernard}
\def\KZ/{{\sl KZ\/}} \def\qKZ/{{\sl qKZ\/}}
\def\KZB/{{\sl KZB\/}} \def\qKZB/{{\sl qKZB\/}}
\def\qKZo/{\qKZ/ operator} \def\qKZc/{\qKZ/ connection}
\def\KZe/{\KZ/ \eq/} \def\qKZe/{\qKZ/ \eq/} \def\qKZBe/{\qKZB/ \eq/}

\def\h@ph{\discretionary{}{}{-}} \def\$#1$-{\,\^{$#1$}\h@ph}

\def\TFT/{Research Insitute for Theoretical Physics}
\def\HY/{University of Helsinki} \def\AoF/{the Academy of Finland}
\def\CNRS/{Supported in part by MAE\~MICECO\~CNRS Fellowship}
\def\LPT/{Laboratoire de Physique Th\'eorique ENSLAPP}
\def\ENSLyon/{\'Ecole Normale Sup\'erieure de Lyon}
\def\LPTaddr/{46, All\'ee d'Italie, 69364 Lyon Cedex 07, France}
\def\enslapp/{URA 14\~36 du CNRS, associ\'ee \`a l'E.N.S.\ de Lyon,
au LAPP d'Annecy et \`a l'Universit\`e de Savoie}
\def\ensemail/{vtarasov\@ enslapp.ens-lyon.fr}
\def\DMS/{Department of Mathematics, Faculty of Science}
\def\DMO/{\DMS/, Osaka University}
\def\DMOaddr/{Toyonaka, Osaka 560, Japan}
\def\dmoemail/{vt\@ math.sci.osaka-u.ac.jp}
\def\MPI/{Max\)-Planck\)-Institut} \def\MPIM/{\MPI/ f\"ur Mathematik}
\def\MPIMaddr/{P\]\&O.\ Box 7280, D\~-\]53072 \,Bonn, Germany}
\def\mpimemail/{tarasov\@ mpim-bonn.mpg.de}
\def\SPb/{St\&Peters\-burg}
\def\home/{\SPb/ Branch of Steklov Mathematical Institute}
\def\homeaddr/{Fontanka 27, \SPb/ \,191011, Russia}
\def\homemail/{vt\@ pdmi.ras.ru}
\def\absence/{On leave of absence from \home/}
\def\support/{Supported in part by}
\def\UNC/{Department of Mathematics, University of North Carolina}
\def\ChH/{Chapel Hill}
\def\UNCaddr/{\ChH/, NC 27599, USA} \def\avemail/{av\@ math.unc.edu}
\def\grant/{NSF grant DMS\~9501290}	
\def\Grant/{\support/ \grant/}

\def\Aomoto/{K\&Aomoto}
\def\Cher/{I\&Cherednik}
\def\Dri/{V\]\&G\&Drin\-feld}
\def\Etingof/{P\]\&Etin\-gof}
\def\Fadd/{L\&D\&Fad\-deev}
\def\Feld/{G\&Felder}
\def\Fre/{I\&B\&Fren\-kel}
\def\Gustaf/{R\&A\&Gustafson}
\def\Kazh/{D\&Kazhdan} \def\Kir/{A\&N\&Kiril\-lov}
\def\Kor/{V\]\&E\&Kore\-pin}
\def\Lusz/{G\&Lusztig}
\def\MN/{M\&Naza\-rov}
\def\Resh/{N\&Reshe\-ti\-khin} \def\Reshy/{N\&\]Yu\&Reshe\-ti\-khin}
\def\SchV/{V\]\&\]V\]\&Schecht\-man} \def\Sch/{V\]\&Schecht\-man}
\def\Skl/{E\&K\&Sklya\-nin}
\def\Smirn/{F\]\&Smirnov} \def\Smirnov/{F\]\&A\&Smirnov}
\def\Takh/{L\&A\&Takh\-tajan}
\def\VT/{V\]\&Ta\-ra\-sov} \def\VoT/{V\]\&O\&Ta\-ra\-sov}
\def\Varch/{A\&\]Var\-chenko} \def\Varn/{A\&N\&\]Var\-chenko}

\def\Astq/{Ast\'erisque}
\def\AMS/{Amer\. Math\. Society}
\def\CMP/{Comm\. Math\. Phys.{}}
\def\DMJ/{Duke\. Math\. J.{}}
\def\FAA/{Func\. Anal\. Appl.{}}
\def\IJM/{Int\. J\. Math.{}}
\def\IMRN/{Int\. Math\. Res\. Notices}
\def\Inv/{Invent\. Math.{}} 
\def\JPA/{J\. Phys\. A{}}
\def\JSM/{J\. Soviet\ Math.{}}
\def\JSP/{J\. Stat\. Phys.{}}
\def\LMP/{Lett\. Math\. Phys.{}}
\def\LMJ/{Leningrad Math\. J.{}}
\def\LpMJ/{\SPb/ Math\. J.{}}
\def\SIAM/{SIAM J\. Math\. Anal.{}}
\def\SMNS/{Selecta Math., New Series}
\def\SPbMJ/{\SPb/ Math\. J.{}}
\def\TMP/{Theor\. Math\. Phys.{}}
\def\ZNS/{Zap\. nauch\. semin. LOMI}

\def\ASMP/{Advanced Series in Math\. Phys.{}}

\def\AMSa/{AMS \publaddr Providence}
\def\Birk/{Birkh\"auser}
\def\CUP/{Cambridge University Press} \def\CUPa/{\CUP/ \publaddr Cambridge}
\def\Spri/{Springer-Verlag} \def\Spria/{\Spri/ \publaddr Berlin}
\def\WS/{World Scientific} \def\WSa/{\WS/ \publaddr Singapore}

\newbox\lefthbox \newbox\righthbox

\let\sectsep. \let\labelsep. \let\contsep. \let\labelspace\relax
\let\sectpre\relax \let\contpre\relax
\def\sf@rm{\the\Sno} \def\sf@rm@{\sectpre\sf@rm\sectsep}
\def\c@l@b#1{\contpre#1\contsep}
\def\l@f@rm{\ifd@bn@\sf@rm\labelsep\fi\labelspace\the\n@@}

\def\sectformdef{\def\sf@rm}

\let\DoubleNum\d@bn@true \let\SingleNum\d@bn@false

\def\NoNewNum{\let\writeldf\relax\def\l@b@l##1##2{\if*##1%
 \@ft@\xdef\csname @##1@##2@\endcsname{\mbox{*{*}*}}\fi}}
\def\NoNewTime{\def\todaydef##1{\def\today{##1}}
 \def\nowtimedef##1{\def\nowtime{##1}}}
\def\NoInput{\let\Input\input\let\writeldf\relax}
\def\Fixed{\NoNewTime\NoInput}

\def\sectfont#1{\def\s@cf@nt{#1}} \sectfont\bf
\def\subsectfont#1{\def\s@bf@nt{#1}} \subsectfont\it
\def\Entcdfont#1{\def\entcdf@nt{#1}} \Entcdfont\relax
\def\labelcdfont#1{\def\l@bcdf@nt{#1}} \labelcdfont\relax
\def\pagecdfont#1{\def\p@g@cdf@nt{#1}} \pagecdfont\relax
\def\subcdfont#1{\def\s@bcdf@nt{#1}} \subcdfont\it
\def\applefont#1{\def\@ppl@f@nt{#1}} \applefont\bf
\def\Refcdfont#1{\def\R@fcdf@nt{#1}} \Refcdfont\bf

\ifamsppt\else\issuefont{no.\kern.25em}\fi

\tenpoint

\Fixed

\loadbold

\Magset
\ifUS
\PaperUS
\else
\PaperA4
\fi

\font@\Beufm=eufm10 scaled 1440
\font@\Eufm=eufm8 scaled 1440
\newfam\Frakfam \textfont\Frakfam\Eufm
\def\Frak{\fam\Frakfam}

\newif\ifMPIM

\let\Sh S

\def\xspace{\kern.34em}

\def\Th#1{\pr@cl{Theorem\xspace\l@L{#1}}\ignore}
\def\Lm#1{\pr@cl{Lemma\xspace\l@L{#1}}\ignore}
\def\Cr#1{\pr@cl{Corollary\xspace\l@L{#1}}\ignore}
\def\Df#1{\pr@cl{Definition\xspace\l@L{#1}}\ignore}
\def\Cj#1{\pr@cl{Conjecture\xspace\l@L{#1}}\ignore}
\def\Prop#1{\pr@cl{Proposition\xspace\l@L{#1}}\ignore}
\def\Pf#1.{\demo{Proof #1}\bgroup\ignore}
\def\epf{\qed\par\egroup\enddemo}

\def\Vert{\ \,\big|\ \,}
\def\cddt{\,{\cdot}\>}

\let\pim\varpi

\def\PP{\Bbb P}
\def\QQ{\Bbb Q}

\def\ab{\Rlap{\bar{\phantom a\)}\]}a}
\def\aab{\Rlap{\>\overline{\}\phantom{\aa}\}}\>}\aa}
\def\ccb{\Rlap{\>\overline{\]\phantom{\cc}\}}\)}\cc}

\def\etb{\bar\eta}
\def\fb{\bar f}
\def\Fb{\Rlap{\,\)\overline{\!\]\phantom F\]}\)}F}
\def\gb{\Rlap{\)\bar{\]\phantom g\)}\]}g}

\def\Rb{\Rlap{\,\)\overline{\!\]\phantom R\]}\)}R}
\def\tb{\Rlap{\)\bar{\]\phantom t\)}\]}t}
\def\Tcb{\Rlap{\,\overline{\!\phantom\Tc\]}\)}\Tc}

\def\aa{\boldkey a}
\def\bb{\boldkey b}
\def\cc{\boldkey c}
\def\ib{\boldkey i}
\def\jb{\boldkey j}
\def\Sb{\bold S}

\def\Bc{\Cal B}
\def\Cc{\Cal C}
\def\Dc{\Cal D}
\def\Ec{\Cal E}
\def\Fc{\Cal F}
\def\Gc{\Cal G}
\def\Jc{\Cal J}

\def\Mc{\Cal M}
\def\Nc{\Cal N}
\def\Rc{\Cal R}
\def\Tc{\Cal T}

\def\phc{\check\phi}

\def\bg{\frak b}
\def\dg{\frak d}
\def\eg{{\Frak e}}
\def\g{\frak g}
\def\gl{\frak{gl}}

\def\hg{\frak h}

\def\ng{\frak n}

\def\ellh{\hat\ell}

\let\qh q	
\def\Qh{\Hat Q}
\def\tth{\Rlap{\)\hat{\]\phantom t\)}\]}t}

\def\I{\roman I}
\def\II{\roman{II}}
\def\III{\roman{III}}
\def\IV{\roman{IV}}

\def\one{\slanted 1}
\def\onb{\Rlap{\)\bar{\]\phantom{\one}}}\one}

\def\nut{\Rlap{\)\tilde{\]\phantom\nu\)}\]}\nu}

\def\elti{\tilde\ell}
\def\Ect{\Tilde\Ec}

\def\Jti{\Rlap{\,\,\)\Tilde{\!\!\]\phantom J}}J}

\def\mti{\Tilde m}
\def\Qti{\Rlap{\>\Tilde{\}\phantom Q\>}\}}Q}
\def\Rti{\Rlap{\>\Tilde{\}\phantom R\>}\}}R}

\def\Vti{\Tilde V}
\def\vti{\tilde v}

\def\phib{\phi_{\sss\bullet}}
\def\phiv{{\sss\]}\vec{{\sss\)}\phi}}

\def\dag{{\sscr\dagger}}
\def\dagg{{\sss\dagger}}

\def\Tps{T^{\smash{\)\pss}}} \def\Tms{T^{\smash{\)\mss}}}
\def\Tpms{T^{\smash{\)\pmss}}}

\def\hga{\hg^*} 

\def\pii{\pi i}  \def\ehg{e\)(\hg)}
 \def\ZtZ{\Z+\}\tau\)\Z} \def\QtQ{\QQ+\]\tau\)\QQ}
 
\def\bgp{\bg_{\sss+}} \def\bgm{\bg_{\sss-}} \def\bgpm{\bg_{\sss\pm}}
\def\ngp{\ng_{\sss+}} \def\ngm{\ng_{\sss-}} \def\ngpm{\ng_{\sss\pm}}

\def\dimC{\dim{\vp q}_{\]\C\>}}
\def\dimF{\dim{\vp q}_{\!\italic{Fun}(\C{\sss\>})\>}}
\def\dimR{\dim{\vp q}_{\!\italic{Rat}(\C{\sss\>})\>}}
\def\oxC{\mathrel{\ox{\vp q}_{\)\C\}}}}
\def\+#1{^{\lb#1\rb}}

\def\Ad{\mathop{\roman{Ad}}\nolimits}
\def\Det{\mathop{\roman{Det}}\nolimits}
\def\Diff{\mathop{\slanted{D\/}}\nolimits}

\def\Fun{\mathop{\slanted{Fun\/}\)}\nolimits}
\def\Funb{\mathop{\slanted{Fun\/\)}_{\sss\bullet}}\nolimits}
\def\Rat{\mathop{\slanted{Rat\/}\)}\nolimits}
\def\Mor{\mathop{\roman{Mor}\>}}

\def\wsl{\mathop{\slanted{wt\/}}}

\def\CN{\C^{\)N}} \def\NC{\Nc_\C\:} \def\DiffC{\Diff(\C\))}
\def\FunC{\Fun(\C\))} \def\Funt{\Fun^{\ox\)2\]}(\C\))}
\def\FunU{\Fun(U)} \def\FunV{\Fun(V)} \def\FunVa{\Fun(V^*)} \def\FunW{\Fun(W)}
\def\FunM{\Fun(M_\mu)} \def\Fend#1{\Fun\bigl(\End(#1)\bigr)}
\def\FVmu{\Fun\bigl(\Vmu\)\bigr)} \def\FVamu{\Fun\bigl(\Vamu\)\bigr)}
\def\RVmu{\Rat\bigl(\Vmu\)\bigr)}

\def\RatC{\Rat(\C\))} \def\Ratt{\Rat^{\ox\)2\]}(\C\))}
\def\RatV{\Rat(V)} \def\Rend#1{\Rat\bigl(\End(#1)\bigr)}

\def\elhV{\Rlap{\phantom{\ell}^{\>\sscr V}}\ellh}

\def\vmk#1{v_{\mu,\Qh\)}[\)#1\)]}
\def\vmQ{v_{\mu,\Qh}} \def\MmQ{M_{\mu,\Qh}} \def\NmQ{N_{\}\mu,\Qh}}
\def\SmQ{S_{\]\mu,\Qh}} \def\VmQ{V_{\]\mu,\Qh}} \def\chimQ{\chi_{\mu,\Qh}}
\def\BmQ{B_{\]\mu,\Qh}} \def\CmQ{C_{\]\mu,\Qh}} \def\pimQ{\pi_{\mu,\Qh}}
\def\MmQa{\Rlap{{\phantom M}^*}\MmQ} \def\vmQa{\Rlap{{\phantom v}^*}\vmQ}
\def\MmQt{M_{\mu,\Qti}} \def\vmQt{v_{\mu,\Qti}} \def\SmQt{S_{\]\mu,\Qti}}
 \def\CmQt{C_{\]\mu,\Qti}}
\def\VmQt{V_{\]\mu,\Qti}} 
\def\MmQta{\Rlap{{\phantom M}^*}\MmQt} \def\vmQta{\Rlap{{\phantom v}^*}\vmQt}

\def\SmQV{\Rlap{\phantom{S}^{{\sss\)}V}}\SmQ}
\def\Mrm{\Rlap{M_\mu}M^{\sscr{\)\italic{rat}}}} \def\RatMm{\Rat(\Mrm)}
\def\Srm{\Rlap{S_\mu}S^{\sscr{\)\italic{rat}}}} \def\RatSm{\Rat(\Sslm)}
\def\Zrm{\Rlap{Z_\mu}Z^{\sscr{\)\italic{rat}}}}  \def\RatZm{\Rat(\Zslm)}
\def\Mslm{\Rlap{M_\mu}M^{\sscr{\ssize\frak{\)s\]l}}}} \def\Msmn{\Mslm[\)\nu\)]}
\def\Sslm{\Rlap{S_\mu}S^{\sscr{\tsize\)\ssize\frak{s\]l}}}}
\def\Zslm{\Rlap{Z_\mu}Z^{\sscr{\ssize\frak{\)s\]l}}}}
\def\Vslm{\Rlap{V_\mu}V^{\sscr{\ssize\frak{\)s\]l}}}} \def\dmn{d_\mu[\)\nu\)]}
\def\Mmn{M_\mu[\)\nu\)]} \def\Umn{U_\mu[\)\nu\)]}
\def\Vmn{V_\mu[\)\nu\)]} \def\Vsmn{\Vslm[\)\nu\)]}
\def\KSmn{(\Ker S_\mu)[\)\nu\)]} \def\Zmn{Z_\mu[\)\nu\)]}
\def\Mcb{\Mc\)[\)\bt\)]} \def\MCb{\Mc_{\)\C\)}\:[\)\bt\)]}
\def\Mcmb{\Mc_\mu[\)\bt\)]}
\def\Agm{A^{\gm\]}_{\vp1}} \def\Bgm{B^{\gm\]}_{\vp1}}
\def\Ar{A^{\sscr{\italic{rat}\}}}} \def\Br{B^{\sscr{\italic{rat}\}}}}

\def\suan{\sum_{a=1}^N} \def\tsuan{\tsum_{a=1}^N} \def\pran{\prod_{a=1}^N}
\def\sumib{\sum_{\ib\)\in\)\Sb_N\!\]}} \def\sumibk{\sum_{\ib\)\in\)\Sb_k\!\]}}
\def\sumcc{\sum_{\cc\)\in Y_k\!\]}}

\def\LN{L^{\wedge N}} \def\RN{R^{\wedge N}}
\def\TN{T^{\wedge N}} \def\tN{t^{\)\wedge N}} 
\def\Tw#1_#2{T\_{#2}\vpb{\,\)\wedge #1\}}}
\def\Tws#1{\Rlap{\phantom{T}\vpb{\>\wedge #1}}T}

\def\vox{v_1\]\lox v_N}
\def\voxi{v_{i_1}\!\lox v_{i_N}} \def\voxj{v_{j_1}\!\lox v_{j_N}}

\def\ano{a=1\lc N} \def\bno{b=1\lc N} 
  \def\fN{f_1\lc f_N}
\def\QN{Q_1\lc Q_N} 

 \def\laN{\la_1\lc\la_N}

\def\abn{1\le a<b\le N}

\def\taak{t_{a_1a_1}\}\ldots\)t_{a_ka_k}}
\def\tabk{t_{a_1b_1}\}\ldots\)t_{a_kb_k}}
\def\tabii{t_{a_ib_i}\)t_{a_{i+1}b_{i+1}}\)t_{a_{i+2}b_{i+2}}}
\def\tabt{t_{a_1b_1}\)t_{a_2b_2}\)t_{a_3b_3}}
\def\thabk{\tth_{a_1b_1}\}\ldots\)\tth_{a_kb_k}}
\def\thcdl{\tth_{c_1d_1}\}\ldots\)\tth_{c_ld_l}}
\def\eabk{e_{a_1b_1}\}\ldots\)e_{a_kb_k}}
\def\ecdl{e_{c_1d_1}\}\ldots\)e_{c_ld_l}}

\def\egz{\eg\)[\)0\)]}
 \def\Umu{U[\)\mu\)]}
\def\Ulmu{U_\la[\)\mu\)]} \def\Uimu{U_\8[\)\mu\)]}
\def\Vmu{V[\)\mu\)]} \def\Vnu{V[\)\nu\)]}
\def\Vamu{V^*[\)\mu\)]} \def\Vmua{\bigl(\Vmu\)\bigr)\vpa}

\def\sltwo{\frak{sl}_2}
\def\gln{\gl_N} \def\sln{\frak{sl}_N}
 \def\Usl{U(\sln)}
\def\Ygl{Y\](\gln)} \def\Ysl{Y\](\sln)}

\def\Eqg{E_{\tau\],\)\gm}(\sln)} \def\Eqgt{E_{\tau\],\)\gm}(\sltwo)}
\def\esl{e_{\tau\],\)\gm}(\sln)} 
\def\eslt{e_{\tau\],\)\gm}(\sltwo)}

\def\els{e^{\sss\Cal O}_{\tau\],\)\gm}(\sln)}
\def\elsp{e^{\sss\Cal O}_{\tau+1,\)\gm}(\sln)}
\def\elsm{e^{\sss\Cal O}_{-1/\]\tau\],\)-\)\gm/\]\tau}(\sln)}
\def\elst{e^{\sss\Cal O}_{\tau\],\)\gm}(\sltwo)}

\def\eslr{e_{\italic{rat}}(\sln)}
\def\elsr{e^{\sss\Cal O}_{\italic{rat}}(\sln)}

 \def\Fsl{F\bigl(SL(N)\bigr)}

\def\Ro{\mathchoice{{\%{\vru1.08ex>\smash{R}}^{\>\smash{\sss\o}\}}}}
 {{\%{\vru1.08ex>\smash{R}}^{\>\smash{\sss\o}\}}}}
 {{\%{\vru.75ex>\smash{R}}^{\)\smash{\sss\o}\]}}}{\@@PS}}

\def\wtd/{\wt/ decomposition} \def\DYB/{dynamical \YB/}
\def\qlo/{quantum loop algebra} \def\eqg/{elliptic \qg/}
\def\oalg/{operator algebra} \def\dqg/{dynamical \qg/}
\def\phf/{phase \fn/} \def\wtf/{\wt/ \fn/} \def\pert/{perturbation}
\def\diag/{diagonalizable} \def\dwt/{dynamical \wt/} \def\dhw/{dynamical \hw/}

\def\Vval/{\$V\}$-\)valued} 

\def\endo/{endomorphism}
\def\gmod/{\$\gln$-module} \def\smod/{\$\sln$-module}
\def\Ymod/{\$\Ygl$-module}
\def\hmod/{\]\$\hg\)$-module} \def\dhmod/{\diag/ \hmod/}
\def\Emod/{\$\Eqg$-module}
\def\emod/{\$\esl$-module} \def\ermod/{\$\eslr$-module}

\def\asmod/{admissible \$\sln$-module}
\def\aemod/{admissible \$\esl$-module}
\def\amod/{admissible module} \def\ahmod/{admissible \hmod/}

\def\eVmod/{\eval/ \Vmod/} \def\epoint/{\eval/ point} \def\ehom/{\eval/ \hom/}
\def\evmod/{\eval/ module} \def\evrep/{\eval/ \rep/}
\def\mult/{multiplicative} \def\mform/{\mult/ form}
\def\mccl/{\mult/ cocycle} \def\mcclo/{\mult/ \$1$-cocycle}
\def\mkform/{\mult/ \$k\)$-form}
\def\mcbd/{\mult/ coboundary}
\def\singv/{singular vector} \def\regsv/{regular \singv/}
\def\nord/{normally ordered} \def\norm/{\nord/ monomial}
\def\obr/{ordered by rows} \def\obrm/{\obr/ monomial}
\def\orm/{ordered monomial} \def\orule/{ordering rule}
\def\trans/{transformation} \def\regtr/{regular \trans/}
\def\Sform/{Shapovalov form} \def\Spair/{Shapovalov pairing}
\def\dSform/{dynamical \Sform/} \def\dSpair/{dynamical \Spair/}
\def\cmod/{contragradient module} \def\cVmod/{contragradient \Vmod/}
\def\cform/{contravariant form} \def\fdim/{finite-dimen\-sional}
\def\HC/{Harish-Chandra} \def\dint/{dominant integral \wt/}

\let\goodbm\relax  \let\mmgood\relax 
   
\let\mline\relax \def\vvm#1>{\ignore} \def\vvnm#1>{\ignore} \def\cnnm#1>{}
\def\cnnu#1>{} \def\vvu#1>{\ignore} \def\vvnu#1>{\ignore} 

\ifMag \ifUS   
  \let\vvu\vv \let\vvnu\vvn \let\cnnu\cnn  \else
 \let\goodbm\goodbreak  \let\mmgood\vvgood \let\cnnm\cnn
 \let\mline\nl \let\vvm\vv \let\vvnm\vvn  \fi
 \let\goodbreak\relax  \let\vvgood\relax
  \def\nngood{}  \fi

\def\wwgood#1:#2>{\vv#1>\vvgood\vv#2>\vv0>}
\def\vskgood#1:#2>{\vsk#1>\goodbreak\vsk#2>\vsk0>}

\def\wwmgood#1:#2>{\ifMag\vv#1>\mmgood\vv#2>\vv0>\fi}
\def\vskmgood#1:#2>{\ifMag\vsk#1>\goodbm\vsk#2>\vsk0>\fi}
\def\vskm#1:#2>{\ifMag\vsk#1>\else\vsk#2>\fi}
\def\vvmm#1:#2>{\ifMag\vv#1>\else\vv#2>\fi}
\def\vvnn#1:#2>{\ifMag\vvn#1>\else\vvn#2>\fi}
\def\nnm#1:#2>{\ifMag\nn#1>\else\nn#2>\fi}

\newif\ifMPIM
\ifx\MPIMpreprint\undefined\else\MPIMtrue\fi \ifMag\MPIMfalse\fi

  \let\mmpgood\relax \def\vskmp#1>{}
\def\vskmpgood#1:#2>{} \def\vvnmp#1>{}

\ifMPIM   \let\mmpgood\vvgood
 \let\vvnmp\vvn \let\vskmp\vsk \let\vskmpgood\vskgood\fi

\def\AMS/{Amer\. Math\. Soc.{}}

\whattime\readldf\writeldf

\csname beta.def\endcsname

\labeldef{F} {1\labelsep \labelspace 1}  {fh}
\labeldef{F} {1\labelsep \labelspace 2}  {albt}
\labeldef{F} {1\labelsep \labelspace 3}  {R}
\labeldef{F} {1\labelsep \labelspace 4}  {Rhh}
\labeldef{F} {1\labelsep \labelspace 5}  {inv}
\labeldef{F} {1\labelsep \labelspace 6}  {DYB}

\labeldef{S} {2} {E}
\labeldef{F} {2\labelsep \labelspace 1}  {Tph}
\labeldef{F} {2\labelsep \labelspace 2}  {Thh}
\labeldef{F} {2\labelsep \labelspace 3}  {RTT}
\labeldef{F} {2\labelsep \labelspace 4}  {Tmu}
\labeldef{F} {2\labelsep \labelspace 5}  {detT}
\labeldef{L} {2\labelsep \labelspace 1}  {qdet}
\labeldef{F} {2\labelsep \labelspace 6}  {TL}
\labeldef{F} {2\labelsep \labelspace 7}  {Dl}

\labeldef{S} {3} {e}
\labeldef{F} {3\labelsep \labelspace 1}  {la12}
\labeldef{F} {3\labelsep \labelspace 2}  {tf}
\labeldef{F} {3\labelsep \labelspace 3}  {tbc}
\labeldef{F} {3\labelsep \labelspace 4}  {tab}
\labeldef{F} {3\labelsep \labelspace 5}  {tac}
\labeldef{F} {3\labelsep \labelspace 6}  {Tcu}
\labeldef{L} {3\labelsep \labelspace 1}  {Rtt}
\labeldef{F} {3\labelsep \labelspace 7}  {RTTc}
\labeldef{F} {3\labelsep \labelspace 8}  {detTc}
\labeldef{F} {3\labelsep \labelspace 9}  {dett}
\labeldef{L} {3\labelsep \labelspace 2}  {qdets}
\labeldef{F} {3\labelsep \labelspace 10} {faut}
\labeldef{L} {3\labelsep \labelspace 3}  {Caut}
\labeldef{F} {3\labelsep \labelspace 11} {Cart}
\labeldef{F} {3\labelsep \labelspace 12} {ell}
\labeldef{F} {3\labelsep \labelspace 13} {vect}
\labeldef{L} {3\labelsep \labelspace 4}  {ehom}
\labeldef{L} {3\labelsep \labelspace 5}  {submod}
\labeldef{L} {3\labelsep \labelspace 6}  {inj}

\labeldef{S} {4} {H}
\labeldef{L} {4\labelsep \labelspace 1}  {PBW}
\labeldef{F} {4\labelsep \labelspace 1}  {tth}
\labeldef{F} {4\labelsep \labelspace 2}  {tQ}
\labeldef{F} {4\labelsep \labelspace 3}  {Detv}
\labeldef{L} {4\labelsep \labelspace 2}  {hwred}
\labeldef{L} {4\labelsep \labelspace 3}  {hwDet}
\labeldef{F} {4\labelsep \labelspace 4}  {Dethw}
\labeldef{L} {4\labelsep \labelspace 4}  {hwirr}
\labeldef{F} {4\labelsep \labelspace 5}  {BNc}
\labeldef{L} {4\labelsep \labelspace 5}  {idb}
\labeldef{L} {4\labelsep \labelspace 6}  {idn}
\labeldef{L} {4\labelsep \labelspace 7}  {Mmu}
\labeldef{L} {4\labelsep \labelspace 8}  {MVv}
\labeldef{L} {4\labelsep \labelspace 9}  {hwqu}
\labeldef{L} {4\labelsep \labelspace 10} {irrmu}
\labeldef{L} {4\labelsep \labelspace 11} {irrV}

\labeldef{S} {5} {S}
\labeldef{L} {5\labelsep \labelspace 1}  {HC}
\labeldef{F} {5\labelsep \labelspace 1}  {qf}
\labeldef{F} {5\labelsep \labelspace 2}  {pimth}
\labeldef{L} {5\labelsep \labelspace 2}  {etapim}
\labeldef{F} {5\labelsep \labelspace 3}  {Smm}
\labeldef{F} {5\labelsep \labelspace 4}  {SmQv}
\labeldef{L} {5\labelsep \labelspace 3}  {SmQ}
\labeldef{L} {5\labelsep \labelspace 4}  {Smv}
\labeldef{L} {5\labelsep \labelspace 5}  {NmQ}
\labeldef{L} {5\labelsep \labelspace 6}  {VmQ}
\labeldef{L} {5\labelsep \labelspace 7}  {Virr}
\labeldef{L} {5\labelsep \labelspace 8}  {fdim}
\labeldef{L} {5\labelsep \labelspace 9}  {dimu}

\labeldef{S} {6} {cg}
\labeldef{L} {6\labelsep \labelspace 1}  {vmQa}
\labeldef{L} {6\labelsep \labelspace 2}  {MM}
\labeldef{L} {6\labelsep \labelspace 3}  {kerpi}
\labeldef{F} {6\labelsep \labelspace 1}  {CB}
\labeldef{L} {6\labelsep \labelspace 4}  {Cvv}
\labeldef{L} {6\labelsep \labelspace 5}  {kerC}
\labeldef{L} {6\labelsep \labelspace 6}  {SQQ}
\labeldef{F} {6\labelsep \labelspace 2}  {SQS}
\labeldef{L} {6\labelsep \labelspace 7}  {CSQ}
\labeldef{F} {6\labelsep \labelspace 3}  {CS}
\labeldef{F} {6\labelsep \labelspace 4}  {CSt}
\labeldef{F} {6\labelsep \labelspace 5}  {Cfg}

\labeldef{S} {7} {R}
\labeldef{F} {7\labelsep \labelspace 1}  {tthr}
\labeldef{F} {7\labelsep \labelspace 2}  {ellr}
\labeldef{F} {7\labelsep \labelspace 3}  {Rb}
\labeldef{F} {7\labelsep \labelspace 4}  {RTTb}
\labeldef{F} {7\labelsep \labelspace 5}  {dlr}
\labeldef{F} {7\labelsep \labelspace 6}  {tN}
\labeldef{F} {7\labelsep \labelspace 7}  {eab}
\labeldef{F} {7\labelsep \labelspace 8}  {sln}
\labeldef{F} {7\labelsep \labelspace 9}  {ellh}
\labeldef{L} {7\labelsep \labelspace 1}  {ptndg}
\labeldef{L} {7\labelsep \labelspace 2}  {pertox}
\labeldef{L} {7\labelsep \labelspace 3}  {CVsv}
\labeldef{L} {7\labelsep \labelspace 4}  {CcV}
\labeldef{L} {7\labelsep \labelspace 5}  {EcV}
\labeldef{L} {7\labelsep \labelspace 6}  {EcUV}
\labeldef{L} {7\labelsep \labelspace 7}  {kerS}
\labeldef{L} {7\labelsep \labelspace 8}  {S8}

\labeldef{S} {8} {Fd}
\labeldef{L} {8\labelsep \labelspace 1}  {abba}
\labeldef{L} {8\labelsep \labelspace 2}  {ne0}
\labeldef{L} {8\labelsep \labelspace 3}  {cdba}
\labeldef{L} {8\labelsep \labelspace 4}  {regsi}
\labeldef{F} {8\labelsep \labelspace 1}  {aa1}
\labeldef{F} {8\labelsep \labelspace 2}  {a1a}
\labeldef{L} {8\labelsep \labelspace 5}  {infdim}
\labeldef{L} {8\labelsep \labelspace 6}  {aksub}
\labeldef{L} {8\labelsep \labelspace 7}  {ZiS}
\labeldef{L} {8\labelsep \labelspace 8}  {ZS}
\labeldef{L} {8\labelsep \labelspace 9}  {te}
\labeldef{L} {8\labelsep \labelspace 10} {ZSr}
\labeldef{F} {8\labelsep \labelspace 3}  {dims}
\labeldef{F} {8\labelsep \labelspace 4}  {Agm}
\labeldef{F} {8\labelsep \labelspace 5}  {AB}
\labeldef{L} {8\labelsep \labelspace 11} {kerim}
\labeldef{F} {8\labelsep \labelspace 6}  {diml}
\labeldef{F} {8\labelsep \labelspace 7}  {dmn1}

\labeldef{S} {9} {J}
\labeldef{L} {9\labelsep \labelspace 1}  {abrr}
\labeldef{F} {9\labelsep \labelspace 1}  {Jc}
\labeldef{F} {9\labelsep \labelspace 2}  {Jch}
\labeldef{L} {9\labelsep \labelspace 2}  {Jccl}
\labeldef{F} {9\labelsep \labelspace 3}  {2ccl}
\labeldef{F} {9\labelsep \labelspace 4}  {lazh}
\labeldef{F} {9\labelsep \labelspace 5}  {Rc}
\labeldef{L} {9\labelsep \labelspace 3}  {Rc3}
\labeldef{L} {9\labelsep \labelspace 4}  {J}
\labeldef{L} {9\labelsep \labelspace 5}  {JUVW}
\labeldef{F} {9\labelsep \labelspace 6}  {J4}
\labeldef{F} {9\labelsep \labelspace 7}  {RJ}
\labeldef{L} {9\labelsep \labelspace 6}  {RUVW}
\labeldef{F} {9\labelsep \labelspace 8}  {RtVW}
\labeldef{L} {9\labelsep \labelspace 7}  {JJ}
\labeldef{F} {9\labelsep \labelspace 9}  {JJt}
\labeldef{L} {9\labelsep \labelspace 8}  {RRt}
\labeldef{F} {9\labelsep \labelspace 10} {Rti}
\labeldef{F} {9\labelsep \labelspace 11} {RUU}
\labeldef{L} {9\labelsep \labelspace 9}  {Ec}

\labeldef{S} {10} {ex}
\labeldef{L} {10\labelsep \labelspace 1}  {ttt}
\labeldef{F} {10\labelsep \labelspace 1}  {LL}
\labeldef{F} {10\labelsep \labelspace 2}  {Lf}
\labeldef{F} {10\labelsep \labelspace 3}  {RLL}
\labeldef{F} {10\labelsep \labelspace 4}  {detL}
\labeldef{F} {10\labelsep \labelspace 5}  {det=1}
\labeldef{F} {10\labelsep \labelspace 6}  {TpL}
\labeldef{L} {10\labelsep \labelspace 2}  {L2t}
\labeldef{F} {10\labelsep \labelspace 7}  {tL}
\labeldef{L} {10\labelsep \labelspace 3}  {t2L}
\labeldef{L} {10\labelsep \labelspace 4}  {cat}
\labeldef{F} {10\labelsep \labelspace 8}  {Gc}
\labeldef{L} {10\labelsep \labelspace 5}  {Ect}

\labeldef{S} {\char 65} {CR}
\labeldef{F} {\char 65\labelsep \labelspace 1}  {tbd}
\labeldef{F} {\char 65\labelsep \labelspace 2}  {thab}
\labeldef{F} {\char 65\labelsep \labelspace 3}  {thac}
\labeldef{F} {\char 65\labelsep \labelspace 4}  {thbd}
\labeldef{F} {\char 65\labelsep \labelspace 5}  {thabcd}
\labeldef{F} {\char 65\labelsep \labelspace 6}  {thcdab}
\labeldef{F} {\char 65\labelsep \labelspace 7}  {thadcb}
\labeldef{F} {\char 65\labelsep \labelspace 8}  {Sab}
\labeldef{F} {\char 65\labelsep \labelspace 9}  {Scd}
\labeldef{L} {\char 65\labelsep \labelspace 1}  {rels}
\labeldef{F} {\char 65\labelsep \labelspace 10} {thabcdk}
\labeldef{F} {\char 65\labelsep \labelspace 11} {thadcbk}

\labeldef{S} {\char 66} {Det}
\labeldef{L} {\char 66\labelsep \labelspace 1}  {kerQ}
\labeldef{F} {\char 66\labelsep \labelspace 1}  {TTR}
\labeldef{F} {\char 66\labelsep \labelspace 2}  {TTS}
\labeldef{F} {\char 66\labelsep \labelspace 3}  {ATT}
\labeldef{L} {\char 66\labelsep \labelspace 2}  {TNl}
\labeldef{F} {\char 66\labelsep \labelspace 4}  {TN}
\labeldef{F} {\char 66\labelsep \labelspace 5}  {Avj}
\labeldef{L} {\char 66\labelsep \labelspace 3}  {Tab}
\labeldef{L} {\char 66\labelsep \labelspace 4}  {TTT}
\labeldef{F} {\char 66\labelsep \labelspace 6}  {TT=T}
\labeldef{L} {\char 66\labelsep \labelspace 5}  {RN}

\labeldef{S} {\char 67} {mform}

\labeldef{S} {\char 68} {norm}
\labeldef{L} {\char 68\labelsep \labelspace 1}  {pbw}
\labeldef{L} {\char 68\labelsep \labelspace 2}  {uni}

\labeldef{S} {\char 69} {sl2}
\labeldef{F} {\char 69\labelsep \labelspace 1}  {Sfg}

\labeldef{S} {\char 70} {B}
\labeldef{L} {\char 70\labelsep \labelspace 1}  {BVW}
\labeldef{L} {\char 70\labelsep \labelspace 2}  {RBB}
\labeldef{L} {\char 70\labelsep \labelspace 3}  {BU}

\document

\hoffset 0pt
\voffset 0pt

\hfuzz15pt

\center
\hrule height 0pt
\vsk.7>

{\twelvepoint\bf \bls1.2\bls
Small Elliptic Quantum Group $\esl$
\par}
\vskm1.5:1.4>
\vskmp.6>
\=
\VT/$^{\,\star}$ \ and \ \Varch/$^{\,*}$
\ifMPIM\else
\vskm1.5:1.4>
{\it
$^\star$\MPIM/, \MPIMaddr/
\vsk.3>
$^*$\UNC/\\ \UNCaddr/
\vskm1.6:1.5>
\sl November \,2000}
\fi
\endcenter

\ftext{\=\bls11pt
$\]^\star\)$\ifMPIM\else\absence/\vv-.06>\nl
\hp{$^*$}\fi
\support/ RFFI grant 99\)\~\)01\~\)00101
\>and \,INTAS grant 99\)\~\)01705\vv.06>\nl
\hp{$^*$}{\tenpoint\sl E-mail\/{\rm:} \mpimemail/\,, \;\homemail/}\vv.1>\nl
${\]^*\)}$\support/ NSF grant DMS\)\~\)9801582\vv.06>\nl
\vv-1.2>
\hp{$^*$}{\tenpoint\sl E-mail\/{\rm:} \avemail/}}

\vskm1.7:1.5>
\vskmp.56>

{\ifMag\ninepoint\fi
\Abstract
The small \eqg/ $\esl$, introduced in the paper, is an elliptic dynamical
analogue of the universal enveloping algebra $\Usl$. We define \hwm/s,
\Vmod/s and \cmod/s over $\esl$, the \dSform/ for $\esl$ and the \cform/ for
\hw/ \emod/s. We show that any \fd/ \smod/ and any \Vmod/ over $\sln$ can be
lifted to the corresponding \emod/ on the same vector space. For the \eqg/
$\Eqg$ we construct the \eval/ morphism ${\Eqg\to\)\esl}$, thus making any
\emod/ into an \eval/ \Emod/.
\endAbs}

\Sno -1

\vskm1.5:1.1>
\vskmp.2>
\vsk0>

\sect{Introduction}
The main purpose of this paper is to define a dynamical \qg/ $\esl$ which is
an elliptic dynamical analogue of the universal enveloping algebra $\Usl$.
We call $\esl$ the small \eqg/, comparing it with the \eqg/ $\Eqg$ introduced
in \Cite{F}\). Our initial motivation to study this object arises from the
wish to understand the structure of \evmod/s over $\Eqg$, which should be
analogous to \evmod/s over the Yangian $\Ysl$.
\vsk.2>
Evaluation modules over $E_{\tau\],\)\gm}(\frak{sl}_2)$ have been defined in
\Cite{FV1}\). They appear naturally in the descripion of transition matrices
for the \tri/ \qKZ/ \deq/ \Cite{TV1}\). They also serve for the definiton of
the \qKZB/ \deq/s and occur in the description of its monodromies,
see \Cite{FTV}\). One should expect \evmod/s over $\Eqg$ for ${N>2}$ to play
a similar role. Symmetric and exterior powers of the vector \rep/ of $\Eqg$,
developed in \Cite{FV2}\), are examples of \evmod/s over $\Eqg$.
In general, \evmod/s over $\Eqg$ arise from \emod/s via the \eval/ morphism
${\Eqg\to\)\esl}$, see Corollary~\[ehom]\), analogous to the \ehom/
$\Ysl\to\)\Usl$.
\vsk.2>
In this paper we prove a PBW type theorem for the small \eqg/ $\esl$.
We define \hwm/s, \Vmod/s and \cmod/s over $\esl$. We show that for any \fd/
\smod/ and any \Vmod/ over $\sln$ one can define the corresponding \emod/ on
the same vector space. Pulling back these \emod/s through the \eval/ morphism
we get \fd/ \evmod/s and \eVmod/s over $\Eqg$. Conjecturally, the same picture
takes place for any \hw/ \smod/.
\vsk.2>
We introduce the \dSform/ for $\esl$, the \dSpair/ and the \cform/ for
the \hwm/s over $\esl$. They play an important role in the construction of
\fd/ \hw/ \emod/s. From another point of view the \cform/ for \$\eslt$-modules
appeared in a disguised form in \Cite{\)TV1\), Appendix C\)}\).
\vsk.2>
The small \eqg/ $\esl$ admits the \tri/ and \rat/ degenerations. They are
closely related to the exchange \qg/s $F_q\bigl(SL(N)\bigr)$ and $\Fsl$
introduced in \Cite{EV2}\). In this paper we consider only the \rat/ \dqg/
$\eslr$ and its relation to the exchange \qg/ $\Fsl$. We contstruct a functor
from a certain category of \smod/s to a category of \ermod/s,
see Theorems~\[EcV]\), \[EcUV]\), and \[Ec]\). We also establish an equivalence
of certain tensor categories of \ermod/s and \rat/ dynamical \rep/s of $\Fsl$,
see Theorem~\[cat]\). In particular, this gives a new construction of \hw/
\rep/s for the exchange \qg/ $\Fsl$.
\vsk.2>
Notice that while $\eslr$ and exchange \qg/s have a coproduct structure,
no coproduct structure is known for the \eqg/ $\esl$. This makes the small
elliptic group similar to the Sklyanin algebra \Cite{S}\), \Cite{HW}\).
\vsk.2>
The paper is organized as follows. After introducing basic notation we recall
the definiton of the \eqg/ $\Eqg$. The small \eqg/ $\esl$ is defined in
Section~\[:e]\). In Section~\[:H] we introduce \hwm/s and \Vmod/s over $\esl$.
The \dSform/ for $\esl$ is defined in Section~\[:S]\). Contragradient modules
and the \cform/ for \hw/ \emod/s are defined in Section~\[:cg]\).
In Section~\[:R] we study the \rat/ \dqg/ $\eslr$. We construct \irr/ \fd/
\emod/s in Section~\[:Fd]\). A functor from a certain category of \smod/s
to a category of \ermod/s is defined in Section~\[:J]\). Relations between
$\eslr$ and $\Fsl$ are studied in Section~\[:ex]\). There are six Appendices in
the paper; they contain useful technichal information, the $\eslt$ example,
and some proofs.
\vsk.3>
The first author would like to thank the \MPIM/ in Bonn,
where he stayed when this paper was being written, for hospitality.

\Sect{Basic notation}
Let $\tau$ be a complex number \st/ $\Im\tau>0$.
Let $\tht(u\);\tau)$ be the Jacobi theta \fn/:
\vvnn-.6:-.2>
$$
\tht(u\);\tau)\,=\,-\!\sum_{m=-\8}^\8\!
\exp\>\bigl(\)\pii\)\tau\)(m+1/2\))^2+2\pii\>(m+1/2\))\>(u+1/2\))\bigr)\,.
\vvmm-.1:-.3>
$$
There is a product formula
\vvnn-.4:-.5>
$$
\tht(u\);\tau)\,=\,2\>e^{\pii\)\tau\]/4}\>\sin\)(\pi\)u)\,\prod_{s=1}^\8\)
(1-e^{2\pii\)s\tau})\>(1-e^{2\pii\)(s\tau+u)})\>(1-e^{2\pii\)(s\tau-u)})\,.
\vv-.1>
\wwmgood-.6:.6>
$$
The \fn/ $\tht(u\);\tau)$ has multipliers $-1$ and
${-\)\exp\>(-2\>\pii\)u-\pii\)\tau)}$ as ${u\)\to\)u+1}$ and
${u\)\to\)u+\tau}$, \resp/. It is an entire \fn/ with only simple zeros lying
on the lattice $\ZtZ$. Usually, we omit the second argument of the theta \fn/,
writing $\tht(u)$ instead of $\tht(u\);\tau)$.
\vsk.3>
Let $\hg$ be a \fd/ commutative Lie algebra, and let $\hga\}$ be the dual
space. An \hmod/ $V\}$ is called \em{\diag/} if it admits a \wtd/
\vvn-.1>
$$
V\)=\>\Plus_{\,\mu\>\in\>\hga\!}\Vmu\,,
\vv-.4>
$$
all \wt/ subspaces $\Vmu$ being \fd/ and the set
\vv.05>
\)$\lb\)\mu\vert\Vmu\ne 0\)\rb$ \)at most countable.
\vsk.3>
Let $V_1\lc V_k$ be \hmod/s. For any \fn/ ${f\]:\hga\to\End(V_1\lox V_k)}$
and any ${i=1\lc k}$ we define an operator $f(h\"i)\in\End(V_1\lox V_k)$
by the rule:
\ifMag
$$
f(h\"i)\>\vox\)=\)f(\)\mu)\>\vox
\Tag{fh}
$$
for any $v\in V_1\lox V_i\)[\)\mu\)]\lox V_k$.
\vsk.2>
\else
\vvn-.2>
$$
f(h\"i)\>\vox\)=\)f(\)\mu)\>\vox\qquad\text{for any}\quad
v\in V_1\lox V_i\)[\)\mu\)]\lox V_k\,.\kern-2em
\Tag{fh}
$$
\par
\fi
{\bls1.06\bls
For a \fd/ vector space $V\}$ over $\C$ denote by $\FunV$ the space of \Vval/
\mef/s on $\hga\}$. If $V\}$ is a \dhmod/, set
$$
\FunV\,=\>\Plus_{\,\mu\>\in\>\hga\!}\FVmu\,.
\vv-.4>
$$
The space $\FunV$ is a vector space over $\FunC$. The space $V\}$ is naturally
embedded in $\FunV$ as the subspace of constant \fn/s. If $V\}$ is an \hmod/,
then $\FunV$ is an \hmod/ with the natural pointwise action of $\hg$ and
$\FunV[\)\mu\)]=\Fun\bigl(\Vmu\bigr)$.
\vsk.2>}
Let $U\]$ be a \dhmod/ and a vector space over $\FunC$. Suppose that the action
of $\hg$ commutes with multiplication by \fn/s. Then each \wt/ subspace
$\Umu$ is a vector space over $\FunC$. Assume that all the \wt/ subspaces
are \fd/ over $\FunC$. Then one can define a \dhmod/ $V\}$, \st/ $U=\FunV$
as \hmod/s, in the following way. For any $\mu$ \st/
\vvm-.08>
$\Umu\ne 0$, pick up a basis $f_1\lc f_k$ of $\Umu$ over $\FunC$ and set
$\Vmu\)=\Plus_{i=1}^k\C f_i$, \,otherwise, set $\Vmu=0$. Then define
\vvmm-.04:-.4>
${V=\!\Plus_{\,\mu\>\in\>\hga\!}\Vmu}$ \)to be the \dhmod/ \st/ $\Vmu$
is a \wt/ subspace of \wt/ $\mu$.
\vskmpgood-.8:.8>
\vsk.2>
Let $V,\)W\}$ be \dhmod/s. The space $\Hom(V,W)$ has the natural \hmod/
structure, but in general the \wt/ subspaces are in\fd/. We set
\vvnm-.15>
$$
\Fun\bigl(\Hom(V,W)\bigr)=\Hom\bigl(V,\FunW\bigr)\,.
\vvm.15>
$$
A \fn/ ${\phi\in\Fun\bigl(\Hom(V,W)\bigr)}$ induces a linear map
$\FunV\to\FunW$, \,acting pointwise: $f(\la)\)\map\)\phi(\la)\)f(\la)$.
This map is usually denoted by the same letter.
\vsk.2>
Denote by $\Diff(V)$ the space of \dif/ operators acting in $\FunV$.
It is spanned over $\C$ by operators of the form
$f(\la)\,\map\,\phi(\la)\>f(\la+\mu)$ where $\phi\in\Fend V$
and $\mu\in\hga\}$.
\vsk.2>
As a rule we do not distinguish a \fn/ $\phi(\la)\in\FunC$ and
the \fn/ $\phi(\la)\}\cdot\]\id\in\Fend V$.
\vskmgood-.8:.8>
\vsk.3>
In this paper we take $\hg$ to be the Cartan subalgebra of the Lie algebra
$\sln$. Fix a basis $h_1\lc h_{N-1}$ of $\hg$. Let $\om_1\lc\om_{N-1}$ be the
fundamental \wt/s: $\bra\)\om_a\>,h_b\)\ket\)=\)\dl_{ab}$. Let $\PP\sub\hga$
be the \wt/ lattice: $\PP=\]\Plus_{a=1}^{N-1}\Z\,\om_a$.
For any $\ano$, set $\eps_a\]=\om_a\]-\om_{a-1}$, where by convention
\vvmm0:-.6>
$\om_{\)0}\]=\om_N=0$. Let $\al_1\lc\al_{N-1}$ be the simple roots:
$\al_a\]=\eps_a\]-\eps_{a+1}$. For $\la\),\mu\in\hga\}$ say that $\la\ge\mu$
if $\la-\mu\)\in\}\Plus_{a=1}^{N-1}\Zp\,\al_a$.
\vsk.2>
{\bls 1.1\bls
Define a bilinear form $(\,{,}\,)$ on $\hga\}$ by the rule
$(\al_a\),\om_b)=\dl_{ab}$ for any $a\),\bno-1$. For any $\la\in\hga\}$ set
\vvm-.5>
$\la_a=(\la\>,\eps_a)$. It is easy to see that $\la=\!\suan\]\la_a\>\eps_a$,
\,$\suan\]\la_a=0$ \,and \,$(\la\>,\mu)=\}\suan\]\la_a\)\mu_a$.
The Weyl group $\)W\}$ acts on $\hga\}$ as the \symg/ $\Sb_N$ permuting
the coordinates $\laN$.
\vskm.3:.2>
Let ${\rho\>=\]\sum_{a=1}^{N-1}\]\om_a=\)-\}\suan\]a\>\eps_a}$
\vvm.1>
be the half-sum of positive roots. For any ${w\in W}$ and ${\la\in\hga\}}$
set $w\cddt\la=w\)(\la+\rho)-\rho$. Notice that
$(\la\>,\rho)\)=\>-\,\la_1-2\)\la_2\lsym-N\la_N$.
\Par}
Let $E_{ab}\in\End(\CN)$ be the matrix with the only nonzero entry equal to $1$
at the intersection of the \text{$a\}$-th} row and \text{$b\]$-th} column.
The assignment \,${h_a\map E_{aa}\]-E_{a+1,\)a+1}}$, \,$\ano-1$, makes $\CN\!$
into an \hmod/, called the \em{vector \rep/} of $\hg$. Henceforth, we always
consider $\CN\!$ as the vector \rep/ of $\hg$.
\vsk.4>
Let $\gm$ be a nonzero complex number. Introduce \fn/s $\al(u,\xi)$ and
$\bt\)(u,\xi)$ as follows:
\vvnn-.5:.1>
$$
\al(u\),\xi\))\,=\;{\tht(u)\>\tht(\xi+\gm)\over\tht(u-\gm)\>\tht(\xi\))}\;,\qqq
\bt\)(u\),\xi\))\,=\,
-\;{\tht(u+\xi\))\>\tht(\gm)\over\tht(u-\gm)\>\tht(\xi\))}\;.\kern-1em
\Tag{albt}
$$
Let $R(u\),\la)$ be the \em{elliptic dynamical \Rm/} \cite{F}\):
\vvn-.1>
$$
R(u\),\la)\,=\}\tsuan E_{aa}\ox E_{aa}\>+\!
\tsum_{\tsize{a,b=1\atop a\ne b}}^N\}\bigl(\)
\al(u\),\la_{ab})\>E_{aa}\ox E_{bb}\>+\>
\bt\)(u\),\la_{ab})\>E_{ab}\ox E_{ba}\bigr)\,,\kern-1em
\vv-.2>
\Tag{R}
$$
where ${\la\in\hga}$ and $\la_{ab}=\la_a\]-\la_b$.
The dynamical \Rm/ has zero \wt/:
$$
\bigl[\)R(u\),\la)\,,\)h\"1\}+\)h\"2\)\bigr]\,=\,0\,,\kern-1em
\vv-.1>
\Tag{Rhh}
$$
satisfies the inversion relation:
\vvn-.3>
$$
R(u\),\la)\>R\"{21}(-u\),\la)\,=\,1\,,\kern-1.3em
\vv-.1>
\Tag{inv}
$$
and the \DYB/:
$$
\align
R\"{12}(u-v,\la-\gm\)h\"3)\> & R\"{13}(u\),\la)\>R\"{23}(v,\la-\gm\)h\"1)\,={}
\Tag{DYB}
\\
\nn6>
{}=\,{} & R\"{23}(v,\la)\>R\"{13}(u\),\la-\gm\)h\"2)\>R\"{12}(u-v,\la)\,.
\kern-1.5em
\endalign
$$
The last equality holds in $\End(\CN\!\ox\CN\!\ox\CN)$. By standard convention,
we assume that $R\"{ij}(u\),\la)$ acts as $R(u\),\la)$ on the \}\$i$-th and
\vvgood
\}\$j$-th tensor factors and as the identity operator on the remaining factors.
\mmgood
For instance, in formula \(DYB) we have ${R\"{12}\]=\)R\ox\id}$ \)and
${R\"{23}\]=\>\id\ox R}$. Notice that for the \Rm/ \(R) in addition we have
\vvnn-.1:0>
\vvnmp.2>
$$
R\bigl(u\),\la-\gm\>(h\"1\!+h\"2)\bigr)\,=\,R(u\),\la)\,.
\mmpgood
\vv.1>
$$

\Sect[E]{Elliptic \qg/ $\{\Eqg$}
A \em{module over the \eqg/} $\Eqg$ is a \dhmod/ $V\}$ together with
\$\Diff(V)\)$-valued \mef/s \,$T_{ab}(u)$, \,$a\),\bno$, \,in a complex \var/
$u$, subject to relations \(Tph)\,--\,\(RTT)\). We combine the \fn/s
$T_{ab}(u)$ into a matrix $T(u)$ with noncommuting entries:
$$
T(u)\,=\>\sum_{a,b}\>E_{ab}\ox T_{ab}(u)\,.
\vv-1.2>
$$
The defining relations are:
\vvn-.3>
$$
T_{ab}(u)\,\phi(\la)\,=\,\phi(\la-\gm\>\eps_b)\,T_{ab}(u)
\vv-.3>
\Tag{Tph}
$$
for any $\phi\in\FunC$,
\vvn-.7>
$$
\gather
\bigl[\)T(u)\,,\)h\"1\}+\)h\"2\)\bigr]\,=\,0\,,
\Tag{Thh}
\\
\nn10>
R\"{12}(u-v,\la-\gm\)h\"3)\>T\"{13}(u)\>T\"{23}(v)\,=\,
T\"{23}(v)\>T\"{13}(u)\>R\"{12}(u-v,\la)\,.\kern-1em
\Tag{RTT}
\\
\cnn-.1>
\endgather
$$
The last equality holds in $\End\bigl(\CN\!\ox\CN\!\ox\FunV\bigr)$. Here
\vv-.4>
\,$T\"{13}(u)=\sum_{a,b}\)E_{ab}\ox\id\ox T_{ab}(u)$ \,and
\,$T\"{23}(u)=\botsmash{\sum_{a,b}\>\id\ox E_{ab}\ox T_{ab}(u)}$.
\vsk.9>
Relations \(Tph) can be written as
\,${T(u)\,\phi(\la+\gm\)h\"1)\>=\,\phi(\la)\,T(u)}$
\,\)for any ${\phi\in\FunC}$. Formula \(Thh) means that for any $\mu\in\hga$
$$
T_{ab}(u)\,\FVmu\,\sub\,
\Fun\bigl(V[\)\mu-\eps_a\]+\eps_b\)]\)\bigr)\,.
\Tag{Tmu}
$$
\par
Introduce the \em{quantum determinant} $\Det\)T(u)$,
\cf. \Cite{FV1}\), \Cite{FV2}\), by the rule
$$
\Det\)T(u)\,=\;{\Tht(\la)\over\Tht(\la-\gm\)h)}\,
\sumib\sign(\)\ib\))\;T_{i_N\],\)N}\bigl(u+(N-1)\>\gm\bigr)
\>\ldots\,T_{i_2,\)2}(u+\gm)\,T_{i_1,\)1}(u)\kern-2em
\vv-.2>
\Tag{detT}
$$
where \,$\dsize\Tht(\la)\>=\!\!\prod_{\abn}\!\tht(\la_a\]-\la_b)$,
\vv.05>
\;the sum is taken over all \perm/s $\ib=(i_1\lc i_N)$, and $\sign(\)\ib\))$
is the sign of the \perm/. It is clear that $\Det\)T(u)$ commutes with
\vv.03>
multiplication by any \fn/ $\phi(\la)\in\FunC$ and with the action of $\hg$.
Hence, $\Det\)T(u)$ acts on $\FunV$ as multiplication by an \$\End(V)$-valued
\mef/ of $u$ and $\la$. We denote this \fn/ by $\Det L(u\),\la)$, \cf. \(TL)\).
\Prop{qdet}
$\bigl[\>\Det\)T(u)\>,\)T_{ab}(v)\)\bigr]\,=\,0\,$ for any $a\),\bno$.
\endpro
\nt
The proposition is proved in Appendix~\[:Det]\).
\Par
According to \(Tph), $T_{ab}(u)$ is a \dif/ operator; for any $v\in\FunV$
we have
$$
\bigl(\)T_{ab}(u)\,v\bigr)(\la)\,=\,L_{ab}(u\),\la)\,v(\la-\gm\>\eps_b)
\Tag{TL}
$$
where $L_{ab}(u\),\la)\in\Fun\bigl(\End(V)\bigr)$.
Set \,$L(u\),\la)\,=\)\sum_{a,b}\>E_{ab}\ox L_{ab}(u\),\la)$.
\vskmgood-.8:.8>
\vsk-.3>
\vsk0>
\Ex
For any ${x\in\C}$ the assignment ${L(u\),\la)\>\map R(u-x,\la)}$ \,makes
$\CN\!$ into an \Emod/. This module is called the \em{vector \rep/} of $\Eqg$
with the \em{\epoint/} $x$. The quantum determinant in this \Emod/ is
\vvnn-.3:-.5>
$$
\Det L(u\),\la)\,=\;{\tht\bigl(u-x+(N-1)\>\gm\bigr)\over\tht(u-x)}\;.
$$
\enddemo
\vsk-.2>
Let $V,\)W\}$ be \Emod/s. An element ${\phi\in\Fun\bigl(\Hom_\hg(V,W)\bigr)}$
is called a \em{morphism} of \Emod/s if the induced map ${\FunV\to\FunW}$
satisfies
$$
\phi(\la)\,T_{ab}(u)\vst{\italic{Fun}(V)}=
\,T_{ab}(u)\vst{\italic{Fun}(W)}\>\phi(\la)
\vvm-.1>
\wwgood-.6:.7>
$$
for any $a\),\bno$. Denote by $\Mor(V,W)$ the space of all morphisms from $V\}$
to $W\}$. A morphism $\phi$ is called an \em{\iso/} if the map $\phi(\la)$ is
bijective for generic $\la$.
\par
An \Emod/ $V\}$ is called \em{\irr/} if for any nontrivial morphism
$\phi\in\Mor(W,V)$ the map map $\phi(\la)$ is surjective for generic $\la$,
and \em{reducible} otherwise.
\vsk.3>
Let $V,\)W\}$ be \Emod/s. Then the \hmod/ $V\]\ox W\}$ is made into an \Emod/
by the rule
\vvnn-.8:-.9>
$$
L_{ab}(u\),\la)\vst{V\ox\)W}\,=
\>\sum_{c=1}^N\,L_{ac}(u\),\la-\gm\)h\"2)\ox L_{cb}(u\),\la)\,,
\vv-.1>
\Tag{Dl}
$$
and $T_{ab}(u)$ acts on $\Fun(V\]\ox W)$ according to \(TL)\). Triple tensor
products ${(U\ox V)\ox W}\}$ and ${U\ox(V\]\ox W)}$ are canonically isomorphic
as \Emod/s. The quantum determinant is group-like, it acts on the \Emod/
$V\]\ox W\}$ by
\vvnn-.2:-.4>
$$
\Det L(u\),\la)\vst{V\ox\)W}\>=\,
\Det L(u\),\la-\gm\)h\"2)\)\ox\>\Det L(u\),\la)\,.
\vv-.1>
$$
If $\phi_1,\phi_2$ are morphisms of \Emod/s, $\phi_1\in\Mor(V_1,W_1)$,
$\phi_2\in\Mor(V_2,W_2)$, then
\vvnm-.2>
$$
\phi_1(\la-\gm\)h\"2)\)\ox\)\phi_2(\la)
\vvm.2>
$$
is a morphism of \Emod/s $V_1\ox V_2$ and $W_1\ox W_2$.
\Rem
Notice that the given definition of modules over the \eqg/ $\Eqg$
is slightly different from the corresponding definition in \Cite{FV1}\).
The definition of morphisms of \Emod/s is also suitably modified.
We choose the present version in order to simplify the exposition.
\enddemo

\Sect[e]{Small \eqg/ $\{\esl$}
Let \>$\Funt$ be the ring of \mef/s $f(\la\+1\},\la\+2)$
\vv.1>
on $\hga\!\oplus\hga\}$ \st/ location of singularities of
$f(\la\+1\},\la\+2)$ in $\la\+1\}$ does not depend on $\la\+2\}$
and vice versa.
\vv.1>
For brevity, we write $f(\la\+1)$ or $f(\la\+2)$ instead of
$f(\la\+1\},\la\+2)$ if the \fn/ does not depend on the other \var/.
\vsk.4>
Given a \dhmod/ $V\}$ we define an action of $\Funt$ in $\FunV$:
for any $f\}\in\Funt$ set $f:\)v(\la)\>\map\>f(\la\>,\la-\gm\)h)\,v(\la)$.
For instance, for any $\phi\in\FunC$ we have
\vvnn-.1:.1>
$$
\phi(\la\+1)\):\>v(\la)\>\map\>\phi(\la)\,v(\la)\,,\qqq
\phi(\la\+2)\):\>v(\la)\>\map\>\phi(\la-\gm\)h)\,v(\la)\,.\kern-2em
\wwmgood-.3:.3>
\vv.1>
\Tag{la12}
$$
We always assume that $\Funt$ acts on $\FunV$ in this way.
\Par
The \em{\oalg/} $\els$ is a unital associative algebra over $\C$ \gby/
elements $t_{ab}$, $a\),\bno$, and \fn/s $f\}\in\Funt$ subject to relations
$$
t_{ab}\,f(\la\+1\},\la\+2)\,=\>
f(\la\+1\!-\gm\>\eps_a\>,\la\+2\!-\gm\>\eps_b)\,t_{ab}\kern-2em
\vv-.3>
\Tag{tf}
$$
for any $f\}\in\Funt$,
\vvn-.5>
$$
\gather
t_{ab}\>t_{ac}\,=\,t_{ac}\>t_{ab}\,,
\Tag{tbc}
\\
\nn10>
t_{ac}\>t_{bc}\,=\,\){\tht(\la\+1_{ab}\!+\gm)\over\tht(\la\+1_{ab}\!-\gm)}
\ t_{bc}\>t_{ac}\,,\rlap{\kern4em for\quad $a\ne b$\,,}
\Tag{tab}
\\
\nn12>
{\tht(\la\+2_{bd}\!+\gm)\over\tht(\la\+2_{bd})}\ t_{ab}\,t_{cd}\,-\,\)
{\tht(\la\+1_{ac}\!+\gm)\over\tht(\la\+1_{ac})}\ t_{cd}\,t_{ab}\,\)=\;
{\tht(\la\+1_{ac}\!+\la\+2_{bd})\,\tht(\gm)\over
\tht(\la\+1_{ac})\,\tht(\la\+2_{bd})}\ t_{ad}\,t_{cb}\,,\kern-2em
\Tag{tac}
\\
\cnn.2>
\endgather
$$
for \;$a\ne c$ \;and \;$b\ne d$. \;Here $\la\+i_{ab}=\la\+i_a\!-\la\+i_b\!$.
\vsk.3>
The ring $\Funt$ is embedded into $\els$ as a commutative subalgebra.
It acts on $\els$ by left multiplication. In this paper we consider $\els$
as the corresponding \$\Funt$-module.
\goodbreak
\vsk.3>
The \oalg/s $\els$, $\elsp$ and $\elsm$ are isomorphic.
\vv.06>
The \iso/ $\els\to\elsp$ corresponds to the property
\vv.1>
$\tht(u\);\tau)=\)e^{\pii/\]4}\,\tht(u\);\tau\]+1)$ of the theta \fn/ and
is tautological. The isomorphism $\els\to\elsm$ corresponds to the equality
\vvnn-.2:-.1>
$$
\tht(u\);-1/\]\tau)\,=\,(i\)\tau)^{1/2}\)
\exp\)(\)\pii\)\tau\)u^2\))\,\tht(-\)\tau\)u\);\tau)\,,\qqq
\Im(i\)\tau)^{1/2}>0\,,\kern-3em
$$
and is given by the following formulae
\vvn-.1>
$$
\gather
f(\la\+1\},\la\+2)\,\)\map\,f(-\tau\)\la\+1\},-\)\tau\)\la\+2)\,,
\\
\nn6>
t_{ab}\;\map\,\exp\)\Bigl(\>{\pii\)\tau\over 2\)N}\,
\bigl(\)(N+1)\tsum_{c=1}^N\)(\la\+1_{ca})\vpb2 -\,
(N-1)\tsum_{c=1}^N\)(\la\+2_{cb})\vpb2 -\,
4\>N\)\la\+1_a\la\+2_b\>\bigr)\]\Bigr)\ t_{ab}\,.
\endgather
$$
\par
Introduce a matrix $\Tc(u)$ with noncommuting entries:
$$
\Tc(u)\,=\>\sum_{a,b}\>E_{ab}\ox\)\Tc_{ab}(u)\,,\kern 4em
\Tc_{ab}(u)\>=\>\tht(u-\la\+1_b\!+\la\+2_a)\,t_{ba}\,.\kern-2em
\Tag{Tcu}
$$
\vsk->
\vsk0>
\Th{Rtt}
The commutation relations \(tbc)\,--\,\(tac) are \eqv/ to
$$
R\"{12}(u-v,\la\+2)\>\Tc\"{13}(u)\>\Tc\"{23}(v)\,=\,
\Tc\"{23}(v)\>\Tc\"{13}(u)\>R\"{12}(u-v,\la\+1)\,,\kern-1em
\Tag{RTTc}
$$
\endpro
\nt
The proof is straightforward and is based on summation formulae for
the theta \fn/. Notice that formula \(RTTc) is similar to \(RTT)\).
\Par
Introduce the \em{quantum determinant} $\Det\Tc(u)$ like in \(detT)\):
$$
\Det\Tc(u)\,=\;{\Tht(\la\+1)\over\Tht(\la\+2)}\,
\sumib\sign(\)\ib\))\;\Tc_{i_N\],\)N}\bigl(u+(N-1)\>\gm\bigr)
\>\ldots\,\Tc_{i_2,\)2}(u+\gm)\,\Tc_{i_1,\)1}(u)\kern-1em
\vv-.2>
\Tag{detTc}
$$
where \,$\dsize\Tht(\la)\>=\!\!\prod_{\abn}\!\tht(\la_a\]-\la_b)$.
By \(Tcu) and \(tf) we obtain
\vvn-.7>
\vvnm.4>
$$
\Det\Tc(u)\,=\;{\Tht(\la\+1)\over\Tht(\la\+2)}\,\sumib\sign(\)\ib\))\,
\pran\>\tht\bigl(u+(a-1)\>\gm-\la\+1_a\!+\la\+2_{i_a}\bigr)\;
t_{N,\)i_N}\>\ldots\,t_{1,\)i_1}\,.\kern-2em
\Tag{dett}
$$
It is clear that \,$\Det\Tc(u)\,f(\la\+1\},\la\+2)\)=
f(\la\+1\},\la\+2)\,\Det\Tc(u)$ \,for any $f\}\in\Funt$.
\vsk.7>
\Prop{qdets}
$\bigl[\>\Det\Tc(u)\>,\)t_{ab}\)\bigr]\,=\,0\,$ for any $a\),\bno$.
\endpro
\nt
The proof is similar to the proof of Proposition~\[qdet]\).
\Par
Let $\fb\)=\)(\fN)$ be a \mccl/ with coefficients in $\FunC$,
\vv.05>
that is, the \fn/s $\fN\in\FunC$ satisfy the condition
$$
f_a(\la)\,f_b(\la-\gm\>\eps_a)\,=\,f_a(\la-\gm\>\eps_b)\,f_b(\la)
$$
for any $a\),\bno$, \cf. Appendix~\[:mform]\). Then the assignments
$$
t_{ab}\>\map\>f_a(\la\+1)\,t_{ab}\qqq\text{and}\qqq
t_{ab}\>\map\>f_b(\la\+2)\,t_{ab}
\Tag{faut}
$$
define two \endo/s of the \oalg/ $\els$.
\vv-.03>
The \endo/s are \aut/s if $\fb$ is \ndeg/, that is, if ${f_a\ne\)0}$
\vv.05>
for any $\ano$. The \aut/s are inner if $\fb$ is a \mcbd/:
$f_a(\la)=\phi(\la-\gm\>\eps_a)\)/\)\phi(\la)$
for a certain \fn/ $\phi\in\FunC$, see \(tf)\).
\Rem
In this paper the inequality ${f\ne\)0}$ for a \mef/ $f$ means that
the \fn/ $f$ is not identically zero.
\vvgood
\enddemo
\Prop{Caut}
The assignment
\vvn-.75>
\vvnm.75>
$$
\gather
f(\la\+1\},\la\+2)\;\map\,f(\>-\>\la\+2\},\)-\>\la\+1)\,,
\Tag{Cart}
\\
\nn7>
t_{ab}\;\map\,\prod_{\tsize{c\)=1\atop c\)\ne b}}^N\>\tht(\la\+1_{cb})\)
\prod_{\tsize{c\)=1\atop c\)\ne a}}^N\)
\bigl(\)\tht(\la\+2_{ca}\!+\gm)\bigr)\vpb{-1}\,t_{ba}
\\
\cnnm.2>
\cnn-.35>
\endgather
$$
defines an involutive anti\aut/ of the \oalg/ $\els$.
\endpro
A \em{module over the small \eqg/} $\esl$ is a \dhmod/ $V\}$ endowed with
an action of the operator algebra $\els$ in the space $\FunV$. Relations
\(la12) and \(tf) mean that $t_{ab}$ acts on $\FunV$ as a \dif/ operator:
$$
(t_{ab}\,v)\)(\la)\,=\,\ell_{ab}(\la)\>v(\la-\gm\>\eps_a)
\qqq\text{for any}\quad v\in\FunV\kern-2em
\Tag{ell}
$$
where $\ell_{ab}(\la)\in\Fend V$ is a suitable \fn/.
In $\Diff(V)$ relations \(la12)\), \(tf) are \eqv/ to
\vvnn0:-.2>
$$
t_{ab}\,\phi(\la)\>=\>\phi(\la-\gm\>\eps_a)\,t_{ab}
\vvmm-.4:-.8>
$$
for any $\phi\in\FunC$, and
\vvm-.4>
$$
t_{ab}\>\FVmu\,\sub\,
\Fun\bigl(V[\)\mu+\eps_a\]-\eps_b\)]\)\bigr)
\vvmm-.1:.2>
$$
for any $\mu\in\hga\}$, which are similar to \(Tph) and \(Tmu)\).
For the quantum determinant we have
$$
\bigl(\)\Det\Tc(u)\,v\bigr)\)(\la)\,=\,\Dc(u\),\la)\>v(\la)
\vvm.2>
$$
for any $v\in\FunV$ where $\Dc(u\),\la)\in\Fend V$ is a suitable \fn/.
\Ex
The assignment
\vvnn-1.1:-1.5>
$$
\alignat2
\ell_{aa}(\la)\, &{}\map\,E_{aa}\>+\sum_{\tsize{b=1\atop b\ne a}}^N\>
{\tht(\la_{ab}\]+\gm)\over\tht(\la_{ab})}\;E_{bb}\,, &&
\Tag{vect}
\\
\nn7>
\ell_{ab}(\la)\, &{}\map\;{\tht(\gm)\over\tht(\la_{ab})}\;E_{ab}\,, &&
\Llap{a\ne b\,,}
\\
\cnn-.15>
\endalignat
$$
$a\),\bno$, makes $\CN\!$ into an \emod/.
\vv.1>
The module is called the \em{vector \rep/} of $\esl$. In the vector \rep/
$$
\Dc(u\),\la)\,=\,\tht(u-\gm)\,\tht(u+\gm)\,\tht(u+2\gm)\>\ldots\>
\tht\bigl(u+(N-1)\>\gm\bigr)\,.
$$
\enddemo
\Cr{ehom}
For any \emod/ $V\}$ and any $x\in\C$ the rule
\,$T_{ab}(u)\>=\)\Tc_{ab}(u-x)\vst{V}$ makes $V\}$ into an \Emod/ called
the \em{\evmod/} $V(x)$ over $\Eqg$.
\endpro
By abuse of notation we call the assignment \)${T(u)\>\map\)\Tc(u)}$
\)the \em{\eval/ morphism} $\Eqg\)\to\>\esl$. It is analogous to the \ehom/
from the Yangian $\Ysl$ to $\Usl$.
\Rem
Corollary~\[ehom] was our main motivation to discover and study the small
\eqg/ $\esl$.
\enddemo
Let $V,\)W\}$ be \emod/s. An element ${\phi\in\Fun\bigl(\Hom_\hg(V,W)\bigr)}$
is a \em{morphism} of \emod/s if the induced map intertwines the corresponding
actions of $\els$:
$$
\phi(\la)\,t_{ab}\vst{\italic{Fun}(V)}=\,t_{ab}\vst{\italic{Fun}(W)}\>\phi(\la)
\vvm.2>
$$
for any $a\),\bno$. Denote by $\Mor(V,W)$ the space of all morphisms from $V\}$
to $W\}$. A morphism $\phi$ is called an \em{\iso/} if the map $\phi(\la)$ is
bijective for generic $\la$.
\goodbreak
\vskm0:.2>
An \emod/ $V\}$ is called \em{\irr/} if for any nontrivial morphism
$\phi\in\Mor(W,V)$ the map $\phi(\la)$ is surjective for generic $\la$,
and \em{reducible} otherwise.
\vsk.2>
Say that an \emod/ $W\}$ is a \em{submodule} of $V\}$ if there is a morphism
$\phi\in\Mor(W,V)$ \st/ the map $\phi(\la)$ is injective for generic $\la$.
The morphism $\phi$ is called an \em{embedding}. The submodule $W\}$ is called
\mmgood
\em{proper} if $\phi$ is not an \iso/. Any \emod/ $V\}$ has at least two
submodules: $V\}$ itself and the \em{trivial submodule} $\lb\)0\)\rb$
with obvious embeddings.
\vsk.2>
{\bls1.06\bls
Let $W\}$ be a submodule of $V\}$. Then one can define the \em{quotient} \emod/
$V\}/\)W\}$ as follows. Fix an embedding $\phi$. The subspace
\vv.04>
$\phi\bigl(\FunW\bigr)\sub\FunV$ is \inv/ \wrt/ the action of $\els$, hence
$\els$ acts on $\FunV\big/\phi\bigl(\FunW\bigr)$. Let ${\la_0\in\hga}\}$
be \st/ the map $\phi(\la_0)$ is injective. Take a complement $U\}$ of
$\phi(\la_0)\>W\}$ in $V\}$, that is, $V\}=U\oplus\phi(\la_0)\>W\}$ as a vector
space. Notice that $V\}=U\oplus\phi(\la)\>W\}$ for generic $\la$ as well.
Then $\FunV=\FunU\)\oplus\phi\bigl(\FunW\bigr)$ and, therefore,
$\FunU=\FunV\big/\phi\bigl(\FunW\bigr)$, which induces an action of $\els$ on
$\FunU$ and makes $U\]$ into an \emod/. The constructed \emod/ does not depend
on a choice $\phi$, $\la_0$ and $U\]$ up to an \iso/ of \emod/s and is denoted
by $V\}/\)W\}$.
\par}
\Lm{submod}
An \emod/ $V\}$ is reducible if and only if it has a nontrivial proper
submodule.
\endpro
\Lm{inj}
An \emod/ $V\}$ is \irr/ if and only if for any nontrivial morphism
$\phi\in\Mor(V,W)$ the map $\phi(\la)$ is injective for generic $\la$.
\endpro

\Sect[H]{Highest \wt/ modules over $\{\esl$}
For any monomial $\tabk$ set \,${\deg\>(\)\tabk)=\)k}$.
\vv.1>
For ${k=0}$ we assume that the monomial equals $1$. As a \$\Funt$-module
the \oalg/ $\els$ is \gby/ all monomials $\tabk$, $k=0\),1\),\ldots{}$.
\vsk.2>
For any \fn/ ${\phi\in\Funt}$ set $\deg\>(\phi)=0$. Since relations
\vv.1>
\(tf)\,--\,\(tac) are homogeneous, the algebra $\els$ is \$\Zp\>$-graded by
$\deg$. Let \)$\eg_k=\lb\)x\in\els\vert \deg\>(x)=k\)\rb$ be the homogeneous
subspace of degree $k$. Each subspace \)$\eg_k\]$ is finitely generated over
$\Funt$.
\vsk.3>
Consider the \em{normal ordering} of generators: $t_{ab}\]<t_{cd}$
if $a-b<c-d$, or $a-b=c-d$ and $a<c$. Say that the monomial $\tabk$
is \em{\nord/} if $t_{a_ib_i}\!<t_{a_jb_j}$ for any $i<j$, or $k=0$.
\Th{PBW}
For any $k\in\Zp$ the \norm/s of degree $k$ form a basis of \,$\eg_k\]$
over $\Funt$.
\endpro
\Pf.
For $k=0$ and $k=1$ the claim is immediate. Let $k>1$. Here we prove that the
\norm/s of degree $k$ span \)$\eg_k$ over $\Funt$. The linear independence of
the \norm/s is proved in Appendix~\[:norm]\).
\vsk.2>
It is clear from relations \(tbc)\,--\,\(tac) that any product $t_{ab}\>t_{cd}$
can be written as a linear combination of \nord/ products. Given a monomial
$\tabk$ we take any disordered product of adjacent factors and replace it by
a suitable sum of \nord/ products, then do the same for each of the obtained
monomials. To see that the procedure terminates and, hence, produces a linear
combination of \norm/s, introduce auxilary gradings on monomials by the rule
\vvn-.3>
$$
r\)(\tabk)\>=\>\tsum_{i=1}^k\>i\>(a_i\]-b_i)\,,\kern4em
r'(\tabk)\>=\>\tsum_{i=1}^k\>i\>b_i\,,
\vv-.2>
$$
and observe that at each nontrivial step of the procedure we replace a monomial
by a sum of monomials of either less degree $r$, or the same degree $r$ and
less degree $r'\]$.
\epf
Introduce modified generators of the algebra $\els$. For any $a\),\bno$ set
\vvnn-.8:-.2>
$$
\tth_{ab}\,=\prod_{1\le c<a}\}\tht(\la\+1_{ca})
\prod_{1\le c<b}\}\bigl(\)\tht(\la\+2_{cb})\bigr)\vpb{-1}\,t_{ab}\,.\kern-1em
\vv-.5>
\Tag{tth}
$$
\par
Let $V\}$ be an \emod/. A nonzero \fn/ ${v\in\FunV}$ is called a \em{\singv/}
if ${t_{ab}\,v\)=\)0}$ for any ${\abn}$. Say that $v$ is a \em{\regsv/} if,
in addition, $v$ is a \wt/ vector \wrt/ the action of $\hg$ and
\vvn-.3>
$$
(\tth_{aa}\)v)\)(\la)\,=\,Q_a(\la)\,v(\la)\,,\qqq\ano\,,\kern-2em
\vvgood
\vvm-.2>
\Tag{tQ}
$$
\vv.1>
for certain \fn/s $\QN\in\FunC$. We call $\Qh=(\)\QN)$
the \em{\dwt/} of $v$. Relation \(thab) implies that $\Qh$ is a \mccl/:
$$
Q_a(\la)\,Q_b(\la-\gm\>\eps_a)\,=\,Q_a(\la-\gm\>\eps_b)\,Q_b(\la)
$$
for any $a\),\bno$. If ${f\}\in\FunC}$, then the \fn/
$\vti(\la)\)=f(\la)\,v(\la)$ is a \regsv/ of \dwt/ $(\Qti_1\lc\Qti_N)$ where
\vvn-.4>
$$
\Qti_a(\la)\,=\,Q_a(\la)\;{f(\la-\gm\>\eps_a)\over f(\la)}\;,\qqq\ano\,.
\kern-2em
\vvmm-.4:-.2>
$$
Hence, the subspace $\FunC\,v$ determines the \dwt/ up to a \mcbd/.
\vsk.2>
Say that $\Qh$ and $v$ are \em{\ndeg/} if ${Q_a\ne\)0}$ for any $\ano$.
\vv.06>
Say that $\Qh$ and $v$ are \em{standard} \>if ${Q_a=1}$ for any $\ano$.
\vsk.3>
By formula \(dett) the quantum determinant acts on a \regsv/ $v$ of \wt/ $\mu$
and \dwt/ $\Qh$ as follows:
\vvnn-.2:-.5>
$$
\bigl(\)\Det\Tc(u)\,v\bigr)\)(\la)\,=\,
\pran\>\tht\bigl(u-\gm\>(\)\mu_a\}-a+1)\bigr)\>
\pran\>Q_a\bigl(\la-\!\tsum_{a<b\le N}\!\gm\>\eps_b\bigr)\;v(\la)\,.\kern-1em
\vvm-.2>
\Tag{Detv}
$$
\par
{\bls1.1\bls
An \emod/ $V\}$ is called a \em{\hwm/} with \em{\hw/} $\mu$, \em{\dhw/} $\Qh$
and \em{\hwv/} $v$ if $v$ is a \regsv/ of \wt/ $\mu$ and \dwt/ $\Qh$ generating
$\FunV$ over $\els$. If $\Qh$ is standard (\ndeg/), then $V\}$ is called
a \em{standard} (\em{\ndeg/}) \emod/ of \hw/ $\mu$. For example, the vector
\rep/ is a standard \emod/ of \hw/ $\om_1$.
\vsk.17>}
It is clear that any \ndeg/ \hwm/ is isomorphic to a pullback of a standard
\hwm/ of the same \hw/ through a suitable \aut/ of the form \(faut)\).
\Prop{hwred}
Let $V\}$ be a \hwm/ with \hw/ $\mu$. Then
\vsk.2>
\atem
$V=\Plus_{\nu\le\mu}\]\Vnu$ \,and \,$\dimC\Vmu=1${\rm;}
\vsk.1>
\bitem
$V\}$ is reducible if and only if it has a \singv/ of \wt/ $\nu<\mu$.
\mmgood
\endpro
\Lm{hwDet}
Let $V\}$ be a \hwm/ of \hw/ $\mu$ and \dhw/ $(\QN)$.
Then for any $v\in\FunV$
\vvn-.2>
$$
\bigl(\)\Det\Tc(u)\,v\bigr)\)(\la)\,=\,
\pran\>\tht\bigl(u-\gm\>(\)\mu_a\}-a+1)\bigr)\>
\pran\>Q_a\bigl(\la-\!\tsum_{a<b\le N}\!\gm\>\eps_b\bigr)\;v(\la)\,.
\vv-.2>
\Tag{Dethw}
$$
\endpro
\nt
The statement follows from Proposition~\[qdets] and formula \(Detv)\).
\Cr{hwirr}
Let ${\gm\)\nin\)\QtQ}$. Then a \ndeg/ \hwm/ $V\}$ of \hw/ $\mu$ is \irr/
\vv.1>
unless $V[\)w\cddt\mu\)]\)\ne\)0$ \,for some $w\in W\}$ \st/ $w\cddt\mu<\mu$.
Here $W\}$ is the Weyl group.
\endpro
\Pf.
\bls 1.06\bls
If $V\}$ is reducible, then it has a \singv/ $v$ of \wt/ $\nu<\mu$. Comparing
formulae \(Detv) and \(Dethw) for the action of $\Det\Tc(u)$ on $v$ we obtain
that $(\)\nu_1-1\)\lc\nu_N-N)=(\)\mu_{i_1}\!-i_1\)\lc\mu_{i_N}\!-i_N)$
for some \perm/ $(i_1\lc i_N)$, since ${\gm\)\nin\)\QtQ}$. That is,
$\nu=w\)(\)\mu+\rho)-\rho=w\cddt\mu$ for some $w\in W\}$ because
$\rho\>=\>-\,\eps_1-2\>\eps_2\lsym-N\eps_N$.
\epf
Let $e(\bgp)$ and $e(\ngp)$ be the left ideals in $\els$ \gby/
\vv.1>
the elements $t_{ab}$ with ${a\le b}$ and ${a<b}$, \resp/.
Let $\Bc\,,\,\Bc'\],\,\Nc,\,\Nc'$ be the following sets of \norm/s
\vvnn-.1:.1>
$$
\align
\Bc\, &{}=\>\lb\>\thabk\vert k\ge 0\,,\ \ a_i\ge b_i\,,\ \;i=1\lc k\)\rb\,,
\Tag{BNc}
\\
\nn6>
\Bc'\]\, &{}=\>\lb\>\thabk\vert k> 0\,,\ \ a_i\le b_i\,,\ \;i=1\lc k\)\rb\,,
\\
\nn6>
\Nc\> &{}=\>\lb\>\thabk\vert k\ge 0\,,\ \ a_i>b_i\,,\ \;i=1\lc k\)\rb\,,
\\
\nn6>
\Nc'\] &{}=\>\lb\>\thabk\vert k>0\,,\ \ a_i<b_i\,,\ \;i=1\lc k\)\rb\,.
\\
\nngood
\endalign
$$
\Lm{idb}
The \norm/s of the form \,${m\>m'}$ where \,$m\in\Nc$ and
\vv.1>
\,$m'\!\in\Bc'\]$ form a basis of \,$e(\bgp)$ over $\Funt$.
\endpro
\Lm{idn}
The \norm/s of the form \,${m\>m'}$ where \,$m\in\Bc$ and
\vv.1>
\,$m'\!\in\Nc'\]$ form a basis of \,$e(\ngp)$ over $\Funt$.
\endpro
\nt
The statements easily follows from relations \(tbc)\,--\,\(tac) and
Theorem~\[PBW]\).
\vsk.2>
For any monomial $\tabk$ set
\vv-.2>
\>${\wsl\)(\tabk)\,=\)\tsum_{i=1}^k\>(\)\eps_{a_i}\!-\)\eps_{b_i})}$,
\vv.1>
and for any \fn/ $\phi\in\Funt$ set $\wsl\)(\phi)=0$.
Since relations \(tf)\,--\,\(tac) are homogeneous,
the algebra $\els$ is \$\PP$-graded by $\wsl$.
\vsk.3>
Let $\mu\in\hga\}$ and $\Qh$ be a \mccl/.
\vv-.15>
Below we define a \em{\Vmod/} $\MmQ$ of \hw/ $\mu$ and \dhw/ $\Qh$ over $\esl$.
\vsk.4>
{\bls1.16\bls Let \>${\NC\)=\]\!\Plus_{m\>\in\>\Nc}\]\C\>m}$ be a \dhmod/ \st/
a monomial $m$ has \wt/ $\wsl\)(m)$, and let $\C\>\vmQ$ be a \onedim/ \hmod/
\st/ $\vmQ$ has \wt/ $\mu$.
\vv.06>
Then $\MmQ=\NC\]\ox\C\>\vmQ$ as an \hmod/. We define an action of
$\els$ in $\Fun(\MmQ)$ by the rule: $1\ox\vmQ$ is a \regsv/ of
weight $\mu$ and \dwt/ $\Qh$, and
\vvn-.3>
$$
m\,(\)1\ox\vmQ\))\,=\,m\ox\vmQ
\vv-.15>
$$
for any $m\in\Nc\]$. This determines an action on $\vmQ$ by any \norm/ and,
hence, by any element of $\els$, \cf. Theorem~\[PBW]\). Finally, for any
$x\in\els$ set
\vvn-.2>
$$
x\,(\)m\ox\vmQ\))\,=\,(x\)m)\,(\)1\ox\vmQ\))
\vvmm.1:-.15>
$$
where the product $x\)m$ should be represented as a linear combination
of \norm/s.
\par}
\Prop{Mmu}
$\MmQ$ is a well-defined \emod/ with \hw/ $\mu$, \dhw/ $\Qh$ and \hwv/
$1\ox\vmQ$.
\endpro
\nt
The statement follows from Lemmas~\[idb]\), \[idn] \)and Theorem~\[PBW]\).
\Par
{}From now on we suppress the symbol of tensor product in the definition
of the \Vmod/ $\MmQ$, for instance, we write $\vmQ$ instead of $1\ox\vmQ$.
\vsk.3>
If $\Qh$ is standard, then the \emod/ ${M_\mu=\MmQ}$ is called
the \em{standard} \Vmod/ with \hw/ $\mu$ and the \hwv/ $v_\mu=\vmQ$.
\Lm{MVv}
Let $V\}$ be an \emod/. For any $v\in\FunV$ which is a \regsv/ of \wt/ $\mu$
and \dwt/ $\Qh$, there is a unique morphism $\phi\in\Mor(\MmQ\,,V)$ which sends
$\vmQ$ to $v$.
\endpro
\Prop{hwqu}
Any \hw/ \emod/ is isomorphic to a suitable quotient of the \Vmod/ over $\esl$
of the same \hw/ and \dhw/.
\endpro
\Prop{irrmu}
For any ${\mu\in\hga}\}$ and a \dwt/ $\Qh$ there exists a unique
up to an \iso/ \irr/ \hw/ \emod/ with \hw/ $\mu$ and \dhw/ $\Qh$.
\endpro
\Prop{irrV}
Let ${\gm\)\nin\)\QtQ}$. Then a \hw/ \emod/ with \hw/ $\mu$ and a \ndeg/ \dhw/
$\Qh$ is isomorphic to the Verma module $\MmQ$ unless $w\cddt\mu<\mu$ for some
$w\in W\}$.
\endpro

\Sect[S]{Dynamical \Sform/}
Let $e(\ngm)$ be the right ideal in $\els$ \gby/ the elements $t_{ab}$ with
${a>b}$, let $\dg$ be the \$\Funt$-submodule \gby/ \norm/s of the form $\taak$,
and let
\vvnm-.3>
$$
\egz\,=\,\lb\)x\in\els\)\vert{\wsl\)(x)=0}\)\rb
\vvgood
\vv-.2>
$$
be the subalgebra of zero \wt/ elements in $\els$. Consider the quotient
\vvn-.1>
$$
\ehg\,=\,\els\)\big/\bigl(\)e(\ngp)+e(\ngm)\bigr)\,.
\vv-.1>
$$
Let $\eta\):\)\els\>\to\>\ehg$ be the natural projection. By Theorem~\[PBW]
\vv.1>
the restriction of $\eta$ to $\dg$ is a bijection. Denote by
\vv.07>
$\etb\):\)\ehg\>\to\>\dg$ the inverse map. For any $x\),y\in\ehg$ define their
product by the rule $x\)y=\eta\)\bigl(\etb\)(x)\>\etb\)(y)\bigr)$. It is easy
to see that this defines an algebra structure on $\ehg$.
\vskmgood-.7:.7>
\Lm{HC}
The restriction of $\eta$ to \)$\egz$ is a \hom/.
\endpro
Set $\qh_a=\eta\)(\)\tth_{aa})$, $\ano$. It follows from \(thab) that $\ehg$ is
\gby/ \fn/s $f\}\in\Funt$ and the pairwise commuting elements $\qh_1\lc\qh_N$
subject to relations
\vvm-.6>
$$
\qh_a\,f(\la\+1\},\la\+2)\,=\>
f(\la\+1\!-\gm\>\eps_a\>,\la\+2\!-\gm\>\eps_a)\,\qh_a\,.\kern-1.4em
\Tag{qf}
$$
\vsk.2>
The assignment
\hfill $\pim\>:\>f(\la\+1\},\la\+2)\;\map\,f(\>-\>\la\+2\},\)-\>\la\+1)$,
\hp{The assignment\kern\parindent\kern-.9em}\nl
\vv.3>
$$
\pim\>:\,t_{ab}\;\map\,
\prod_{1\le c<b}\}\tht(\la\+1_{cb})\,\tht(\la\+1_{bc}\!-\gm)
\prod_{1\le c<a}\}\bigl(\)\tht(\la\+2_{ca})\,\tht(\la\+2_{ac}\!-\gm)\bigr)
\vpb{-1}\,t_{ba}
$$
defines an involutive anti\aut/ of the \oalg/ $\els$, which differs from
the anti\aut/ \(Cart) by a suitable \aut/ of the form \(faut). We have
\vvm-.6>
$$
\pim(\tth_{ab})\,=\,\)\tth_{ba}\,.
\Tag{pimth}
$$
For any $x\in\ehg$ set $\pim(x)=\eta\)\bigl(\pim(\)\etb\)(x))\bigr)$.
\Lm{etapim}
For any $m\in\els$ we have
$\eta\)\bigl(\pim(m)\bigr)\>=\>\pim\bigl(\)\eta\)(m)\bigr)$.
\endpro
For any $m_1\),m_2\in\els$ set
\vvnn0:-.4>
$$
S(m_1\>,\)m_2)\,=\,\eta\)\bigl(\pim(m_1)\>m_2\bigr)\,\in\,\ehg\,.
\vvmm0:-.2>
\Tag{Smm}
$$
Since $\pim$ is involutive,
\>$S(m_1\>,\)m_2)\>=\>\pim\bigl(S(m_2\>,\)m_1)\bigr)$.
Moreover, \>$S(m_1\)m_2\>,\)m_3)\>=\>S\bigl(m_2\>,\)\pim(m_1)\>m_3\bigr)$
for any $m_1\),m_2\),m_3\in\els$. $S$ is called the \em{\dSform/} on $\els$.
\vskm-.5:.1>\vsk0>
\Ex
Let $a<b$. Then \,$\dsize S(\)\tth_{ba}\>,\)\tth_{ba})\,=\,
-\,{\tht(\la\+1_{ab}\!-\la\+2_{ab})\,\tht(\gm)\over
\tht(\la\+1_{ab})\,\tht(\la\+2_{ab})}\ \qh_a\>\qh_b$.
\enddemo
Let $\mu\in\hga\}$ and let $\Qh\)=(\QN)$ be a \mccl/.
\vv.15>
Then there is an algebra \hom/ $\chimQ\>:\)\ehg\>\to\DiffC$:
\vvnn-.2:-.8>
$$
\gather
\chimQ\)\bigl(\)f(\la\+1\},\la\+2))\,:\,\phi(\la)\,\map\,
f(\la\>,\la-\gm\>\mu)\,\phi(\la)\,,
\\
\nn5>
\chimQ\)(\)\qh_a)\,:\,\phi(\la)\,\map\,Q_a(\la)\,\phi(\la-\gm\>\eps_a)\,.
\endgather
$$
\vskm0:-.2>
Consider the \Vmod/ $\MmQ$.
\vvn-.1>
For any $m\in\els$ and $v\in\Fun(\MmQ)$ set
$$
\SmQ\)(m\>,v)\,=\,\chimQ\)\bigl(S(m\>,m'\))\bigr)\]\cddt\]1\in\FunC\,.
\Tag{SmQv}
$$
Here $m'\}\in\els$ is determined by $v=m'\)\vmQ$ and $\)1\in\FunC$ is
\vv.1>
the constant \fn/. It is easy to see that $\SmQ\)(m\>,v)$ does not depend
\vv.1>
on the choice of $m'\}$. We call $\SmQ$ the \em{\dSpair/} for $\MmQ$.
\goodbm
\Lm{SmQ}
For any \>$m_1\),m_2\in\els$ and \>$v\in\Fun(\MmQ)$ we have
\vvn-.1>
$$
\SmQ\)(m_1\>,m_2\>v)\,=\>\SmQ\)\bigl(\pim(m_2)\>m_1\>,v\bigr)\,=\>
\SmQ\)\bigl(\)1\>,\pim(m_1)\>m_2\>v\bigr)\,.
\vv-.1>
$$
Here $1\in\els$ is the identity element.
\wwgood-:>
\endpro
\Lm{Smv}
Let $m\in\els$ be a \$\wsl\)$-homogeneous element, and $v\in\MmQ[\)\nu\)]$.
\vv.1>
If \)$\wsl\)(m)\ne\mu-\nu$, then \)$\SmQ\)(m\>,v)\)=\)0$.
Otherwise, \>$\pim(m)\,v\)=\)\SmQ\)(m\>,v)\,\vmQ$.
\endpro
\Ex
$\SmQ\)(\)1\>,\vmQ)\)=\)1$. \;If $a<b$, then
\vvn-.3>
$$
\SmQ\)(\)\tth_{ba}\>,\)\tth_{ba}\)\vmQ)\,=\,
-\,{\tht(\gm\>\mu_{ab})\,\tht(\gm)\over\tht(\la_{ab})\,
\tht(\la_{ab}\]-\gm\>\mu_{ab})}\ Q_a(\la)\>Q_b(\la-\gm\>\eps_a)\,.
\vvmm-.6:-.1>
$$
\enddemo
Set
\vvnn-.4:-.5>
$$
\Ker\SmQ\,=\>\lb\)v\in\Fun(\MmQ)\)\vert\SmQ\)(m\>,v)\)=\)0\quad
\text{for any}\ \ \,m\in\els\)\rb\,.
$$
The subspace $\Ker\SmQ$ is \inv/ under the action of $\els$ and defines
a proper submodule $\NmQ$ of $\MmQ$.
\Prop{NmQ}
$\NmQ$ is the maximal proper submodule of $\MmQ$, that is, for any proper
\vv.1>
submodule $U\]$ of $\MmQ$ the embedded image of $\FunU$ is contained in
$\Fun(\NmQ)=\Ker\SmQ$.
\endpro
\Pf.
\bls1.1\bls
Let $\phi\in\Mor(U,\MmQ)$ be the embedding. $\phi\bigl(\FunU\bigr)$ is
a direct sum of its \wt/ components and $\phi\bigl(\FunU\bigr)[\)\mu\)]\)=\)0$,
since $U\]$ is a proper submodule.
Therefore, if $v\in\phi\bigl(\FunU\bigr)[\)\nu\)]$ and $\wsl\)(m)=\mu-\nu$,
then $\pim(m)\>v=0$. Hence, by Lemma~\[Smv] $\SmQ\)(m\>,v)=0$ for any $m\in\els$
and $v\in\phi\bigl(\FunU\bigr)$.
\epf
\Cr{VmQ}
The quotient module ${\VmQ\)=\MmQ\)/\NmQ}$ is
\vv.05>
the \irr/ \hw/ \emod/ with \hw/ $\mu$ and \dhw/ $\Qh$.
\endpro
For any \hwm/ $V\}$ with \hw/ $\mu$, \dhw/ $\Qh$ and \hwv/ $v$ one can define
the \Spair/ similarly to \(SmQv): for any ${m\in\els}$ and $v'\}\in\FunV$ set
\vvn-.1>
$$
\SmQV\)(m\>,v')\,=\,\chimQ\)\bigl(S(m\>,m'\))\bigr)\]\cdot\]1\in\FunC
\vvmm0:.2>
$$
where $m'\}\in\els$ is determined by $v'\]=m'\)v$. Propositions~\[hwqu] and
\[NmQ] \)imply that $\SmQV\)(m\>,v)$ does not depend on the choice of $m'\}$.
\Prop{Virr}
The module $V\}$ is \irr/ if and only if \>$\Ker\SmQV$ is trivial.
Otherwise, $\Ker\SmQV$ defines the maximal proper submodule of $V\}$.
\endpro
\vskm-.1:0>
Let $\Vslm$ be the \irr/ \smod/ of \hw/ $\mu$. Set $\dmn=\dimC\Vsmn$.
\Th{fdim}
Let ${\mu\in\hga}\}$ be a \dint/. Then \>$\dimC\VmQ[\)\nu\)]\>\le\)\dmn$
for any $\nu\in\hga\}$. In particular, the \emod/ $\VmQ$ is \fd/.
\mmgood
\endpro
\Th{dimu}
Let ${\gm\)\nin\)\QtQ}$. The \emod/ $\VmQ$ for
\vv-.14>
a \ndeg/ \dwt/ $\Qh$ is \fd/ if and only if $\mu$ is a \dint/.
\vv.1>
Moreover, \>$\dimC\VmQ[\)\nu\)]\>=\)\dmn$ for any $\nu\in\hga\}$.
\endpro
\nt
Theorems~\[fdim] and \[dimu] are proved in Section~\[:Fd]\).
\vsk.4>
If \,$\Qh$ \)is standard,
we set \>$\chi_\mu=\chimQ$, $S_\mu=\SmQ$, $N_\mu=\NmQ$ \,and \>$V_\mu=\VmQ$.

\Sect[cg]{Contragradient modules over $\{\esl$ and \cform/}
For any \dhmod/ $V\}$ define an involutive linear map $\psi:\FunV\to\FunV$
by the rule: if $f\in\FVmu$, then $(\psi f\))(\la)\)=\)f(-\)\la+\gm\)\mu)$.
\vsk.3>
Let $V\}$ be a \dhmod/ and let ${V^*\]=\!\}\Plus_{\,\mu\>\in\>\hga\!}\Vmua}\!$
be its restricted dual space. We consider $V^*\!$ as a \dhmod/ \st/ $\Vmua\!$
\vv.05>
is a \wt/ subspace of \wt/ $\mu$. For any $B\in\End(V)$ we denote by
\vv.05>
$B^*\}\in\End(V^*)$ the dual map. For a \dif/ operator $A\in\Diff(V)$
we define the dual operator $A'\}\in\Diff(V^*)$ by the rule:
\vvm-.3>
$$
\text{if}\quad (A\>v)(\la)=B(\la)\>v(\la+\mu)\,,\quad\text{then}\quad
(A'\phi\))(\la)\,=\,\bigl(B(\la-\mu)\bigr)\vpa\phi(\la-\mu)\,,
$$
and the operator ${A^\dag\}\in\Diff(V^*)}$: ${A^\dag\]=\)\psi\>A'\psi}$.
\vv.06>
The assignment ${A\>\map\]A^\dag}\]$ is an involutive anti\iso/
$\Diff(V)\alb\to\)\Diff(V^*)$.
\vskgood-:>
\Ex
Let $V\}$ be an \emod/. Then
$$
\bigl((t_{ab})^{\]\dag}\phi\)\bigr)(\la)\,=\,
\bigl(\ell_{ab}\bigl(-\la+(\mu+\eps_b)\>\gm\bigr)\bigr)\vpa\)
\phi(\la-\gm\>\eps_b)
$$
for any $\phi\in\FVamu$, \cf. \(ell)\).
\enddemo
Given an \emod/ $V\}$ we make the \hmod/ $V^*\!$ into an \emod/ as follows:
the action of an element $m\in\els$ in $\FunVa$ is given by
\vv.06>
$\bigl(\pim(m)\bigr)^{\]\dagg}\!$ where $\pim(m)$ is understood as a \dif/
operator acting in $\FunV$.
\vv.06>
It is easy to check that the definition is consistent, that is,
the ring $\Funt$ acts on $\FunVa$ in the prescribed way, \cf. \(la12).
The obtained \emod/ $V^*\!$ is called the \em{\cmod/} to the \emod/ $V\}$.
\vsk.3>
Consider the \Vmod/ $\MmQ$. Recall that $\MmQ[\)\mu\)]=\C\>\vmQ$.
Fix $\vmQa\in\MmQa[\)\mu\)]$ by the rule $\bra\)\vmQa\,,\)\vmQ\)\ket=1$.
For any \)$\ano$ set
\vvn-.2>
$$
\Qti_a(\la)\,=\,Q_a\bigl(-\la+(\mu+\eps_a)\>\gm\bigr)\,.
\vv-.2>
$$
Notice that \>$Q_a(\la)=\Qti_a\bigl(-\la+(\mu+\eps_a)\>\gm\bigr)$ \>as well.
\Prop{vmQa}
For the \cVmod/ $\MmQa$ the constant \fn/ $\vmQa$
\vv-.1>
is a \regsv/ of \wt/ $\mu$ and \dwt/ $\Qti=(\Qti_1\lc\Qti_N)$.
\endpro
\Cr{MM}
There is a unique morphism \>$\pimQ\in\Mor(\MmQ\,,\MmQta)$ sending $\vmQ$
to \)$\vmQta$.
\wwmgood-.7:.7>
\endpro
\Th{kerpi}
$\Ker\pimQ\,=\>\Ker\SmQ$.
\endpro
\vskm.2:.1>
The morphism $\pimQ$ induces a \$\FunC$-bilinear map
\ifMag
\vvnm-.4>
$$
\gather
\BmQ\):\)\Fun(\MmQ)\oxC\Fun(\MmQt)\,\to\>\FunC\,,
\\
\nn8>
\BmQ\)(v\>,\vti\))\,=\,\bra\)\pimQ\>v\>,\vti\ket\,.
\endgather
$$
\else
$\BmQ:\Fun(\MmQ)\oxC\Fun(\MmQt)\>\to\)\FunC$:
$$
\BmQ\)(v\>,\vti\))\,=\,\bra\)\pimQ\>v\>,\vti\ket\,.
$$
\fi
Define a bilinear map $\CmQ:\Fun(\MmQ)\oxC\Fun(\MmQt)\>\to\)\FunC$ by the rule
$$
\CmQ\)(v\>,\vti\))\>=\>\BmQ\)(v\>,\psi\)\vti\))\,.
\Tag{CB}
$$
The map $\CmQ$ is called the \em{\cform/}.
\Th{Cvv}
Let $v\in\Fun\bigl(\MmQ[\)\nu\)]\)\bigr)$ and
$\vti\in\Fun\bigl(\MmQt[\)\nut\)]\)\bigr)$. Then \>$\CmQ\)(v\>,\vti\))\)=\)0$
unless $\nu=\nut$. Moreover,
\vvn.2>
$$
\CmQ\)(v\>,\vti\))(\la)\,=\,\CmQt\)(\vti\>,v\))(-\)\la+\gm\)\nu\))\,.
$$
\endpro
\vsk-.6>
\vsk0>
\Th{kerC}
\vsk.1>
\atem
$\CmQ\)(v\>,\vti)\)=\)0$ \>for any \>$\vti\in\Fun(\MmQt)$
\>if and only if \>$v\in\Ker\SmQ$.
\vsk.06>
\bitem
$\CmQ\)(v\>,\vti)\)=\)0$ \>for any \>$v\in\Fun(\MmQ)$
\>if and only if \>$\vti\in\Ker\SmQt$.
\endpro
\nt
Theorems~\[kerpi]\,\)--\,\)\[kerC] are proved at the end of the section.
\Par
Let $V\}$ and $\Vti\}$ be \hw/ \emod/s of the same \hw/ $\mu$ and \dhw/s $\Qh$
\vvmm.08:-.06>
and $\Qti$, \resp/. By the last corollary the form $\CmQ$ descends
to a form ${\FunV\oxC\Fun(\Vti)\>\to\)\FunC}$,
\vv.1>
denoted by the same letter.
After obvious modification Theorem~\[kerC] remains true in this case.
\vv.1>
In particular, for the \irr/ \hw/ \emod/s $\VmQ$ and $\VmQt$
\vv.06>
the corresponding
\fn/-valued bilinear form $\Fun(\VmQ)\oxC\Fun(\VmQt)\>\to\)\FunC$ is \ndeg/.
\Ex
$\CmQ\)(\vmQ\,,\vmQt\))\)=\)1$. \;If $a<b$, then
\vvn-.3>
$$
\CmQ\)(\)\tth_{ba}\)\vmQ\,,\)\tth_{ba}\)\vmQt)\,=\,
-\,{\tht(\gm\>\mu_{ab})\,\tht(\gm)\over\tht(\la_{ab}\]+\gm)\,
\tht(\la_{ab}\]-\gm\>\mu_{ab}\]+\gm)}\ Q_a(\la+\gm\>\eps_a)\>Q_b(\la)\,.
\vv-.1>
$$
\enddemo
For a monomial \>$m=\thabk$ set \>$\zt'(m)=\tsum_{i=1}^k\>\eps_{a_i}$ and
\vv-.2>
\>$\zt''(m)=\tsum_{i=1}^k\>\eps_{b_i}$. Notice that $\wsl(m)=\zt'(m)-\zt''(m)$.
\vskgood-.7:.7>
\Lm{SQQ}
For any $m_1\),m_2\in\els$ we have
$$
\SmQ\)(m_1\>,m_2\>\vmQ)(\la)\,=\>\SmQt\)(m_2\>,m_1\)\vmQt)
\bigl(-\la+\gm\>\bigl(\)\mu+\zt'(m_1)+\zt''(m_2)\bigr)\]\bigr)\,.
\vv-.2>
$$
\endpro
\Pf.
\Wlg/ we can assume that $\zt'(m_1)+\zt''(m_2)=\zt''(m_1)+\zt'(m_2)$
because otherwise the expressions on both sides of the formula equal zero
by Lemma~\[Smv]\).
\vsk.2>
Let $m_1=\thabk$ and $m_2=\thcdl$.
\vvm.1>
Set $(s_1\lc s_{k+l})=(a_1\lc a_k,d_1\lc d_l)$.
By the definition of the \Spair/, \cf. \(SmQv)\), for any \mccl/ $\Qh$ we have
\vvnm-.8>
$$
\SmQ\)(m_1\>,m_2\>\vmQ)(\la)\,=\>S_\mu(m_1\>,m_2\>v_\mu)(\la)\;
\prod_{i=1}^{k+l}\,Q_{s_i}\bigl(\la-\!\tsum_{1\le j<i}\!\gm\>\eps_{s_j}\bigr)
\,.
\vvm-.4>
\Tag{SQS}
$$
Here $S_\mu$ and $v_\mu$ correspond to the so-called standard case,
\cf. Section~\[:S]\). Since $\Qh$ is a \mccl/, the product can be written
\vv-.5>
also as \>$\prod_i\,Q_{s_i}\bigl(\la-\]\tsum_{j>i}\gm\>\eps_{s_j}\bigr)$.
By formula \(SQS) and the last remark it suffices to verify that
\vvn-.1>
$$
S_\mu(m_1\>,m_2\>v_\mu)(\la)\,=\>
S_\mu(m_2\>,m_1\)v_\mu)\bigl(-\la+\gm\>\mu+\gm\>\zt'(m_1)+\zt''(m_2)\bigr)
$$
which follows from the property
\>$S(m_1\>,\)m_2)\>=\>\pim\bigl(S(m_2\>,\)m_1)\bigr)$,
commutation relations \(qf) and formula \(SmQv)\).
\epf
\Prop{CSQ}
For any $m_1\),m_2\in\Nc$ we have
$$
\CmQ\)(m_1\)\vmQ\,,\)m_2\>\vmQt)(\la)\,=\>
\SmQ\)(m_2\>,m_1\)\vmQ)\bigl(\la+\gm\>\zt''(m_2)\bigr)\,.
\Tag{CS}
$$
\endpro
\vsk-.3>\vsk0>
\Pf.
\bls1.06\bls
Recall that by the definition of \Vmod/s $m_1\)\vmQ$ and $m_2\>\vmQt$ are
constant \fn/s, since $m_1\),m_2\in\Nc$. If $\wsl(m_1)\ne\wsl(m_2)$, then the
expressions on both sides of formula \(CS) vanish. For $\wsl(m_1)=\wsl(m_2)$
the straightforward application of the definition of the \cform/ together with
Lemma~\[Smv] gives
$$
\CmQ\)(m_1\)\vmQ\,,\)m_2\>\vmQt)(\la)\,=\>
\SmQt\)(m_2\>,m_1\)\vmQt)\bigl(-\la+\gm\>\mu+\gm\>\zt'(m_1)\bigr)\,,
\vv.1>
\Tag{CSt}
$$
and using Lemma~\[SQQ] we complete the proof.
\epf
\Pf of Theorem~\[Cvv]\).
The first part of the theorem is an easy consequence of the definition of
the \cform/. The second part follows from Proposition~\[CSQ]\), Lemma~\[SQQ]\),
and the property
\vvn.1>
$$
\CmQ\)\bigl(f(\la)\>v\>,\)g(\la)\>\vti\))\,=\,
f(\la)\,g(-\la+\gm\>\nu)\,\CmQ\)(v\>,\vti\))
\vv.1>
\Tag{Cfg}
$$
for any $f,\)g\in\FunC$, $v\in\Fun\bigl(\MmQ[\)\nu\)]\)\bigr)$ and
$\vti\in\Fun\bigl(\MmQt[\)\nu\)]\)\bigr)$, \cf. \(CB)\).
\epf
\Pf of Theorems~\[kerpi] and \[kerC]\).
It is clear that
\ifMag
\vvm>\nl\vvm>
\cline{
$\Ker\pimQ\>=\>\bigl\lb\)v\in\Fun(\MmQ)\vert\CmQ\)(v\>,\vti)\)=\)0$ \>for any
\>${\vti\in\Fun(\MmQt)\)\bigr\rb}$.}
\else
\>$\Ker\pimQ=\bigl\lb\)v\in\Fun(\MmQ)\vert\CmQ\)(v\>,\vti)\)=\)0$ \>for any
\>${\vti\in\Fun(\MmQt)\)\bigr\rb}$.
\fi
So Theorem~\[kerpi] is \eqv/ to claim \)a)
\vv.06>
of Theorem~\[kerC]\). Claim \)b) of the latter follows from claim \)a) and
Theorem~\[Cvv]\).
\vsk.2>
Since $\pimQ\>\vmQ\)=\)\vmQta\ne 0$, the subspace $\Ker\pimQ\sub\Fun(\MmQ)$
\vv.1>
defines a proper submodule of $\MmQ$, and by Proposition~\[NmQ] we have that
$\Ker\pimQ\sub\Ker\SmQ$.
\vsk.3>
Let ${v\in\Ker\SmQ}$. We write it out as a linear combination of basis vectors:
\ifMag
$$
v\,=\!\tsum_{m\>\in\>\Nc}\]\phi_m(\la)\,m\>\vmQ\,.
\mmgood
$$
\else
\vv-.4>
${v\>=\!\]\sum_{m\>\in\>\Nc}\]\phi_m(\la)\,m\>\vmQ}$.
\fi
Then by Proposition~\[CSQ] for any $\mti\in\Nc$ we have
$$
\align
\CmQ\)(v\,,\)\mti\>\vmQt)(\la)\,=\sum_{m\>\in\>\Nc}\phi_m(\la)\>
\SmQ\)(\mti\>\vmQt\,,m\>\vmQ)\bigl(\la+\gm\>\zt''(\mti)\bigr)\,={} &
\\
{}=\>\SmQ\)(\mti\>\vmQt\,,\)v\))\bigl(\la+\gm\>\zt''(\mti)\bigr)\,={} &\,0\,.
\endalign
$$
Therefore, $\CmQ\)(v\>,\vti)\)=\)0$ \>for any \>$\vti\in\Fun(\MmQt)$
by property \(Cfg).
\vvgood
\epf

\Sect[R]{Rational dynamical \qg/ $\{\eslr$}
Introduce the spaces $\RatC$, $\RatV$ and $\Rat\bigl(\Hom(V,W)\bigr)$ similar
to the spaces $\FunC$, $\FunV$ and $\Fun\bigl(\Hom(V,W)\bigr)$, replacing in
the definitions \mef/s by \raf/s. Let $\Ratt=\RatC\oxC\RatC$.
\vsk.3>
{\bls1.06\bls We define the \em{\oalg/} $\elsr$ and \em{modules over the \rat/
\dqg/} $\eslr$ similar to the elliptic case with the following modification:
we replace the spaces of \mef/s by the respective spaces of \raf/s, substitute
the theta \fn/ $\tht(u)$ by the linear \fn/ $u\map u$, and set $\gm=1$.
For instance, formulae \(tf)\), \(tab)\), \(tth) and \(ell) become
\vvn-.1>
\vvnm.1>
$$
\gather
t_{ab}\,f(\la\+1\},\la\+2)\,=\>
f(\la\+1\!-\eps_a\>,\la\+2\!-\eps_b)\,t_{ab}\kern-2em
\\
\nn8>
t_{ac}\>t_{bc}\,=\,\){\la\+1_{ab}\!+1\over\la\+1_{ab}\!-1}
\ t_{bc}\>t_{ac}\,,\rlap{\kern4em for\quad $a\ne b$\,,}
\\
\nn12>
\tth_{ab}\,=\prod_{1\le c<a}\}\la\+1_{ca}
\prod_{1\le c<b}\}(\la\+2_{cb})\vpb{-1}\,t_{ab}\,,
\Tag{tthr}
\\
\ald
\nn11>
(t_{ab}\,v)\)(\la)\,=\,\ell_{ab}(\la)\>v(\la-\eps_a)
\qqq\text{for any}\quad v\in\RatV\,.\kern-4em
\Tag{ellr}
\\
\cnn-.1>
\endgather
$$
In the last formula $V\}$ is an \ermod/. The definitions of \hwm/s, \dwt/s,
etc.\ can be obviously transfered to the \rat/ case.
\vsk.2>}
The \rat/ case can be considered as a degeneration of the elliptic case
obtained by rescaling \var/s: $u\to\gm\>u$, $\la\to\gm\)\la$ and taking
the limit $\gm\to 0$.
\Par
Consider the limit $\Rb(\la)$ of the \rat/ version of the \Rm/ \(R)
as $u\to\8$:
\vvn-.2>
$$
\Rb(\la)\,=\sum_{a,b=1}^N\}E_{aa}\ox E_{bb}\,
+\]\sum_{\tsize{a,b=1\atop a\ne b}}^N\!
{E_{aa}\ox E_{bb}-E_{ab}\ox E_{ba}\over\la_{ab}}\;.
\vv-.6>
\Tag{Rb}
$$
It is a constant \sol/ of the \DYB/:
$$
\Rb\"{12}(\la-h\"3)\>\Rb\"{13}(\la)\>\Rb\"{23}(\la-h\"1)\,={}
\Rb\"{23}(\la)\>\Rb\"{13}(\la-h\"2)\>\Rb\"{12}(\la)\,.
$$
$\Rb(\la)$ is the simplest example of the Hecke type dynamical \Rm/,
see \Cite{EV1}.
\goodbm
\Par
Let $\Tcb\)=\>\sum_{a,b}\>E_{ba}\ox\)t_{ab}$, where $t_{ab}$ are the generators
\vvn-.1>
of $\elsr$ obeying the \rat/ version of commutation relations
\(tf)\,--\,\(tac)\). Then one can write relations \(tbc)\,--\,\(tac)
in the \Rm/ form:
\vvnn0:.2>
$$
\Rb\"{12}(\la\+2)\>\Tcb\"{13}\>\Tcb\"{23}\)=\,\)
\Tcb\"{23}\>\Tcb\"{13}\>\Rb\"{12}(\la\+1)\,.
\Tag{RTTb}
$$
\vsk-.1>
Let $V,\)W\}$ be \ermod/s. Then the \hmod/ $V\]\ox W\}$ is made into
an \ermod/ by the rule
\vvnn-.5:-.7>
$$
\ell_{ab}(\la)\vst{V\ox\)W}\,=
\>\tsum_{c=1}^N\,\ell_{cb}(\la-h\"2)\ox \ell_{ac}(\la)\,,
\vv-.1>
\Tag{dlr}
$$
and $t_{ab}$ acts on $\Rat(V\]\ox W)$ according to \(ellr)\).
\Par
Consider the following element in $\elsr$:
\vvnn0:-.1>
$$
\tsize\tN\>=\>\sumib\sign(\)\ib\))\;t_{N,\)i_N}\>\ldots\,t_{1,\)i_1}\,.
\vvmm0:-.5>
\Tag{tN}
$$
The product
\,${D\>=\!\!\prod_{\abn}\!\bigl(\la\+1_{ab}\}/\la\+2_{ab}\bigr)\;\tN}$
coincides with the top coefficient of the \rat/ version of the quantum
determinant $\Det\Tc(u)$, \cf. \(detTc). Hence, $D$ is a central element
in $\elsr$.
\goodbreak
\vsk.3>
Let $V\}$ be an \ermod/. The action of $D$ on $\RatV$ commutes with
multiplication by any \fn/ $\phi(\la)\in\RatC$ and, therefore, is given
by multiplication by a certain \fn/ $D(\la)\in\Rend V\}$. The module $V\}$
is called \em{\ndeg/} if $D(\la)$ is invertible for generic $\la$, and
\em{semi\-standard} if $D(\la)=1$. For instance, any \ndeg/ (standard) \hw/
\ermod/ is \ndeg/ (semistandard) in the new sense.
\vsk.2>
One can check that the element $D$ is group-like, it acts on the \Emod/
$V\]\ox W\}$ by
\vvnm-.1>
$$
D(\la)\vst{V\ox\)W}\>=\,D(\la-h\"2)\)\ox\>D(\la)\,.
\vvm.2>
$$
Therefore, a tensor product of \ndeg/ (semistandard) \ermod/s is \ndeg/
(semistandard)\).
\vsk.3>
Let $e_{ab}$, $a\),\bno$, be the standard generators of the Lie algebra $\gln$:
$$
[\)e_{ab}\,,\)e_{cd}\)]\,=\,\dl_{bc}\>e_{ad}-\dl_{ad}\>e_{cb}\,.
\Tag{eab}
$$
We identify the Lie algebra $\sln$ with the subalgebra of traceless elements
in $\gln$:
$$
\sln\,=\,
\bigl\lb\>\tsum_{a,b}x_{ab}\>e_{ab}\Vert\}\tsum_a\>x_{aa}=0\>\bigr\rb\,,
\vvmm-.4:-.5>
\Tag{sln}
$$
and $\hg$ with the subalgebra of diagonal elements in $\sln$:
\,$\hg=\bigl\lb\>\sum_a x_{aa}\>e_{aa}\Vert\}\tsum_a\>x_{aa}=0\>\bigr\rb$.
\,The standard basis of $\hg$ is $h_a\]=e_{aa}\]-e_{a+1,\)a+1}$, \,$\ano-1$.
\vvn.1>
The assignment \,$e_{ab}\map E_{ab}$, $a\),\bno$, makes $\CN\!$ into
the \em{vector \rep/} of $\gln$ and $\sln$.
\vsk.4>
Let $V\}$ be an \ermod/. The elements $\tth_{ab}$ act on $\RatV$
as \dif/ operators:
\vvnm-.3>
$$
(\tth_{ab}\,v)\)(\la)\,=\,\ellh_{ab}(\la)\>v(\la-\eps_a)\,,
\Tag{ellh}
$$
with coefficients $\ellh_{ab}(\la)\in\Rend V$. The module $V\}$ is called
\em{perturbative} if these coefficients have the following behaviour as $\la$
goes to infinity in a generic direction:
\vsk.1>
\atem
for any $\ano$ the \fn/ $\ellh_{aa}(\la)$ has a limit $\elti_{aa}$ which is
an invertible operator;
\bitem for any $a\),\bno$, ${a\ne b}$, the \fn/ $\la_{ab}\)\ellh_{ab}(\la)$
has a limit $\elti_{ab}$.
\vsk.3>\nt
The operators $\elti_{ab}$, $a\),\bno$, satisfy the following commutation
relations:
\vvn-.1>
$$
[\)x\>,\)\elti_{bc}\)]\>=\>(x\>,\)\eps_b\]-\eps_c)\,\elti_{bc}\,,\kern4em
[\>\elti_{aa}\,,\)\elti_{bc}\)]\>=\>0\,,
\vv-.1>
$$
for any $x\in\hg$ and $a,b,c=1\lc N$,
\vvn-.2>
$$
[\>\elti_{ab}\,,\)\elti_{ba}\)]\,=\,
(e_{aa}\]-e_{bb})\,\elti_{aa}\>\elti_{bb}
\vv-.9>
$$
for $a\ne b$, \>and
$$
[\>\elti_{ab}\,,\)\elti_{bc}\)]\,=\,\elti_{ac}\>\elti_{bb}
$$
for \pd/ $a,b,c$.
Hence, the assignment \,${e_{ab}\map\elti_{aa}\1\>\elti_{ab}}$
\vv.07>
\>for ${a\ne b}$, \,supple\-mented by the action of $\hg$, makes $V\}$ into
an \smod/ which we denote by $\Cc(V)$. We say that $V\}$ is a \em{\pert/} of
$\Cc(V)$. It is clear that $V\}$ coincide with $\Cc(V)$ as a vector space.
\Lm{ptndg}
Let $V\}$ be a perturbative \ermod/. Then $V\}$ is \ndeg/.
\endpro
\Pf.
The operator $D(\la)$ is invertible for generic $\la$ because it has
an invertible limit $\elti_{11}\]\ldots\)\elti_{NN}$ as $\la$ goes
to infinity in a generic direction.
\epf
\Lm{pertox}
A tensor product of perturbative \ermod/s is perturbative.
\endpro
\Ex
The assignment
\vvnn-.6:-1.5>
$$
\alignat2
\ellh_{aa}(\la)\, &{}\map\,1\,-\sum_{a<b\le N}\){E_{bb}\over\la_{ab}^2}\;, &&
\\
\nn7>
\ellh_{ab}(\la)\, &{}\map\;{E_{ab}\over\la_{ab}}\;, && \}a\ne b\,,
\\
\cnn-.15>
\nngood
\endalignat
$$
$a\),\bno$, makes $\CN\!$ into an \ermod/ $V\}$, which is a \pert/ of
the vector \rep/ of $\sln$. The module $V\}$ is isomorphic to the vector
\rep/ $U\]$ of $\eslr$, \cf. \(vect)\), \>by the following \iso/:
\vvn-.9>
\vvnm.9>
$$
\phi(\la)\,=\,\suan\,\tprod_{1\le b<a}\!\la_{ba}\,E_{aa}\,\in\,\Mor(V,U)\,.
$$
\enddemo
\Lm{CVsv}
Let $V\}$ be a perturbative \ermod/. If $V\}$ has a \wt/ \singv/,
then $\Cc(V)$ has a \singv/ of the same \wt/\/{\rm;}
\endpro
\Pf.
Let $v\in\RVmu$ be a \singv/.
\vv.1>
This means that for generic $\la$ the value $v(\la)$ belongs to the subspace
${K_\la=\Cap_{a,\)b}\)\Ker\]\ellh_{ab}(\la+\eps_a)\vst{\Vmu}\sub V\}}$.
\vv-.04>
It is clear that $\dim\)K_\la$ does not depend on $\la$ for generic $\la$.
\vv.08>
Moreover, $K_\la$ has a limit $K_\8$ as $\la$ goes to infinity in a certain
generic direction, and $\dim\)K_\8=\dim\)K_\la\ge 1$. To complete the proof
we observe that the subspace of \singv/s in $\Cc(V)[\)\mu\)]$ contains $K_\8$.
\mmgood
\epf
\Lm{CcV}
Let $V\}$ be a perturbative \ermod/. Then
\vsk.1>
\atem
if \,$\Cc(V)$ is \irr/, then $V\}$ is \irr/\/{\rm;}
\vsk.1>
\bitem
if \,$\Cc(V)$ is a \hwm/, then $V\}$ is a \hwm/ of the same \hw/\/{\rm;}
\vsk.1>
\bitem
if \,$\Cc(V)$ is a \Vmod/, then $V\}$ is isomorphic to a \Vmod/ of
the same \hw/.
\endpro
\Pf.
If $V\}$ is reducible, then $\FunV$ has a nontrivial proper \inv/
\$\FunC$-linear subspace $U\]$, which is a direct sum of its \wt/ components.
For each \wt/ component $\Umu$ consider the subspace $\Ulmu\sub V$ spanned
\vvm.08>
by values of \fn/s from $\Umu$ regular at $\la$. It is clear that
\vv.1>
$\dimC\Ulmu=\dimF\Umu$ for generic $\la$. Moreover, $\Ulmu$ has a limit
\vv.1>
$\Uimu$ as $\la$ goes to infinity in a certain generic direction,
\vv.1>
$\dim\)\Uimu=\dim\)\Ulmu$, and the direction can be taken the same
for all $\mu$. Then $U_\8=\Plus_\mu\)\Uimu$ is a nontrivial proper
\vvmm0:-.1>
\inv/ subspace in $\Cc(V)$, since $\dim\)\Uimu<\dim\)\Vmu$ at least
for some $\mu$. Claim \)a) \)is proved.
\vskm.3:.2>
Let $v$ be the \hwv/ of $\Cc(V)$. Then the constant \fn/ ${v\in\RatV}$ is
a \regsv/ by the \wt/ reason. Any \wt/ subspace $\Vmu$ has a basis given by
vectors of the form $\eabk\)v$. Then the corresponding \fn/s $\thabk\)v$ span
$\RVmu$, which proves claim \)b)\). Claim \)c) \)follows from claim \)b)
\)and comparison of dimensions of \wt/ subspaces.
\epf
We say that an \smod/ $V\}$ is \em{admissible} if $V\}$ is a \dhmod/. Notice
that any \hw/ \smod/ is admissible. In Section~\[:J] we define a functor $\Ec$
from the category of \asmod/s to the category of semistandard \ermod/s,
\cf. Theorem~\[Ec]\). We summarize the properties of this functor
in the next two theorems.
\Th{EcV}
Let $V\}$ be an \asmod/. Then
\vsk.1>
\atem
$\Ec(V)$ coincides with $V\}$ as an \hmod/\/{\rm;}
\vsk.1>
\bitem
$\Ec(V)$ is a \pert/ of $V\}$, moreover, $\elti_{aa}=1$ for any $\ano${\rm;}
\vsk.1>
\bitem
$\Ec(V)$ is a semistandard \ermod/\/{\rm;}
\vsk.1>
\bitem
if \>${v\in V\}}$ is a \wt/ \singv/, then \>${v\in\FunV}$ considered as
a constant \fn/, is a standard \singv/ of the same \wt/ for $\Ec(V)${\rm;}
\vsk.1>
\bitem
if \>$V\}$ is a \hwm/, then $\Ec(V)$ is a standard \hwm/ with the same \hw/
and \hwv/\/{\rm;}
\vsk.1>
\bitem
if \>$V\}$ is a \Vmod/, then $\Ec(V)$ is isomorpic to the standard \Vmod/
with the same \hw/ and \hwv/.
\endpro
\Th{EcUV}
Let \)$U,\)V\}$ be \hw/ \smod/s. Then
\vsk.1>
\atem
any element of \,$\Mor\bigl(\Ec(U)\),\Ec(V)\bigr)$ is a constant \fn/\/{\rm;}
\vsk.1>
\bitem
the map \,${\Hom_{\sln}\](\)U,V)\)\to\>\Mor\bigl(\Ec(U)\),\Ec(V)\bigr)}$
defined by $\Ec$ is an \iso/\/{\rm;}
\vsk.1>
\bitem
the above \iso/ coincides with the restriction of the natural embedding
of \,$\Hom(\)U,V)$ into $\Fun\bigl(\Hom(\)U,V)\bigr)$.
\endpro
\vsk-.1>
\nt
The theorems are proved in Section~\[:J]\).
\Par
Let $V\}$ be a \hw/ \smod/ with \hw/ $\mu$ and \hwv/ $v$. Let $\Sh$ be
the $\sln$ \Sform/ on $V\}$, and let $S_\mu$ be the \dSpair/ for $\Ec(V)$.
\Prop{kerS}
$\Ker S_\mu\)=\>\Fun(\)\Ker\Sh\))$.
\endpro
\Pf.
\bls1.06\bls
Since $\Ker\Sh$ is an $\sln$ submodule of $V\}$, then $\Ec(\Ker\Sh)$ is
\vv-.04>
an $\eslr$ submodule of $\Ec(V)$, and $\Fun(\)\Ker\Sh\))\)\sub\)\Ker S_\mu$
by Proposition~\[Virr]\). On the other hand, $V\}/\]\Ker\Sh$ is an \irr/
\smod/, therefore, $\Ec(V)\big/\Ec(\Ker\Sh)=\Ec\bigl(V\}/\]\Ker\Sh\bigr)$
is an \irr/ \ermod/ by Lemma~\[CcV]\).
Hence, $\Ker S_\mu\)\sub\)\Fun(\)\Ker\Sh\))$ by Proposition~\[Virr]\).
\epf
\Lm{S8}
Let $a_i\ne b_i$ for any $i=1\lc k$ and let $c_j\ne d_j$ for any $j=1\lc l$.
Then
\ifMag
\vvnm-.9>
$$
\align
(-1)^k\,\tprod_{i=1}^k\)\la_{a_ib_i}\>\tprod_{j=1}^l\)\la_{c_jd_j}\;
& S_\mu(\thabk\>,\)\thcdl\>v_\mu)\,\to{}
\\
\nn2>
& {}\!\]\to\,\Sh(\eabk\)v_\mu\),\)\ecdl\>v_\mu)\kern-1em
\\
\cnn-.2>
\endalign
$$
\else
\vvn-.3>
$$
(-1)^k\,\tprod_{i=1}^k\)\la_{a_ib_i}\>\tprod_{j=1}^l\)\la_{c_jd_j}\;
S_\mu(\thabk\>,\)\thcdl\>v_\mu)\,\to\,
\Sh(\eabk\)v_\mu\),\)\ecdl\>v_\mu)\kern-1em
\vv-.2>
$$
\fi
as $\la$ goes to infinity in a generic direction.
\vvgood 
\endpro

\Sect[Fd]{Finite-dimensional \hwm/s over $\{\esl$}
In this section we assume that ${\gm\)\nin\)\QtQ}$ and consider only \ndeg/
\dwt/s. We do not mention this assumption explicitly. To save space we
usually formulate the results only for standard \hw/ \emod/s if they can be
generalized to arbitrary \hwm/s by pulling back through \aut/s \(faut)\).
\Par
Let $v$ be a standard \singv/ of \wt/ $\mu$, and $k$ be a nonnegative integer.
\vsk.8>
\Lm{abba}
Let $a<b$. Then \;$\dsize\tth_{ab}\,\tth_{ba}^{\>k}\)v\,=\,
-\,{\tht\bigl((\)\mu_{ab}\]-k+1)\>\gm\bigr)\,\tht(k\)\gm)\over\tht(\la_{ab})\,
\tht\bigl(\la_{ab}\]-(\)\mu_{ab}\]-2\)k+2)\>\gm\bigr)}\ \tth_{ba}^{\>k-1}\)v$.
\endpro
\vsk-.3>
\vsk0>
\Pf.
Take formula \(thadcbk) for $a=b$, $c=d$, and replace $c$ by $b$.
Since $\tth_{aa}\)v\)=\)\tth_{bb}\>v\)=\)v$, \,we have
$$
\tth_{ab}\,\tth_{ba}^{\>k}\)v\,=\,-\,
{\tht\bigl(\la_{ab}\+1\!-\la_{ab}\+2\!+(k-1)\>\gm\bigr)\,
\tht(k\)\gm)\over\tht(\la\+1_{ab})\,\tht(\la\+2_{ab})}
\ \tth_{ba}^{\>k-1}\)v\,.
\vv-.1>
$$
Now apply convention \(la12) and observe that the vector
$\tth_{ba}^{\>k-1}\)v$ has \wt/ $\mu-(k-1)\>(\eps_a\]-\eps_b)$.
\vv.06>
Since $\nu_{ab}=(\)\nu\),\eps_a\]-\eps_b)$ for any $\nu\in\hga\}$
and $(\eps_a\]-\eps_b\>,\eps_a\]-\eps_b)=2$, the lemma is proved.
\epf
\Cr{ne0}
Let $a<b$ and \)$\mu_{ab}\nin\)\lb\)0\),1\)\lc k-1\)\rb$.
Then \,$\tth_{ba}^{\>k}\)v\)\ne 0$.
\endpro
\Pf.
By Lemma~\[abba]
\vv-.8>
$$
\tth_{ab}^{\>k}\,\tth_{ba}^{\>k}\)v\,=\,(-1)^k\,\prod_{j=0}^{k-1}\,
{\tht\bigl((\)\mu_{ab}\]-j)\>\gm\bigr)\,\tht\bigl((j+1)\>\gm\bigr)
\over\tht(\la_{ab}\]-j\)\gm)\,
\tht\bigl(\la_{ab}\]-(\)\mu_{ab}\]-j)\>\gm\bigr)}\ v\ne 0\,.
\vv-1.5>
$$
\epf
\Lm{cdba}
\;$\tth_{cd}\,\tth_{a+1,\)a}^{\>k}\)v\)=\)0$ \;for any $c<d$,
$(c\>,d\))\ne(a\>,a+1)$. \;$\tth_{cc}\,\tth_{a+1,\)a}^{\>k}\)v\)=\)v$
\;for any $c\ne a\>,a+1$.
\endpro
\vsk->\vsk0>
\Lm{regsi}
If \)${(\)\mu\),\al_a)\)=\)k-1}$ \)and \>${\tth_{a+1,\)a}^{\>k}\)v\)\)\ne\)0}$,
\>then ${\prod_{j=1}^k\tht(\la_{a,\)a+1}\]+j\)\gm)\;\tth_{a+1,\)a}^{\>k}\)v}$
\vv-.4>
\>is a standard \singv/ of \wt/ $\mu-k\)\al_a$.
\mmgood
\endpro
\Pf.
By Lemmas~\[abba] and \[cdba] the \fn/ $\tth_{a+1,\)a}^{\>k}\)v$ is a \singv/
and
\ifMag
$$
\tth_{cc}\,\tth_{a+1,\)a}^{\>k}\)v\,=\,v
\vvm.1>
$$
\else
$\tth_{cc}\,\tth_{a+1,\)a}^{\>k}\)v=v$
\fi
for any $c\ne a\>,a+1$.
On the other hand, it follows from formulae \(thac) and \(thbd) that
\ifMag
\<aa1>
$$
\gather
\tth_{aa}\>\tth_{a+1,\)a}^{\>k}\)v\,=
\,\){\tht(\la_{a,\)a+1}\]+k\)\gm)\over\tht(\la_{a,\)a+1})}
\ \tth_{a+1,\)a}^{\>k}\)v\,,\kern-1em
\Tag{aa1}
\\
\nn13>
{\align
\tth_{a+1,\)a+1}\>\tth_{a+1,\)a}^{\>k}\)v\,=\,\)
{\tht(\la_{a,\)a+1}\]-\gm\)\mu_{a,\)a+1}\]+k\)\gm)\over
\tht(\la_{a,\)a+1}\]-\gm\)\mu_{a,\)a+1}\]+2\)k\)\gm)}\ \tth_{a+1,\)a}^{\>k}\)v
\,={}\!\}\kern-1em &
\Tag{a1a}
\\
\nn8>
{}={\tht(\la_{a,\)a+1}\]+\gm)\over\tht\bigl(\la_{a,\)a+1}\]+(k+1)\>\gm\bigr)}
\ \tth_{a+1,\)a}^{\>k}\)v\kern-1em &
\endalign}
\endgather
$$
\else
\vvn-.2>
$$
\gather
\tth_{aa}\>\tth_{a+1,\)a}^{\>k}\)v\,=
\,\){\tht(\la_{a,\)a+1}\]+k\)\gm)\over\tht(\la_{a,\)a+1})}
\ \tth_{a+1,\)a}^{\>k}\)v\,,\kern-1.4em
\Tag{aa1}
\\
\nn15>
\tth_{a+1,\)a+1}\>\tth_{a+1,\)a}^{\>k}\)v\,=\,\)
{\tht(\la_{a,\)a+1}\]-\gm\)\mu_{a,\)a+1}\]+k\)\gm)\over
\tht(\la_{a,\)a+1}\]-\gm\)\mu_{a,\)a+1}\]+2\)k\)\gm)}\ \tth_{a+1,\)a}^{\>k}\)v
\,={\tht(\la_{a,\)a+1}\]+\gm)\over\tht\bigl(\la_{a,\)a+1}\]+(k+1)\>\gm\bigr)}
\ \tth_{a+1,\)a}^{\>k}\)v\kern-1.6em
\Tag{a1a}
\\
\cnn-.2>
\endgather
$$
\fi
because ${\mu_{a,\)a+1}=(\)\mu\),\al_a)=\)k-1}$.
Hence, multiplying $\tth_{a+1,\)a}^{\>k}\)v$ by the product
\mline
\vv-.5>
${\prod_{j=1}^k\tht(\la_{a,\)a+1}\]+j\)\gm)}$ \)one gets a standard \singv/.
\vvm.5>
\wwgood-.8:.8>
\epf
\Prop{infdim}
An \irr/ standard \hw/ \emod/ with \hw/ $\mu$ is in\fd/
if $\mu$ is not a \dint/.
\endpro
\Pf.
\bls1.06\bls
Let $v$ be the \hwv/. Assume that $\mu$ is not a \dint/,
and let $a$ be \st/ $(\)\mu\),\al_a)\nin\Zp$. Then the \fn/s
$v\,,\;\tth_{a+1,\)a}\)v\,,\;\tth_{a+1,\)a}^{\>2}\)v\,,\alb\;\ldots{}$
are linearly independent over $\FunC$ because all of them are nonzero
by Corollary~\[ne0] and they have distinct \wt/s \wrt/ the action of $\hg$.
\epf
\Prop{aksub}
Let ${(\)\mu\),\al_a)\)=\)k-1}$. Then ${\tth_{a+1,\)a}^{\>k}\)v_\mu}$
generates a submodule of the \Vmod/ $M_\mu$ isomorphic to the \Vmod/
$M_{\mu\)-\)k\al_a}$.
\endpro
\Pf.
By Lemma~\[regsi] there is a nontrivial morphism
${\phi\in\Mor(M_{\mu\)-\)k\al_a},M_\mu)}$ which sends
the \hwv/ ${v_{\mu\)-\)k\al_a}\!\in M_{\mu\)-\)k\al_a}}$ to
$\prod_{j=1}^k \tht(\la_{a,\)a+1}\]+j\)\gm)\;\tth_{a+1,\)a}^{\>k}\)v_\mu$,
and it remains to show that $\phi$ is an embedding. In other words, one has
\vvnm.1>
to prove that the induced map ${\Fun(M_{\mu\)-\)k\al_a})\to\)\FunM}$ is
injective.
\vsk.3>
For any \Vmod/ $M_\nu$ set
\vv.1>
$\Fun_j(M_\nu)=\lb\>x\>v_\nu\vert x\in\els\,,\ \,\deg\>(x)\le j\)\rb$ and
\vvn-.1>
$$
\Funb(M_\nu)\,=\,
\Plus_{j=0}^\8\,\Fun_j(M_\nu)\)\big/\}\Fun_{j-1}(M_\nu)\,.
$$
\par
\bls1.1\bls\nt
Let $\Nc$ be given by \(BNc)\). It is clear that the set
${\lb\>x\>v_\nu\vert x\in\Nc\,,\ \,\deg\>(x)\le j\)\rb}$
is a basis of $\Fun_j(M_\nu)$ over $\FunC$, and the set
${\lb\>x\>v_\nu\vert x\in\Nc\,,\ \,\deg\>(x)=j\)\rb}$ induces
a basis of $\Fun_j(M_\nu)\)\big/\}\Fun_{j-1}(M_\nu)$.
We identify $\Funb(M_\nu)$ with the space of \pol/s in \var/s
$u_{21}\), u_{31}\)\lc u_{N,\)N-1}$ with coefficients in $\FunC$:
for any monomial ${\tth_{b_1c_1}\}\ldots\)\tth_{b_jc_j}\!\in\Nc}$
a class of the \fn/ $\tth_{b_1c_1}\}\ldots\)\tth_{b_jc_j}\)v_\nu$
in the quotient space $\Fun_j(M_\nu)\)\big/\}\Fun_{j-1}(M_\nu)$
is mapped to \)$u_{b_1c_1}\}\ldots\)u_{b_jc_j}$.
\vsk.27>
The map ${\phi\):\Fun(M_{\mu\)-\)k\al_a})\to\)\FunM}$ induces
a map $\phib:\)\Funb(M_{\mu\)-\)k\al_a})\to\>\Funb(M_\mu)$.
Consider $\phib$ as a map from $\FunC\)[\)u_{21}\lc u_{N,\)N-1}\)]$ to itself.
By formula \(thabcdk) we find that
\ifMag
\,${\phib\>(\)u_{b_1c_1}\}\ldots\)u_{b_jc_j})\,=\,
f\bigl(\la-\gm\tsum_{i=1}^j\]\eps_{b_i}\bigr)\,
u_{b_1c_1}\}\ldots\)u_{b_jc_j}\>u_{a+1,\)a}^k}\>$
\vv-.5>
\else
$$
\phib\>(\)u_{b_1c_1}\}\ldots\)u_{b_jc_j})\,=\,
f\bigl(\la-\gm\tsum_{i=1}^j\]\eps_{b_i}\bigr)\,
u_{b_1c_1}\}\ldots\)u_{b_jc_j}\>u_{a+1,\)a}^k
\vv-.5>
$$
\fi
where $f(\la)\)=\prod_{j=1}^k \tht(\la_{a,\)a+1}\]+j\)\gm)$.
Hence, $\phib$ is injective, and so does $\phi$.
\epf
\vskm-.5:0>
{}From now on till the end of the section fix a \dint/ $\mu$ and set
$k_a=(\)\mu\),\al_a)+1$, \,$\ano-1$. Notice that $k_a\}\in\Zpp$ for any $a$.
Denote by $Z_\mu$ the subspace in $\FunM$ generated over $\els$ by \fn/s
\>$\tth_{a+1,\)a}^{\>k_a}\)v_\mu$, \,$\ano-1$. Notice that the \fn/s
$\tth_{a+1,\)a}^{\>k_a}\)v_\mu$ are \regsv/s, \cf. Lemma~\[regsi]\).
\vsk.3>
Let $S_\mu$ be the \Spair/ for $M_\mu$, \cf. \(SmQ)\).
\Prop{ZiS}
$Z_\mu\sub\)\Ker S_\mu$.
\endpro
\Pf.
Lemma~\[regsi] implies that $Z_\mu$ is an \inv/ \$\FunC$-linear subspace
in $\FunM$, therefore it defines a submodule of $M_\mu$. Hence, the statement
follows from Proposition~\[NmQ]\).
\epf
\Th{ZS}
$\Ker S_\mu\)=\)Z_\mu$ for generic \)$\gm$.
\endpro
\Pf.
Since both $\Ker S_\mu$ and $Z_\mu$ are direct sums of their \wt/ components,
we have to prove that $\KSmn\)=\)\Zmn$ for any $\nu\le\mu$.
Notice that $(\Ker S_\mu)[\)\mu\)]=\FunC\>v_\mu=Z_\mu[\)\mu\)]$.
\vsk.3>
By Proposition~\[ZiS] it suffices to prove that $\KSmn$ and $\Zmn$ have
same dimensions. For doing this we employ the deformation argument.
\vsk.2>
Consider the \Vmod/ $\Mslm\}$ of \hw/ $\mu$ over $\sln$ and the \Sform/
$\Sslm\}$ on it. Let $\Mrm\}=\)\Ec(\Mslm)$, see Theorem~\[EcV]\), and let
$\Srm\}$ be the corresponding \dSpair/. Recall that, $\Mrm\}$ is isomorphic to
the \Vmod/ of \hw/ $\mu$ over $\eslr$. By abuse of notation we denote
generators of $\els$ and $\elsr$ by the same letters and write $v_\mu$ for
the \hwv/ of each of the modules $M_\mu$, $\Mslm\}$ and $\Mrm\}$.
\vsk.2>
Let $\Zslm$ be the $\sln$ submodule in $\Mslm$ \gby/ \singv/s
\>$e_{a+1,\)a}^{\>k_a}\)v_\mu$, \,$\ano-1$, and let $\Zrm$ be the subspace
in $\RatMm$ generated over $\elsr$ by \fn/s \>$\tth_{a+1,\)a}^{\>k_a}\)v_\mu$,
\,$\ano-1$. It is known that $\Ker\Sslm\}=\Zslm$. By Theorem~\[EcV] and
Proposition~\[kerS] we have that $\Ker\Srm\}=\RatSm$.
\vskm-.4:0>
\vskgood-:>
\Lm{te}
${\prod_{j=1}^{k_a}(\la_{a,\)a+1}\]+j\))\;\tth_{a+1,\)a}^{\>k_a}\)v_\mu\>=\>
(-1)^{\)k_a}\>e_{a+1,\)a}^{k_a}\)v_\mu}$ \,\)for any $\ano-1$.
\endpro
\vskm-1.1:-.9>\vsk0>
\Pf.
\bls1.1\bls
Both vectors
$\prod_{j=1}^{k_a}(\la_{a,\)a+1}\]+j\))\;\tth_{a+1,\)a}^{\>k_a}\)v_\mu$
\vvn-.1>
and $e_{a+1,\)a}^{k_a}\)v_\mu$ have \wt/ $\mu-k_a\al_a$. Thus they are
pro\-portional over $\RatC$ because $\dimC\Mslm[\)\mu-k_a\al_a]\)=\)1$.
The proportionality coefficient is a constant \fn/, since both of them are
standard \singv/s by Theorem~\[EcV] and the \rat/ version of Lemma~\[regsi]\).
The constant equals $1$ because $\la_{a,\)a+1}^{k_a}\)
\tth_{a+1,\)a}^{\>k_a}\)v_\mu\)\to\>(-1)^{\)k_a}\>e_{a+1,\)a}^{k_a}\)v_\mu$
\>as $\la$ goes to infinity in a generic direction.
\epf
\Cr{ZSr}
$\Ker\Srm\}=\)\Zrm\}=\)\RatZm$.
\endpro
\nt
Since the elliptic case is a deformation of the \rat/ one we obtain that
\mmgood
for generic~$\gm$
\ifMag
$$
\align
\dimF\]\KSmn\, &{}\le\,\dimR\](\Ker\Srm)[\)\nu\)]\,={}
\Tag{dims}
\\
\nn4>
{}=\,\dimC(\Ker\Sslm)[\)\nu\)]\, &
{}=\,\dimC\Zslm[\)\nu\)]\,=\,\dimR\Zrm[\)\nu\)]\,\le{}
\\
\nn5>
& {}=\,\dimF\Zmn\,\le\,\dimF\]\KSmn\,.\kern-2em
\endalign
$$
\else
$$
\align
\dimF\]\KSmn\> &{}\le\>\dimR\](\Ker\Srm)[\)\nu\)]\>=\>\dimC(\Ker\Sslm)[\)\nu\)]
\>={}
\Tag{dims}
\\
\nn4>
{}=\>\dimC\Zslm[\)\nu\)]\> &{}=\>\dimR\Zrm[\)\nu\)]\>\le\>\dimF\Zmn\>\le\>
\dimF\]\KSmn\,.\kern-2.5em
\endalign
$$
\fi
Here the last inequality is due to Proposition~\[ZiS]\). Therefore,
all the dimensions in \(dims) are the same, which proves the theorem.
The rest of the proof is a justification of this informal reasoning.
\Par
For a pair of sequences $\aa=(a_1\lc a_j)$ and $\bb=(b_1\lc b_j)$ let
\>$\tth_{\aa\bb}=\tth_{a_1b_1}\}\ldots\)\tth_{a_jb_j}$. \,Set
\vvn-.4>
$$
\wsl\)(\aa\),\]\bb\))\>=\>\wsl\)(\)\tth_{\aa\bb})\>=
\tsum_{i=1}^j\>(\)\eps_{a_i}\!-\)\eps_{b_i})\,,\qqq
\zt\)(\aa\),\]\bb\))\>=\tsum_{i=1}^j\>\eps_{b_i}\,.
\vvmm0:-.3>
$$
Let \,$\Mc=\lb\)(\aa\),\]\bb\))\vert\tth_{\aa\bb}\in\Nc\)\rb$,
\>the set $\Nc$ being defined in \(BNc)\). Set
\vvn-.4>
$$
\Mcb\>=\>\lb\)m\in\Mc\>\vert\]\wsl(m)=\bt\)\rb\qquad\text{and}\qquad
\tsize\Mcmb\>=\>\Cup_{a=1}^{N-1}\>\Mc\)[\)\bt+k_a\al_a]\,.
\vv-.2>
$$
\par
Let $\bt=\nu-\mu$. Introduce a matrix $\Agm(\la)$ with entires labeled
by pairs of elements of $\Mcb$:
\vvm-.6>
$$
A^\gm_{m,m'}(\la)\,=\,
S_\mu(\)\tth_m\>,\)\tth_{m'}\)v_\mu)\bigl(\la+\gm\>\zt(m)\bigr)\,,
\Tag{Agm}
$$
and a matrix $\Bgm(\la)$ with entries labeled by pairs
$(m\),m')\in\Mcb\)\x\Mcmb$ and given by the rule:
\vvn.1>
$$
m'\>\tth_{a+1,\)a}^{\>k_a}\)v_\mu\,=
\sum_{m\>\in\>\Mcb\!}\!B^\gm_{m,m'}(\la)\,\)\tth_m\)v_\mu\,,\kern4em
\Rlap{m'\}\in\Mc\)[\)\bt+k_a\al_a]\,.}
\vvmm0:-.4>
$$
Both $\Agm(\la)$ and $\Bgm(\la)$ are \mef/s of $\gm\>,\la$.
Moreover, if $\gm\to 0$, then
$$
\Agm(\gm\)\la)\to\Ar(\la)\,,\kern4em \Bgm(\gm\)\la)\to\Br(\la)\,,
\Tag{AB}
$$
where $\Ar(\la)$ and $\Br(\la)$ are defined in the same way for the \rat/ case
\vv.07>
(\)in formula \(Agm) for the \rat/ case
the argument in \rhs/ is $\la+\zt(m)$)\).
\vsk.3>
Each matrix naturally defines a linear map: $\Agm(\la)$ and $\Ar(\la)$ act on
\vvmm.2:-.6>\mline
${\MCb\>=\]\!\Plus_{m\>\in\>\Mcb\!}\!\C\>m}$, whilst $\Bgm(\la)$ and $\Br(\la)$
map ${\Plus_{m\>\in\>\Mcmb}\!\!\C\>m}$ \)to \)$\MCb$.
\vskm.4:.5>
Given \fn/s $\phi_m(\la)\in\FunC$, $m\in\Mcb$, consider a vector
\vvmm.2:-.2>\mline
\>$\phc\)(\la)\>=\!\]\sum_{m\>\in\>\Mcb\!}\phi_m(\la)\,m\in\MCb$.
\,Set \,$\phiv\)(\la)\>=\!
\sum_{m\>\in\>\Mcb\!}\!\phi_m(\la)\,\tth_m\)v_\mu\)\in\>\FunM$.
\Lm{kerim}
\vvn.2>
\atem
$\phiv\in\KSmn$ \,iff \,$\phc\)(\la)\in\Ker\]\Agm(\la)$ for generic $\la${\rm;}
\vsk.2>
\bitem
$\phiv\in(\Ker\Srm)[\)\nu\)]$ \,iff \,$\phc\)(\la)\in\Ker\]\Ar(\la)$
for generic $\la${\rm;}
\vsk.3>
\bitem
$\phiv\in\Zmn$ \,iff \,$\phc\)(\la)\in\Im\]\Bgm(\la)$ for generic $\la${\rm;}
\vsk.2>
\bitem
$\phiv\in\Zrm[\)\nu\)]$ \,iff \,$\phc\)(\la)\in\Im\]\Br(\la)$
for generic $\la$.
\endpro
\Pf.
Claims \)a) and \)b) follow from formulae \(pimth)\,--\,\(SmQv) and their \rat/
versions, \resp/. Claims \)c) and \)d) are straightforward.
\wwgood-.8:.8>
\epf
Corollary~\[ZSr] implies that $\Ker\]\Ar(\la)=\)\Im\]\Br(\la)$ for generic
$\la$. The standard deformation reasoning, \cf. \(AB)\), shows that
\vvn.1>
\ifMag
\vvn-.2>
$$
\align
\dimC\]\Ker\]\Agm(\la)\,\le\,\dimC\]\Ker\]\Ar(\la)\,=\,
\dimC\]\Im\]\Br(\la)\,={}\!\}\kern-2em &
\Tag{diml}
\\
\nn5>
{}=\,\dimC\]\Im\]\Bgm(\la)\,\le\,\dimC\]\Ker\]\Agm(\la)\kern-2em &
\\
\cnn-.1>
\endalign
$$
\else
$$
\dimC\]\Ker\]\Agm(\la)\>\le\>\dimC\]\Ker\]\Ar(\la)\>=\>\dimC\]\Im\]\Br(\la)
\>=\>\dimC\]\Im\]\Bgm(\la)\>\le\>\dimC\]\Ker\]\Agm(\la)\kern-.2em
\vv.1>
\Tag{diml}
$$
\fi
\par
\bls1.09\bls\nt
for generic $\la$, provided $\gm$ is generic. Here the last inequality is due
to Proposition~\[ZiS], which implies that $\Im\]\Bgm(\la)\sub\Ker\]\Agm(\la)$
for generic $\la$. Therefore, all the dimensions in \(diml) are the same.
Hence, $\Ker\]\Agm(\la)\)=\)\Im\]\Bgm(\la)$ for generic $\la$, provided $\gm$
is generic. By Lemma~\[kerim] we see that $\KSmn=\Zmn$ for generic $\gm$.
Theorem~\[ZS] is proved.
\epf
Let $V_\mu$ be the \irr/ standard \hw/ \emod/ of \hw/ $\mu$. Let $N_\mu$ be
the \$\esl$-submodule of $M_\mu$ \st/ $\Fun(N_\mu)=\Ker S_\mu$. We have
$V_\mu=M_\mu\)/N_\mu$, \cf. Corollary~\[VmQ]\). Let $\Vslm$ be the \irr/ \smod/
of \hw/ $\mu$. Set $\dmn=\dimC\Vsmn$.
\Pf of Theorem~\[fdim]\).
It follows from the proof of Theorem~\[ZS] that for generic $\gm$
$$
\align
\dimC\Vmn\> &{}=\>\dimC\Mmn-\dimF\]\KSmn\>={}
\Tagg{dmn1}
\\
\nn4>
&{}=\>\dimC\Msmn-\dimC(\Ker\Sslm)[\)\nu\)]\>=\>\dimC\Vsmn\>=\>\dmn\,.
\ifMag\kern-1em\else\kern-1.4em\fi
\endalign
$$
Since $\dimF\]\KSmn$ can jump only up at a specific $\gm$, we have that
\vvm.1>
$\dimC\Vmn\alb\>\le\>\dmn$ for arbitrary $\gm$.
\epf
\Pf of Theorem~\[dimu]\).
{\bls1.06\bls
We have already proved that ${\dimC\Vmn\>=\>\dmn}$ for any $\nu$, provided
$\gm$ is generic, \cf. \(dmn1)\). Consider an \emod/ $U_\mu$ defined in
the following way. Let $U_\mu=\)V_\mu$ if $\gm$ is generic, that is,
if ${\dimC\Vmn\>=\>\dmn}$. Otherwise, define $U_\mu$ by analytic continuation
from generic $\gm$. It is not difficult to justify the given definition of
$U_\mu$, and to see that $U_\mu$ is a \hw/ \emod/ with \hw/ $\mu$ and
$\dimC\Umn\>=\>\dmn$ for any $\nu$.
\vsk.2>}
If $\mu$ is a \dint/, then $w\cddt\mu<\mu$ and $d_\mu[\)w\cddt\mu\)]\)=\)0$ 
for any nontrivial element $w\in W\}$. Hence, $U_\mu$ is \irr/
by Corollary~\[hwirr] \>if ${\gm\)\nin\)\QtQ}$, that is, $U_\mu=\)V_\mu$.
The theorem is proved.
\epf
\Rem
There is another approach to constructing \fd/ \irrep/ of $\esl$. One can start
from the vector \rep/ of $\esl$ and apply the fusion procedure technique
developed in the nondynamical case, \cf. \Cite{C}\), \Cite{N}\).
If ${\gm\)\nin\)\QtQ}$, then any \irr/ \fd/ standard \hw/ \emod/ can be
obtained in this way. The \sym/ and exterior powers of the vector \rep/
of $\esl$ have been constructed by this technique in \Cite{FV2}\).
We will address this approach elsewhere.
\enddemo

\Sect[J]{Definition of functor $\{\Ec$}
In this section we construct a functor from the category of \asmod/s to
the category of semistandard \ermod/s. The construction is similar to the
construction of the functor from the category of \fd/ \smod/s to the category
of \rat/ \rep/s of the exchange \qg/ $\Fsl$, developed in \Cite{EV2}\).
In Section~\[:ex] we discuss the relation of these two constructions in detail.
\vskm.3:.2>
Let \>$\ngp=\)\bigl\lb\)\tsum_{a<b}x_{ab}\>e_{ab}\>\bigr\rb$ \>and
\>$\ngm=\)\bigl\lb\)\tsum_{a>b}x_{ab}\>e_{ab}\>\bigr\rb$ be the standard
nilpotent subalgebras in $\sln$, and let $\bgpm\]=\)\hg\)\oplus\ngpm$. Set
\vvnm-.4>
$$
\Xi\,=\,\){1\over 2}\>\tsuan e_{aa}^{\)2}\,-\,
{1\over 2N}\>\bigl(\tsuan e_{aa}\bigr)\vpb{2}\in\,U(\)\hg\))\,.
$$
\vskm-.8:-.1>
\vsk0>
\Prop{abrr}
There exists a unique power series $\Jc(\la\);z)$ in $z$ with coefficients in
\$\Usl\ox\Usl$-valued \fn/s of \>$\la$ with the properties\/{\rm:}
\vsk.1>
\atem
$\Jc(\la\);z)$ satisfies the \em{\rat/ ABRR \eq/}
\vsk-.8>
$$
\Jc(\la\);z)\,\bigl(1\ox(\la-z\>\Xi\))\bigr)\,=\,\bigl(1\ox(\la-z\>\Xi\))\>+
\>z\!\tsum_{\abn}\!e_{ab}\ox e_{ba}\bigr)\,\Jc(\la\);z)\,;
\vv-.2>
\Tag{Jc}
$$
\bitem
the coefficients of the series \>$\bigl(\Jc(\la\);z)-1\)\bigr)$ are
\$\,U(\)\bgp)\>\ngp\]\ox\)\ngm U(\)\bgm)$-\)valued \fn/s of $\la$.
\vvgood
\endpro
\Pf.
Let \>$\Jc(\la\);z)\)=\sum_{k=0}^\8 \Jc_k(\la)\>z^k\}$. Equation \(Jc) is \eqv/
to certain recurrence relations for coefficients $\Jc_k(\la)$ with the initial
condition $\Jc_0(\la)=1$. It is straightforward to verify that at each step
the recurrence relations uniquely determine $\Jc_k$ from $\Jc_0\lc\Jc_{k-1}$.
\epf
It follows from the proof of the last proposition that the coefficients
\vv.1>
of the series $\Jc(\la\);z)$ are \raf/s of $\la$, and for any $x\in\hg$
\vvn-.2>
$$
\bigl[\)\Jc(\la\);z)\>,\)x\ox 1+1\ox x\)\bigr]\,=\,0\,.
\Tag{Jch}
$$
\par
Denote by $\Dl:\Usl\)\to\)\Usl\ox\Usl$ the coproduct for $\Usl$.
\Th{Jccl}
The series $\Jc(\la\);z)$ satisfies the \eq/
$$
\Jc\"{(12)3}(\la\);z)\>\Jc\"{12}(\la-z\)h\"3;z)\,=\,
\Jc\"{1(23)}(\la\);z)\>\Jc\"{23}(\la\);z)\,.
\Tag{2ccl}
$$
Here \>$\Jc\"{(12)3}=(\Dl\ox\id\))(\Jc)$, \>$\Jc\"{1(23)}=(\)\id\ox\Dl)(\Jc)$,
\vv.06>
\>$\Jc\"{12}=\Jc\ox 1$, \>$\Jc\"{23}=1\ox \Jc$, and the meaning of
$\Jc\"{12}(\la-z\)h\"3;z)$ is explained below, \cf. \(lazh)\).
\endpro
\Pf.
The statement is a degeneration of Theorem~3.1 in \cite{ESS}\),
and can be proved in the same way.
\epf
\Rem
Let $x_1\lc x_{N-1}$ be a basis of $\hga\}$, and let $x^1\lc x^{\)N-1}$ be
the dual basis of $\hg$. Write $\la=\la^1x_1\]\lsym+\)\la^{N-1}x_{N-1}$. 
For a \raf/ $f(\la)$ we define a series $f(\la-z\)h\"3)$ by the Taylor
series expansion:
\vvnn-.5:-.8>
$$
f(\la-z\)h\"3)\,=\,f(\la\);z)\,-\,z\>
\sum_{a=1}^{N-1}\,{\der f(\la)\over\der\la^a}\;(x^a)\"3\)+\ldots{}\,,
\vv-.2>
\Tag{lazh}
$$
and extend the definition to series in $z$ with coefiicients in \raf/s of $\la$
in the natural way.
\wwmgood-.8:.8>
\enddemo
\Rem
The \eq/ \(2ccl) is usually called the \em{dynamical \$2$-cocycle condition}.
\enddemo
Define a \em{\rat/ exchange matrix} $\Rc(\la\);z)$ by the rule:
\vvnn-.2:-.1>
$$
\Rc(\la\);z)\,=\,\bigl(\Jc(\la\);z)\bigr)\vpb{-1}\Jc\"{21}(\la\);z)\,.
\Tag{Rc}
$$
\Th{Rc3}
$\Rc(\la)$ satisfies the \DYB/:
\ifMag
$$
\align
\Rc\"{12}(\la-z\)h\"3;z)\> &\Rc\"{13}(\la\);z)\>\Rc\"{23}(\la-z\)h\"1;z)\,={}
\\
\nn6>
{}={}\, &\Rc\"{23}(\la\);z)\>\Rc\"{13}(\la-z\)h\"2;z)\>\Rc\"{12}(\la\);z)\,.
\endalign
$$
\else
$$
\Rc\"{12}(\la-z\)h\"3;z)\>\Rc\"{13}(\la\);z)\>\Rc\"{23}(\la-z\)h\"1;z)\,=
\,\Rc\"{23}(\la\);z)\>\Rc\"{13}(\la-z\)h\"2;z)\>\Rc\"{12}(\la\);z)\,.
\vv-.3>
$$
\fi
\endpro
\nt
The statement follows from Theorem~\[Jccl] and cocommutativity of
the coproduct $\Dl$.
\Par
Say that an \$n$-tuple $V_1\lc V_n$ of \smod/s is admissible if for any \pd/
$i_1\lc i_k$ the tensor product $V_{i_1}\!\lox V_{i_k}$ is a \dhmod/; that is,
$V_1\lc V_n$ are \dhmod/s and all \wt/ subspaces of any tensor product
$V_{i_1}\!\lox V_{i_k}$ are \fd/.
\vsk.3>
Let $V,\)W\}$ be an admissible pair of \smod/s.
\vv.04>
Denote by $J_{VW}\:(\la\);z)\in\End(V\]\ox W)$ the action of $\Jc(\la\);z)$ in
the tensor product $V\]\ox W\}$. It follows from the explicit form of recurrence
relations in the proof of Proposition~\[abrr] that there is a unique \fn/
$J_{VW}\:(\la)\in\Rend{V\ox W}$ \st/ the series $J_{VW}\:(\la\);z)$ coincides
with the expansion of $J_{VW}\:(\la/z)$ at ${z=0}$. The \fn/ $J_{VW}\:(\la)$
admits the following description.
\Lm{J}
$J_{VW}\:(\la)$ is the unique \sol/ of the \eq/
$$
J_{VW}\:(\la)\,\bigl(1\ox(\la-\Xi\))\bigr)\vst{V\ox\)W}\)=\,
\bigl(1\ox(\la-\Xi\))\>+\!
\tsum_{\abn}\!e_{ab}\ox e_{ba}\bigr)\vst{V\ox\)W}\)J_{VW}\:(\la)
\vvmm-.3:-.4>
$$
\st/ $\bigl(J_{VW}\:(\la)-1\)\bigr)\in
\bigl(U(\)\bgp)\>\ngp\]\ox\)\ngm U(\)\bgm)\bigr)\vst{V\ox\)W}$. Moreover,
$J_{VW}\:(\la)$ commutes with the action of $\hg$ in $V\]\ox W\}$,
\cf. \(Jch)\).
\endpro
\Prop{JUVW}
For any admissible triple of \smod/s $U,\)V,\)W\}$ we have
$$
J_{U\ox V,W}(\la)\>\bigl(J_{UV}(\la-h\"3)\ox 1\)\bigr)\,=\,
J_{U,V\ox W}(\la)\>\bigl(1\ox J_{VW}(\la)\bigr)\,.
\vvgood
\vv-.5>
\Tag{J4}
$$
\endpro
\Pf.
The \fn/ $\bigl(J_{UV}(\la-h\"3)\ox 1\)\bigr)$ is defined by the rule \(fh).
To get relation \(J4) from formula \(2ccl) one needs to verify that the series
obtained by expansion of $\bigl(J_{UV}(z\1\la-h\"3)\ox 1\)\bigr)$ at ${z=0}$
coincides with the action of $\Jc\"{12}(\la-z\)h\"3;z)$, defined by \(lazh)\),
in $U\]\ox V\]\ox W\}$, which is simple.
\epf
For any $A\in\End(W\]\ox V)$ let $A\"{21}\]=\)P A\,P\1\!\in\End(V\]\ox W)$
where $P\]:W\]\ox V\]\to\)V\]\ox W\}$ is the \perm/ map: $P(x\ox y)=y\ox x$.
Set
\vvnn-.2:-.3>
$$
R_{VW}\:(\la)\)=\)
\bigl(J_{VW}\:(\la)\bigr)\vpb{-1}\bigl(J_{WV}\:(\la)\bigr)\"{21}\,.
\Tag{RJ}
$$
It is clear that the action of the series $\Rc(\la\);z)$ in $V\]\ox W\}$
coincides with the expansion of $R_{VW}\:(\la/z)$ at ${z=0}$. Like in the proof
of Proposition~\[JUVW] we get the following assertion from Theorem~\[Rc3].
\vskmgood-.6:.6>
\Prop{RUVW}
For any admissible triple of \smod/s $U,\)V,\)W\}$ we have
$$
R_{UV}\"{12}(\la-h\"3)\>R_{UW}\"{13}(\la)\>R_{VW}\"{23}(\la-h\"1)\,=\,
R_{VW}\"{23}(\la)\>R_{UW}\"{13}(\la-h\"2)\>R_{UV}\"{12}(\la)\,.
$$
\endpro
Let $\Jti_{VW}\:(\la)$ be the fusion matrix for $\Usl$ defined in \Cite{EV2},
and let
\vvn-.1>
$$
\Rti_{VW}\:(\la)\)=\)
\bigl(\Jti_{VW}\:(\la)\bigr)\vpb{-1}\bigl(\Jti_{WV}\:(\la)\bigr)\"{21}
\vv-.5>
\Tag{RtVW}
$$
be the corresponding dynamical \Rm/.
Let \,$X\)=\suan(-1)^{a\)-1}E_{a,N-\)a\)+1}$, considered as
an element of the group $SL(N)$, and let $w_X\:$ be the longest element
\vv.06>
of the Weyl group. We have $\Ad_X\:(e_{ab})=e_{N-\)a\)+1,\)N-\)b\)+1}\:$
\)and \>$w_X\:(\eps_a)=\eps_{N-\)a\)+1}\:$ for any $a\),\bno$.
\Lm{JJ}
For any \fd/ \smod/s $V,\)W\}$ we have
\vvn-.1>
$$
\Jti_{VW}\:(\la)\,=\,
(X\]\ox\]X)\>J_{VW}\:\bigl(w_X\:(\la+\rho\))\bigr)\)(X\]\ox\]X)\1\).
\vv-.6>
\Tag{JJt}
$$
\endpro
\Pf.
It is shown in \)\Cite{\)EV2\), Section~9\)} \>that
$\Jti_{VW}\:(\la)$ is the only \sol/ of the \eq/
$$
\Jti_{VW}\:(\la)\,\bigl(1\ox(\la+\rho-\Xi\))\bigr)\vst{V\ox\)W}\)=\,
\bigl(1\ox(\la+\rho-\Xi\))\>+\!
\tsum_{\abn}\!e_{ba}\ox e_{ab}\bigr)\vst{V\ox\)W}\)\Jti_{VW}\:(\la)
\ifMag\kern-.7em\fi
\vv-.4>
$$
\st/ \)${\bigl(\Jti_{VW}\:(\la)-1\)\bigr)\in
\bigl(\ngm U(\)\bgm)\ox U(\)\bgp)\>\ngp\bigr)\vst{V\ox\)W}}$. By Lemma~\[J]
\vv.1>
\rhs/ of formula \(JJt) has the same properties, which proves the claim.
\epf
\vsk-.3>
\vsk0>
\Cr{RRt}
$\Rti_{VW}\:(\la)\,=\,
(X\]\ox\]X)\>R_{VW}\:\bigl(w_X\:(\la+\rho\))\bigr)\)(X\]\ox\]X)\1\}$.
\endpro
\vsk.3>
Henceforward, let $U\]$ be the vector representation of $\sln$.
By formula (36) in \Cite{EV2} we have
\vvnn-.6:-.1>
$$
\Rti_{UU}\:(\la)\,=\sum_{a,b=1}^N\}E_{aa}\ox E_{bb}\,
-\!\}\sum_{\abn}{E_{bb}\ox E_{aa}\over(\la_{ab}\]-a+b\))^2}\;
-\]\sum_{\tsize{a,b=1\atop a\ne b}}^N\]{E_{ab}\ox E_{ba}\over
\la_{ab}\]-a+b}\;.\ifMag\kern-1em\else\kern-1.4em\fi
\vv-1.3>
\Tag{Rti}
$$
Therefore,
\vvn-.7>
$$
R_{UU}\:(\la)\,=\sum_{a,b=1}^N\}E_{aa}\ox E_{bb}\,
-\!\}\sum_{\abn}{E_{aa}\ox E_{bb}\over\la_{ab}^2}\;
-\]\sum_{\tsize{a,b=1\atop a\ne b}}^N\]{E_{ab}\ox E_{ba}\over
\la_{ab}}\;.\ifMag\kern-1em\else\kern-1.4em\fi
\vvmm-.5:-.3>
\Tag{RUU}
$$
Let $V\}$ be an \asmod/. Introduce \fn/s \>$\elhV_{ab}\in\Rend V$,
\>$a\),\bno$, by the equality
\vvnn-1:-.2>
$$
R_{UV}\:(\la)\,=\sum_{a,b=1}^N\}E_{ba}\ox \elhV_{ab}(\la)\,.
\mmgood
$$
\vskm-.5:-.6>
\vsk0>
\Th{Ec}
Let $V\}$ be an \asmod/. Then the rule
\ifMag
$$
(\tth_{ab}\,v)\)(\la)\,=\,\elhV_{ab}(\la)\>v(\la-\eps_a)\,,
$$
\else
\,${(\tth_{ab}\,v)\)(\la)\)=\)\elhV_{ab}(\la)\>v(\la-\eps_a)}$, \>
\fi
for any \)$a\),\bno$ and any \)$v\in\RatV$, endows $V\}$ with an \ermod/
structure. The constructed \ermod/ is denoted by $\Ec(V)$.
\endpro
\Pf.
The statement follows from Proposition~\[RUVW]\), and formulae \(tthr)\),
\(RTTb)\), \(ellh) and \(RUU)\).
\vvgood
\epf
We define the functor $\Ec$ from the category of \asmod/s to the category
of semistandard \ermod/s by sending an object $V\}$ to $\Ec(V)$ and a morphism
$\phi\in\Hom(V,W)$ to the corresponding constant \fn/ $\phi\in\Mor(V,W)$.
\Pf of Theorem~\[EcV]\).
Claim \)a) of the theorem is immediate.
Say that $f(\la)=O\bigl(\)|\)\la\)|^k\bigr)$ if $f(s\)\la)=O(s^{\)k})$ as
$s\to+\8$ for generic $\la$. From Lemma \[J] and formula \(RJ) we have that
\vvnn-.8:-.2>
$$
\gather
J_{UV}\:(\la)\,=\,1\,-\!\sum_{\abn}\!\!{E_{ba}\ox e_{ab}\over\la_{ab}}\,\)+\,
O\bigl(\)|\)\la\)|^{-2}\bigr)\,,
\\
\nn8>
R_{UV}\:(\la)\,=\,1\,+\sum_{\tsize{a,b=1\atop a\ne b}}^N
{E_{ba}\ox e_{ab}\over\la_{ab}}\,\)+\,O\bigl(\)|\)\la\)|^{-2}\bigr)\,,
\\
\cnn-.2>
\endgather
$$
which proves claim \)b)\). Claim \)c) follows from Lemma~\[ptndg]\).
Claims \)e) and \)f) follow from claims \)b) and \)d)\), \>and Lemma~\[CcV]\).
\vsk.2>
To prove claim \)d) \)one should show that for any \singv/ ${v\in V}$
we have \>${\elhV_{aa}(\la)\,v\)=\)v}$ \>for any $a$,
\>and \>${\elhV_{ab}(\la)\,v\)=0\)}$ \>for any ${a<b}$.
By Lemma~\[J] and formula \(RJ) we see that
\>$\bigl(\elhV_{aa}(\la)-1\)\bigr)\in\bigl(\ngm\)\Usl\>\ngp\bigr)\vst V$
\>for any $a$, \>and \>$\elhV_{ab}(\la)\in\bigl(\Usl\>\ngp\bigr)\vst V$
for any $a<b$, which implies claim \)d)\).
\vsk.3>
It remains to prove claim \)c\). The element $D$ acts on $\RatV$ as
multiplication by $D(\la)$. It is clear from the definition of $\Ec(V)$ that
there exists some independent of $V\}$ element in a certain completion of
$\Usl$ \st/ its action on $V\}$ coincides with $D(\la)$. Since an element of
$\Usl$ is uniquely determined by its action in \hwm/s, it suffices to prove
claim \)c) \)under the assumption that $V\}$ is a \hwm/. In the last case
claim \)c) \)follows from claim \)d)\).
\epf
The proof of Theorem~\[EcUV] is similar to the proof of Theorem~45
in \Cite{EV2}\).

\Sect[ex]{Exchange quantum group $\{\Fsl$}
For any $\ano$ \>let \>$\ib^a\]=\)(1\lc a-1,\)a+1\lc N)$. Set
$$
\tsize\tb_{ab}\,=\,(-1)^{a+b}\}\tsum_{\jb\)\in\)\Sb_{N-1}\!\!\!}\sign(\)\jb\))
\;t_{i^a_{N-1},\>i^b_{j_{N-1}}}\>\ldots\,t_{i^a_1,\>i^b_{j_1}}\,.
$$
\vsk-1.2>
\vsk0>
\Lm{ttt}
$\sum_{c=1}^N\tb_{ac}\,t_{bc}\,=\,\dl_{ab}\,\tN\}$, \,where $\tN$ is defined
by \(tN)\).
\vvm-.2>\vvm0>
\endpro
\Pf.
The formula coincides with the equality of the top coefficients in the \rat/
version of formula \(TT=T) for $k=1$ and \)$T(u)=\Tc(u)$, \cf. \(Tcu)\).
\mmgood
\epf
Consider the exchange quantum group $\Fsl$ defined in \Cite{EV2}\).
It admits the following description, see \)\Cite{\)EV2\), Section~5.3\)}\).
$\Fsl$ is a unital associative algebra over $\C$ \gby/ \fn/s $f\]\in\Ratt$
and elements $\Tps_{ab}$, $\Tms_{ab}\}$, \,$a\),\bno$, \,subject to relations
\(LL)\,--\,\(RLL) and \(det=1)\).
\vsk.2>
Let ${\Rti(\la)=\Rti_{UU}\:(\la)}$, \cf. \(Rti)\).
Set \,${\Tpms\!=\sum_{a,b=1}^N E_{ab}\ox\Tpms_{ab}}$.
\vvm-.5>
\,The defining relations for $\Fsl$ are
\vvn-.6>
$$
\gather
\Tps\>\Tms=\,\Tms\>\Tps=\,\id\)\ox 1\,,
\Tag{LL}
\\
\nn8>
\Tps_{ab}\,f(\la\+1\},\la\+2)\,=
\>f(\la\+1\!-\eps_b\>,\la\+2\!-\eps_a)\,\Tps_{ab}
\Tag{Lf}
\\
\cnn-.3>
\endgather
$$
for any $f\in\Ratt$,
\vvn-.3>
$$
\Rti\"{12}(\la\+2)\>T\"{13}\)T\"{23}\,=\,T\"{23}\)T\"{13}\>\Rti\"{12}(\la\+1)
\kern-1em
\Tag{RLL}
$$
where \,$T\"{13}=\sum_{a,b}\)E_{ab}\ox\id\ox\Tps_{ab}$ \,and
\,$T\"{23}=\botsmash{\sum_{a,b}\>\id\ox E_{ab}\ox\Tps_{ab}}$,
\,and relation \(det=1) below.
\Rem
In this paper the \var/s $\la\+1,\)\la\+2$ and the generators
$\Tps_{ab}\),\>\Tms_{ab}$ correspond to the \var/s $\la^2,\>\la^1$ and
the generators $L_{ab}\),\>L\1_{ab}$ in \Cite{EV2}\).
\wwgood-.8:.8>
\enddemo
For any \perm/ $\ib\in\Sb_N$ let \,$\dsize\La_{\)\ib}(\la)\>=\!\!
\prod_{\tsize{\abn\atop i_a<i_b}}\!(1+\la_{ab}\1)$. \,Set
\vvn-.1>
$$
\Det\)\Tps=\;{1\over\La_{\)\id}(\la\+1)}\,\sumib\sign(\)\ib\))\,
\La_{\)\ib}(\la\+2)\;\Tps_{i_N\],\)N}\)\ldots\,\Tps_{i_1,\)1}\kern-1em
\Tag{detL}
$$
where $\id=(1\lc N)$. The last defining relation for $\Fsl$ is
$$
\Det\)\Tps=\,1\,.
\vvm-.06>
\Tag{det=1}
$$
\Rem
The element $\Det\)\Tps\!$ corresponds to the element $D$ in \Cite{EV2}\).
Formula \(detL) can be derived from the definition of $D$ therein.
The complete proof of formula \(detL) will appear elsewhere.
\enddemo
Recall that, given a \dhmod/ $V\}$, we assume the following action of $\Ratt$
on \Vval/ \fn/s:
$$
f\):\>v(\la)\,\map\,f(\la\>,\la-\mu)\,v(\la)
\vvmm0:.4>
$$
for any $f\in\Ratt$ and any \fn/ $v(\la)$ with values in $\Vmu$.
\vsk.3>
A \rat/ dynamical \rep/ of $\Fsl$ is a \dhmod/ $V\}$ endowed with an action
of $\Fsl$ on \Vval/ \mef/s by \dif/ operators with \rat/ coefficients:
\vvnn0:-.3>
$$
(\Tps_{ab}\,v)\)(\la)\,=\,L_{ab}(\la)\>v(\la-\eps_b)\,,\qqq\Rlap{a\),\bno\,,}
\vv.1>
\Tag{TpL}
$$
where ${L_{ab}(\la)\in\Rend V}$ are suitable \fn/s.
\Prop{L2t}
Let $V\}$ be a \rat/ dynamical \rep/ of $\Fsl$. Then the rule
\vvn.2>
$$
(t_{ab}\,v)\)(\la)\,=\,\prod_{1\le c<a}\](\la_{ca}\]+1\))\vpb{-1}\!
\prod_{1\le c<b}\](\la_{cb}\]+e_{bb}\]-e_{cc}\]+1\))
\,L_{ba}(\la-\rho)\>v(\la-\eps_a)
\Tag{tL}
$$
for any \)$a\),\bno$ and any $v\in\RatV$, endows $V\}$ with a structure
of a semi\-standard \ermod/.
\endpro
\nt
The proof is straightforward.
\Prop{t2L}
Let $V\}$ be a semistandard \ermod/. Then formulae \(TpL)\), \(tL) make $V\}$
into a \rat/ dynamical \rep/ of $\Fsl$.
\endpro
\Pf.
It is straightforward to verify relations \(Lf)\), \(RLL) and \(det=1).
To complete the proof it remains to find the action of elements $\Tms_{ab}$
to obey relations \(LL). This can be done using Proposition~\[ttt]\), \)since
the \ermod/ $V\}$ is \ndeg/.
\epf
The last two propositions define a functor $\Fc$ from the category of \rat/
dynamical \rep/s of $\Fsl$ to the category of semistandard \ermod/s:
an object $V\}$ goes to itself and a morphism
$\phi(\la)\in\Rat\bigl(\Hom(V,W)\bigr)$ goes
to $\phi(\la-\rho\))\in\Mor(V,W)$. Furthermore,
the propositions imply the following assertion.
\Th{cat}
The functor $\Fc$ is an equivalence of the categories.
\endpro
For both categories involved in the last theorem the subcategories of \fdim/
objects are tensor categories, the tensor product of \rat/ dynamical \rep/s of
$\Fsl$ being defined in \Cite{EV2}\). One can show that the restriction of the
functor $\Fc$ to these subcategories is an equivalence of tensor categories.
\vsk.3>
Let $\Gc$ be the functor from the category of \fd/ \smod/s to the category of
\fd/ dynamical \rep/s of $\Fsl$ defined in \Cite{EV2}\): an \smod/ $V\}$ goes
to the \rep/ $\Gc(V)$ of $\Fsl$ given by the rule
\vvn-.1>
$$
\Rti_{UV}\:(\la)\,=\sum_{a,b=1}^N\!E_{ab}\ox L_{ab}(\la)\,,\kern-1em
\Tag{Gc}
$$
and a morphism $\phi\in\Hom(V,W)$ goes to the corresponding constant \fn/
\vv-.04>
$\phi\in\Rat\bigl(\Hom(V,W)\bigr)$. The composition $\Ect=\Fc\}\o\Gc$ is
a functor from the category of \fd/ \smod/s to the category of semistandard
\ermod/s.
\vskgood-.8:.8>
\Th{Ect}
The functor $\Ect$ is isomorphic to the restriction of the functor $\Ec$
to the category of \fd/ \smod/s.
\endpro
\nt
The theorem is proved in Appendix~\[:B]\).
\Rem
Let $V\}$ be an \irr/ \fd/ \smod/. Then the \ermod/ $\Ec(V)$ is an \irr/
standard \hwm/ over the dynamical \qg/ $\eslr$. Such modules have been
described in Section~\[:Fd]\). Applying the functor inverse to $\Fc$ we get
a new description of the dynamical \rep/s of $\Fsl$ induced from \irr/ \fd/
\smod/s. This new description is a new \hwm/ theory for the exchange dynamical
\qg/ $\Fsl$.
\mmgood
\enddemo

\Appendix
\Sect[CR]{Commutation relations in $\{\els$}
\nt
In this section we collect useful commutation relations which hold in the
\oalg/ $\els$.
\vsk.4>
In the definition of $\els$ one can replace relations \(tac)
by the following relations:
\vvm-.6>
$$
{\tht(\la\+1_{ac}\!-\gm)\over\tht(\la\+1_{ac})}\ t_{ab}\,t_{cd}\,-\,\)
{\tht(\la\+2_{bd}\!-\gm)\over\tht(\la\+2_{bd})}\ t_{cd}\,t_{ab}\,\)=\;
{\tht(\la\+1_{ac}\!+\la\+2_{bd})\,\tht(\gm)\over
\tht(\la\+1_{ac})\,\tht(\la\+2_{bd})}\ t_{cb}\,t_{ad}\,,
\ifMag\kern-1em\else\kern-1.4em\fi
\vv.2>
\Tag{tbd}
$$
for \;$a\ne c$ \;and \;$b\ne d$. \,The last formula implies that
$$
\tth_{aa}\,\tth_{bb}\>-\,\tth_{bb}\,\tth_{aa}\,=\;
{\tht(\la\+1_{ab}\!+\la\+2_{ab})\,\tht(\gm)\over
\tht(\la\+1_{ab})\,\tht(\la\+2_{ab})}\ \tth_{ba}\,\tth_{ab}
\ifMag\kern-1em\else\kern-1.4em\fi
\Tag{thab}
$$
for $a<b$. Under the same assumption we have
\ifMag
$$
\align
{\tht(\la\+2_{ab}\!-\gm)\,\tht(\la\+2_{ab}\!+\gm)
\over\tht(\la\+2_{ab})^2}\ \tth_{aa}\,\tth_{bb}\,-\,\)
{\tht(\la\+1_{ab}\!-\gm)\,\tht(\la\+1_{ab}\!+\gm)
\over\tht(\la\+1_{ab})^2}\ \tth_{bb}\,\tth_{aa}\,={}\!\}&
\\
\nn10>
{}=\;{\tht(\la\+1_{ab}\!+\la\+2_{ab})\,\tht(\gm)\over
\tht(\la\+1_{ab})\,\tht(\la\+2_{ab})}\ \tth_{ab}\,\tth_{ba} & \,,\!
\\
\cnn-.5>
\endalign
$$
\else
\vvn.2>
$$
{\tht(\la\+2_{ab}\!-\gm)\,\tht(\la\+2_{ab}\!+\gm)
\over\tht(\la\+2_{ab})^2}\ \tth_{aa}\,\tth_{bb}\,-\,\)
{\tht(\la\+1_{ab}\!-\gm)\,\tht(\la\+1_{ab}\!+\gm)
\over\tht(\la\+1_{ab})^2}\ \tth_{bb}\,\tth_{aa}\,=\;
{\tht(\la\+1_{ab}\!+\la\+2_{ab})\,\tht(\gm)\over
\tht(\la\+1_{ab})\,\tht(\la\+2_{ab})}\ \tth_{ab}\,\tth_{ba}\,,
\vv-.7>
$$
\fi
$$
\align
\tth_{ab}\,\tth_{ba}\,-\,\){\tht(\la\+1_{ab}\!-\gm)\,\tht(\la\+1_{ab}\!+\gm)
\over\tht(\la\+1_{ab})^2}\ \tth_{ba}\,\tth_{ab}\, &{}=\,-\,
{\tht(\la\+1_{ab}\!-\la\+2_{ab})\,\tht(\gm)\over
\tht(\la\+1_{ab})\,\tht(\la\+2_{ab})}\ \tth_{aa}\,\tth_{bb}\,,
\\
\ald
\nnm11:10>
\tth_{ab}\,\tth_{ba}\,-\,\){\tht(\la\+2_{ab}\!-\gm)\,\tht(\la\+2_{ab}\!+\gm)
\over\tht(\la\+2_{ab})^2}\ \tth_{ba}\,\tth_{ab}\, &{}=\,-\,
{\tht(\la\+1_{ab}\!-\la\+2_{ab})\,\tht(\gm)\over
\tht(\la\+1_{ab})\,\tht(\la\+2_{ab})}\ \tth_{bb}\,\tth_{aa}\,.
\\
\cnn.4>
\endalign
$$
More general relations are listed below.
We assume that $a<c$ and $b<d$ therein.
$$
\alignat2
\tth_{ab}\,\tth_{cb}^{\>k}\, &{}=\,\)
{\tht(\la\+1_{ac}\!+k\)\gm)\over\tht(\la\+1_{ac})}\ \tth_{cb}^{\>k}\>\tth_{ab}
\,, & \tth_{ab}^{\>k}\>\tth_{cb}\, &{}=\,\){\tht(\la\+1_{ac}\!+\gm)
\over\tht\bigl(\la\+1_{ac}\!-(k-1)\>\gm\bigr)}\ \tth_{cb}\,\tth_{ab}^{\>k}\,,
\ifMag\kern-3em\else\kern-1.8em\fi
\Tag{thac}
\\
\nn10>
\tth_{ad}\,\tth_{ab}^{\>k}\, &{}=\,\){\tht(\la\+2_{bd}\!-k\)\gm)\over
\tht(\la\+2_{bd})}\ \tth_{ab}^{\>k}\,\tth_{ad}\,, &
\ifMag\kern3em\else\kern4em\fi
\tth_{ad}^{\>k}\>\tth_{ab}\, &{}=\,\){\tht(\la\+2_{bd}\!-\gm)\over
\tht\bigl(\la\+2_{bd}\!+(k-1)\>\gm\bigr)}\ \tth_{ab}\>\tth_{ad}^{\>k}\,,
\ifMag\kern-3em\else\kern-1.8em\fi
\Tag{thbd}
\\
\cnn-.6>
\endalignat
$$
$$
\tth_{ab}\,\tth_{cd}\>-\,\tth_{cd}\,\tth_{ab}\,=\;
{\tht(\la\+1_{ac}\!+\la\+2_{bd})\,\tht(\gm)\over
\tht(\la\+1_{ac})\,\tht(\la\+2_{bd})}\ \tth_{cb}\,\tth_{ad}\,,\kern-1em
\vv-.4>
\Tag{thabcd}
$$
\ifMag
\mmgood
$$
\align
{\tht(\la\+2_{bd}\!+\gm)\,\tht(\la\+2_{bd}\!-\gm)
\over\bigl(\)\tht(\la\+2_{bd})\bigr)\vpb2}\ \tth_{ab}\,\tth_{cd}\,-\,\)
{\tht(\la\+1_{ac}\!+\gm)\,\tht(\la\+1_{ac}\!-\gm)
\over\bigl(\)\tht(\la\+1_{ac})\bigr)\vpb2}\ \tth_{cd}\,\tth_{ab}\,={}\,
\kern-1em &
\Tag{thcdab}
\\
\nn10>
{}=\;{\tht(\la\+1_{ac}\!+\la\+2_{bd})\,\tht(\gm)\over
\tht(\la\+1_{ac})\,\tht(\la\+2_{bd})}\ \tth_{ad}\,\tth_{cb}\,,\kern-1em &
\endalign
$$
\else
$$
{\tht(\la\+2_{bd}\!+\gm)\,\tht(\la\+2_{bd}\!-\gm)
\over\bigl(\)\tht(\la\+2_{bd})\bigr)\vpb2}\ \tth_{ab}\,\tth_{cd}\,-\,\)
{\tht(\la\+1_{ac}\!+\gm)\,\tht(\la\+1_{ac}\!-\gm)
\over\bigl(\)\tht(\la\+1_{ac})\bigr)\vpb2}\ \tth_{cd}\,\tth_{ab}\,=\;
{\tht(\la\+1_{ac}\!+\la\+2_{bd})\,\tht(\gm)\over
\tht(\la\+1_{ac})\,\tht(\la\+2_{bd})}\ \tth_{ad}\,\tth_{cb}\,,\kern-2em
\vv-.4>
\Tag{thcdab}
$$
\fi
$$
\align
\tth_{ad}\,\tth_{cb}\,-\,\){\tht(\la\+1_{ac}\!+\gm)\,\tht(\la\+1_{ac}\!-\gm)
\over\bigl(\)\tht(\la\+1_{ac})\bigr)\vpb2}\ \tth_{cb}\,\tth_{ad}\,=\,-\,
{\tht\bigl(\la\+1_{ac}\!-\la\+2_{bd})\,\tht(\gm)\over
\tht(\la\+1_{ac})\,\tht(\la\+2_{bd})}\ \tth_{ab}\,\tth_{cd}\,,\kern-2em
\Tagg{thadcb}
\\
\nn12>
\tth_{ad}\,\tth_{cb}\,-\,\){\tht(\la\+2_{bd}\!+\gm)\,\tht(\la\+2_{bd}\!-\gm)
\over\bigl(\)\tht(\la\+2_{bd})\bigr)\vpb2}\ \tth_{cb}\,\tth_{ad}\,=\,-\,
{\tht\bigl(\la\+1_{ac}\!-\la\+2_{bd})\,\tht(\gm)\over
\tht(\la\+1_{ac})\,\tht(\la\+2_{bd})}\ \tth_{cd}\,\tth_{ab}\,.\kern-2em
\\
\nngood
\cnn.25>
\endalign
$$
All these formulae follow from \(tf)\,--\,\(tac) and \(tbd)\),
and the summation formulae for the theta \fn/.
Combining formula \(thabcd) with formulae \(thac) and \(thbd)
for $k=1$ we can obtain the following Serre\)-type relations:
\vvmm0:-.4>
$$
\alignat2
\tht(\la\+1_{ac}\!+\la\+2_{bd})\,\tht(\la\+1_{ac}\!-\gm)\,
\tht(\la\+2_{bd}\!-2\)\gm)\,\tht(\gm)\;\tth_{ab}^{\>2}\, &
\Rlap{\tth_{cd}\,-{}} &&
\Tag{Sab}
\\
\nn6>
{}-\,\tht(\la\+1_{ac}\!+\la\+2_{bd}\!-\gm)\,\tht(\la\+1_{ac})\,
\tht(\la\+2_{bd}\!-\gm)\,\tht(2\)\gm)\; & \tth_{ab}\, &&
\tth_{cd}\,\tth_{ab}\,+{}
\\
\nn6>
&& \Llap{{}+\,\tht(\la\+1_{ac}\!+\la\+2_{bd}\!-2\)\gm)\,\tht(\la\+1_{ac}\!+\gm)
\,\tht(\la\+2_{bd})\,\tht(\gm)\;} & \tth_{cd}\,\tth_{ab}^{\>2}\,=\,0\,,
\kern-1em
\\
\cnnm.5>
\cnn-.9>
\endalignat
$$
$$
\alignat2
\tht(\la\+1_{ac}\!+\la\+2_{bd}\!+2\)\gm)\,\tht(\la\+1_{ac})\,
\tht(\la\+2_{bd}\!-\gm)\,\tht(\gm)\;\tth_{ab}\, & \Rlap{\tth_{cd}^{\>2}\,-{}}&&
\Tag{Scd}
\\
\nn6>
{}-\,\tht(\la\+1_{ac}\!+\la\+2_{bd}\!+\gm)\,\tht(\la\+1_{ac}\!+\gm)\,
\tht(\la\+2_{bd})\,\tht(2\)\gm)\; & \tth_{cd}\, && \tth_{ab}\,\tth_{cd}\,+{}
\\
\nn6>
&& \Llap{{}+\,\tht(\la\+1_{ac}\!+\la\+2_{bd})\,\tht(\la\+1_{ac}\!+2\)\gm)\,
\tht(\la\+2_{bd}\!+\gm)\,\tht(\gm)\;} & \tth_{cd}^{\>2}\,\tth_{ab}\,=\,0\,.
\kern-1em
\endalignat
$$
\vsk0>
\Lm{rels}
Let $a<c$ and $b<d$. Then the following relations hold\/\){\rm:}
\vvnn.4:.2>
$$
\align
\tth_{ab}\,\tth_{cd}^{\>k}\>-\,\tth_{cd}^{\>k}\,\tth_{ab}\, &{}=\;
{\tht\bigl(\la\+1_{ac}\!+\la\+2_{bd}\!+(k-1)\>\gm\bigr)\,\tht(k\)\gm)\over
\tht(\la\+1_{ac})\,\tht\bigl(\la\+2_{bd}\!+(k-1)\>\gm\bigr)}
\ \tth_{cb}\,\tth_{cd}^{\>k-1}\)\tth_{ad}\,,
\ifMag\kern-2em\else\kern-1em\fi
\\
\nn12>
\tth_{ab}^{\>k}\,\tth_{cd}\>-\,\tth_{cd}\,\tth_{ab}^{\>k}\, &{}=\;
{\tht\bigl(\la\+1_{ac}\!+\la\+2_{bd}\!-(k-1)\>\gm\bigr)\,\tht(k\)\gm)\over
\tht\bigl(\la\+1_{ac}\!-(k-1)\>\gm\bigr)\,\tht(\la\+2_{bd})}
\ \tth_{cb}\,\tth_{ab}^{\>k-1}\)\tth_{ad}\,,
\ifMag\kern-2em\else\kern-1em\fi
\Tagg{thabcdk}
\\
\cnnm.2>
\cnn-.35>
\endalign
$$
$$
\align
\tth_{ad}\,\tth_{cb}^{\>k}\,-\,\){\tht(\la\+1_{ac}\!+k\)\gm)\,
\tht(\la\+1_{ac}\!-\gm)\,\tht\bigl(\la\+2_{bd}\!-(k-1)\>\gm\bigr)\over
\bigl(\)\tht(\la\+1_{ac})\bigr)\vpb2\>\tht(\la\+2_{bd})}
\ \tth_{cb}^{\>k}\,\tth_{ad}\,={} &
\Tagg{thadcbk}
\\
\nnm10:12>
{}=\,-\,{\tht\bigl(\la\+1_{ac}\!-\la\+2_{bd}\!+(k-1)\>\gm\bigr)\,\tht(k\)\gm)
\over\tht(\la\+1_{ac})\,\tht(\la\+2_{bd})}
\,\) & \]\!\ \tth_{cb}^{\>k-1}\)\tth_{ab}\,\tth_{cd}\,.
\ifMag\kern-3em\else\kern-2em\fi
\\
\cnn-.3>
\endalign
$$
\endpro
\Pf.
For $k=1$ these formulae coincide with formulae \(thabcd) and \(thadcb),
\resp/. All the proofs for $k>1$ are similar to each other.
So we will prove only formula \(thadcbk).
\vsk.3>
Multiply formula \(thcdab) by $\tth_{cb}^{\>k-1}\}$ from the right, and push
\vv.06>
the factor $\tth_{cb}^{\>k-1}\}$ in \lhs/ through all the products from right
to left using relations \(thac). Taking into account formula \(thabcd) we get
\vvn.1>
$$
\align
\tth_{ad}\,\tth_{cb}^{\>k}\,-\,\){\tht(\la\+1_{ac}\!+k\)\gm)\,
\tht(\la\+1_{ac}\!-\gm)\,\tht\bigl(\la\+2_{bd}\!-(k-1)\>\gm\bigr)\over
\bigl(\)\tht(\la\+1_{ac})\bigr)\vpb2\>\tht(\la\+2_{bd})}
\ \tth_{cb}^{\>k}\,\tth_{ad}\,={} &
\\
\nnm8:10>
{}=\bigl(F(\la\+1_{ac},\la\+2_{bd})+F(\)-\>\la\+2_{bd},-\>\la\+1_{ac})\bigr)\,
& \>\tth_{cb}^{\>k-1}\>\tth_{ab}\,\tth_{cd}\,\kern-1em
\\
\cnn-.2>
\endalign
$$
where
$$
F(u\),v)\,=\;{\tht\bigl(u+(k-1)\>\gm\bigr)\,\tht(v+\gm)\,\tht(v-k\)\gm)\over
\tht(u+v)\,\tht(v)\,\tht(\gm)}\;.
\vv.4>
$$
$\=F(u,v)+F(-v,-u)\,$ is a quasiperiodic \fn/ of $u$ with periods $1$ and
$\tau$, and it has only simple poles. Thus, it is completely determined by
its multipliers and residues. Hence,
$$
F(u,v)+F(-v,-u)\,=\,-\,
{\tht\bigl(u-v+(k-1)\>\gm\bigr)\,\tht(k\)\gm)\over\tht(u)\,\tht(v)}\;.
\vv-1.3>
$$
\epf

\Sect[Det]{Quantum determinant}
The construction of $\Det\)T(u)$ and the proof of Proposition~\[qdet]
can be obtained by extending the standard fusion procedure technique
to the dynamical case.
\vsk.3>
For any $k=2\lc N$ \)let $A_k$ be the complete antisymmetrizer
in $(\CN)\vpb{\ox k}$:
\vvn-.2>
$$
A_k(x_1\lox x_k)\,=\;{1\over k\)!}\,
\sumibk\)\sign(\)\ib\))\;(x_{i_1}\lox x_{i_k})\,.
\vv-.2>
$$
Set $A\)=A_2$ and let $S=1-A$ be the corresponding symmetrizer.
The \Rm/ $R(u\),\la)$ has a simple pole at $u=\gm$.
Denote by$Q(\la)$ the resique of $R(u\),\la)$ at $u=\gm$.
\Lm{kerQ}
$\Ker Q(\la)=\Ker A=\Im S$.
\endpro
Due to the inversion relation \(inv) we can write relation \(RTT)
in the following form:
$$
T\"{13}(u)\>T\"{23}(v)\>R\"{21}(u-v,\la)\,=\,
R\"{21}(u-v,\la-\gm\)h\"3)\>T\"{23}(v)\>T\"{13}(u)\,.\kern-1em
\vv.1>
\Tag{TTR}
$$
Then by Lemma~\[kerQ] we have that
\,$A\"{12}\>T\"{23}(u+\gm)\>T\"{13}(u)\,S\"{12}\)=\,0$, \>which is \eqv/
to each of the following relations:
\vvnn0:-.2>
$$
\align
T\"{23}(u+\gm)\>T\"{13}(u)\,S\"{12}\) &{}=\,
S\"{12}\>T\"{23}(u+\gm)\>T\"{13}(u)\,S\"{12}\>,
\Tag{TTS}
\\
\nn6>
A\"{12}\>T\"{23}(u+\gm)\>T\"{13}(u)\, &{}=\,
A\"{12}\>T\"{23}(u+\gm)\>T\"{13}(u)\,A\"{12}\>.
\endalign
$$
Formula \(TTS) shows that for any $i=1\lc k-1$ we have
\vvn.1>
$$
\align
T\"{k,\)k+1}\bigl(u+(k-1)\>\gm\bigr)\,\ldots\,T\"{1,\)k+1}(u)\, &
S\"{i,\)i+1}\)={}
\\
\nn6>
{}=\,S\"{i,\)i+1}\>T\"{k,\)k+1}\bigl(u+(k-1)\>\gm\bigr)\,\ldots\,{} &
T\"{1,\)k+1}(u)\,S\"{i,i+1}\>,\kern-1em
\\
\cnn-.5>
\endalign
$$
which implies
\vvn-.3>
$$
\align
A_k\"{1\lc k}\, &
T\"{k,\)k+1}\bigl(u+(k-1)\>\gm\bigr)\,\ldots\,T\"{1,\)k+1}(u)\,={}
\Tag{ATT}
\\
\nn6>
{}=\,{} & A_k\"{1\lc k}\>T\"{k,\)k+1}\bigl(u+(k-1)\>\gm\bigr)\,\ldots\,
T\"{1,\)k+1}(u)\,A_k\"{1\lc k}\kern-1em
\\
\cnn-.1>
\endalign
$$
because \>$\Ker A_k\)=\)\sum_{i=1}^{k-1}\)\Im S\"{i,i+1}$ \>and
\>$(\CN)\vpb{\ox k}\}=\>\Im A_k\oplus\)\Ker A_k$.
\vvm.15>
We denote the restriction of
\ifMag
$A_k\"{1\lc k}\,T\"{k,\)k+1}\bigl(u+(k-1)\>\gm\bigr)\,\ldots\,T\"{1,\)k+1}(u)$
\vvm.07>
\else
\vvn-.5>
$$
A_k\"{1\lc k}\,T\"{k,\)k+1}\bigl(u+(k-1)\>\gm\bigr)\,\ldots\,T\"{1,\)k+1}(u)
\vv-.1>
$$
\fi
to \)$\Im A_k\ox\Fun(V)$ by $T^{\wedge k}(u)$ and call it
\vvm.1>
the \em{\$k$-th exterior power} of $T(u)$. Since $\Im A_N$ is \onedim/,
the top exterior power $\TN(u)$ can be considered as an element of
$\End\bigl(\FunV\bigr)$.
\Lm{TNl}
For any \perm/ $\jb\in\Sb_N$ we have
$$
\TN(u)\,=\,
\sign(\)\jb\))\sumib\sign(\)\ib\))\;T_{i_N\],\)j_N}\bigl(u+(N-1)\>\gm\bigr)
\>\ldots\,T_{i_2,\)j_2}(u+\gm)\,T_{i_1,\)j_1}(u)\,.
\vv-.8>
\Tag{TN}
$$
\endpro
\Pf.
Let $v_1\lc v_N$ be the standard basis of $\CN\}$.
Then for any $\jb\in\Sb_N$ we have
$$
A_N\>(\voxj)\,=\,\sign(\)\jb\))\,A_N\>(\vox)\,.
\Tag{Avj}
$$
Let $v\in\Fun(V)$. By the definition of $\TN(u)$ and relation \(ATT) we get
\vvn.2>
$$
\align
A_N\> & (\voxj)\ox \TN(u)\,v\,={}
\\
\nn7>
{}={}\, &\!\sumib A_N\>(\voxi)\ox\)T_{i_N\],\)j_N}\bigl(u+(N-1)\>\gm\bigr)
\>\ldots\,T_{i_2,\)j_2}(u+\gm)\,T_{i_1,\)j_1}(u)\,v\,.\kern-2em
\endalign
$$
By formula \(Avj) the expression in \rhs/ equals
\ifMag
$$
\align
& \]\sign(\)\jb\))\> A_N\>(\voxj)\,\ox{}
\\
\nn7>
& \,\,{}\ox\sumib\sign(\)\ib\))\;T_{i_N\],\)j_N}\bigl(u+(N-1)\>\gm\bigr)
\>\ldots\,T_{i_2,\)j_2}(u+\gm)\,T_{i_1,\)j_1}(u)\,v\,,\kern-1em
\endalign
$$
\else
$$
\sign(\)\jb\))\>A_N\>(\voxj)\,\ox
\sumib\sign(\)\ib\))\;T_{i_N\],\)j_N}\bigl(u+(N-1)\>\gm\bigr)
\>\ldots\,T_{i_2,\)j_2}(u+\gm)\,T_{i_1,\)j_1}(u)\,v\,,
\vv-.5>
$$
\fi
which proves the lemma.
\wwgood-:>
\epf
Like in the ordinary linear algebra there is a connection between the
\vv.06>
complementary exterior powers $T^{\wedge k}(u)$ and $T^{\wedge(N-k)}(u)$
and the top exterior power $\TN(u)$, \cf. Theorem~\[TTT]\).
For any two sequences $\aa=(\)a_1\lc a_k)$ and $\bb=(\)b_1\lc b_k)$ set
$$
\Tw k_{\aa\bb}(u)\,=\,
\sumibk\sign(\)\ib\))\;T_{a_{i_k}\]b_k}\bigl(u+(k-1)\>\gm\bigr)
\>\ldots\,T_{a_{i_2}\]b_2}(u+\gm)\,T_{a_{i_1}\]b_1}(u)\,.
\vv-.2>
$$
\Lm{Tab}
Let $\ib\),\jb\in\Sb_k$, and let $\)\aa'\]=(a_{i_1}\)\lc a_{i_k})$,
$\)\bb'\]=(a_{j_1}\)\lc a_{j_k})$. Then
$$
\Tws k_{\aa'\bb'}(u)\,=\,\sign(\)\ib\))\)\sign(\)\jb\))\;\Tw k_{\aa\bb}(u)\,.
$$
\endpro
\nt
The proof is similar to the proof of Lemma~\[TN]\).
\goodbm
\vsk.3>
Denote by $Y_k$ the set of increasing \)\$k\)$-tuples of integers from
\vv.07>
$\lb\)1\lc N\)\rb$. For any ${\aa\in\]Y_k}$ define ${\aab\in\]Y_{N-k}}$
to be the complement of $\aa$, that is,
$\lb\)a_1\lc a_k\),\ab_1\lc\ab_{N-k}\)\rb\)=\)\lb\)1\lc N\)\rb$. Denote
by $\aa\)\aab$ the \perm/ $(\)a_1\lc a_k\),\ab_1\lc\ab_{N-k}\))\in\Sb_N$.
\Th{TTT}
\vvnn0:-.2>
$$
\sign(\aa\)\aab\))\>\sumcc\sign(\cc\)\ccb\))\;
\Tw{(N-k)}_{\ccb\)\aab}\bigl(u-(k-1)\>\gm\bigr)\,\Tw k_{\cc\)\bb}(u)\,=\,
\dl_{\aa\bb}\,\TN(u)\,.
\Tag{TT=T}
$$
\endpro
\nt
The proof is similar to the proof of the analogous formlula in the ordinary
linear algebra.
\vsk.3>
It is clear from relation \(Tph) and formula \(TN) that the \dif/ operator
$\TN(u)$ commutes with multiplication by scalar \fn/s. So, there exists a \fn/
$\LN(u\),\la)\alb\in\Fend V$ \st/
\vvnm-.4>
$$
\bigl(\)\TN(u)\,v\bigr)(\la)\,=\,\LN(u\),\la)\,v(\la)
$$
for any $v\in\FunV$. Denote by $\RN(u\),\la)$ the \fn/ $\LN(u\),\la)$
for the vector \rep/ of $\Eqg$ with the \epoint/ $x=0$.
\Lm{RN}
\vvnn-.7:-.2>
$$
\RN(u\),\la)\,=\;{\tht\bigl(u+(N-1)\>\gm\bigr)\over\tht(u)}\;
\suan\,\)\prod_{\tsize{b=1\atop b\ne a}}^N\>
{\tht(\la_{ab}\]-\gm)\over\tht(\la_{ab})}\ E_{aa}\,.
\vv-.8>
$$
\endpro
\Pf.
Recall that in the vector \rep/
\vvn-.2>
$$
\bigl(\)T_{aa}(u)\,v\bigr)(\la)\,=\,E_{aa}\)v(\la-\gm\>\eps_a)\>+
\tsum_{\tsize{a,b=1\atop a\ne b}}^N\!
\al(u\),\la_{ab})\,E_{bb}\>v(\la-\gm\>\eps_b)
\vv-.4>
$$
and \,$\bigl(\)T_{ab}(u)\,v\bigr)(\la)\)=\)
\bt\)(u\),\la_{ab})\,E_{ba}\)v(\la-\gm\>\eps_b)$, \,where $\al(u\),\xi)$ and
\vv.1>
$\bt\)(u\),\xi)$ are given by \(albt). By Lemma \(TNl) this implies that
$\RN(u\),\la)$ is a linear combination of the matrices $E_{11}\lc E_{NN}$ with
some \fn/al coefficients. Moreover, taking formula \(TN) for the \perm/ $\jb$
\st/ $j_1=\)a$ we observe that in the sum only the term with $\ib=\jb$
contributes to the coefficient of $E_{aa}$, which, therefore, can be easily
found.
\epf
\Pf of Proposition~\[qdet]\).
Following the definition of $\TN(u)$ we find from relations \(Tph) and \(RTT)
that
\ifMag
\vvn-.1>
$$
\align
(\RN)\vpb{(1)}(u-v,\la-\gm\)h\"2)\, & (\TN)\vpb{(2)}(u)\,T\"{12}(v)\,={}
\\
\nn6>
{}=\,{} & T\"{12}(v)\,(\TN)\vpb{(2)}(u)\,(\RN)\vpb{(1)}(u-v,\la+\gm\)h\"1)\,.
\endalign
$$
\else
\vvn.1>
$$
(\RN)\vpb{(1)}(u-v,\la-\gm\)h\"2)\,(\TN)\vpb{(2)}(u)\,T\"{12}(v)\,=\,
T\"{12}(v)\,(\TN)\vpb{(2)}(u)\,(\RN)\vpb{(1)}(u-v,\la+\gm\)h\"1)\,.
\vv.4>
$$
\fi
By Lemma~\[RN] this is \eqv/ to
\vvn.2>
$$
\prod_{\tsize{c=1\atop c\ne a}}^N\>
{\tht(\la_{ac}\]-\gm\)h_{ac}\]-\gm)\over\tht(\la_{ac}\]-\gm\)h_{ac})}
\ \TN(u)\,T_{ab}(v)\,=\,
\prod_{\tsize{c=1\atop c\ne b}}^N\>{\tht(\la_{bc}\]-\gm)\over\tht(\la_{bc})}
\ T_{ab}(v)\,\TN(u)\,.
\vv-.4>
\mmgood
$$
for any ${\)a\),\bno}$. Since
\,$\dsize{\Det\)T(u)\>=\,{\Tht(\la)\over\Tht(\la-\gm\)h)}\ \TN(u)}$
\vvnn0:-.4>\mline
\,where \,$\dsize{\Tht(\la)\>=\!\]\prod_{\abn}\!\tht(\la_{ab})}$,
\,the proposition is proved.
\wwgood-.6:.6>
\epf

\Sect[mform]{Multiplicative forms}
In this section we essentially follow \Cite{\)EV1, Section 1.4\)}\). Notice
that all over the paper the words \em{cocycle} and \em{coboundary} mean
\em{\$1$-cocycle} and \em{\$1$-coboundary}, \resp/.
\vsk.2>
Let $I_k$ be the set of \)\$k\)$-tuples of \pd/ integers from
$\lb\)1\lc N\)\rb$. A \em{\mkform/} $\fb$ is a map ${I_k\)\to\)\FunC}$,
\,${\fb:\)\aa\,\map\)f_{\aa}}$, \>\st/ for any $\aa\in I_k$ and any
$i=1\lc k$ we have
$$
f_{a_1\lc\>a_k}(\la)\,
f_{a_1\lc\>a_{i-1},\>a_{i+1},\>a_i,\>a_{i+2}\lc\>a_k}(\la)\,=\,1
$$
Let $\Om^{\>k}\}$ be the set of all \mkform/s. If $\fb$ and $\gb$ are
\mkform/s, then \,$\fb\)\gb\>:\)\aa\,\map\)f_{\aa}\>g_{\aa}$ \>and
\,$\fb\)/\)\gb\>:\)\aa\,\map\)f_{\aa}\)/\)\gb_{\aa}$ \>are \mkform/s,
which defines an abelian group structure on $\Om^{\>k}\}$. The neutral
element is the form $\onb$: $\one_{\)\aa}(\la)=1$ for any $\aa\in I_k$.
\vsk.2>
For any nonzero \fn/ $f(\la)$ and any $\ano$ \>set
\>$(\dl_a\)f)(\la)\>=f(\la-\gm\>\eps_a)\)/f(\la)$, \>and
\vv.1>
for any $\fb\in\Om^{\>k}$ define a \mform/ $d\)\fb\in\Om^{\>k+1}\}$ by
\vvnm-.1>
$$
(d\)\fb\))_{\)a_1\lc\>a_{k+1}}(\la)\,=\>\prod_{i=1}^{k+1}\>\bigl((\dl_{a_i}
f_{a_1\lc\>a_{i-1},\>a_{i+1}\lc\>a_{k+1}})(\la)\bigr)\vpb{(-1)^{i-1}}.\kern-1em
\vvm-.1>
$$
For any \mform/ $\fb$ we have $d^{\>2}\}\fb=\onb$. The \mform/ $\fb$ is called
a \em{\mccl/} if $d\)\fb=\onb$, and a \em{\mcbd/} if $\fb=d\>\gb$ for
a suitable \mform/ $\gb$.
\Par
For any \mef/ $f(\xi)$ in one \var/ the \mult/ \$1$-form $\Fb=(F_1\lc F_N)$:
\vvnm.3>
$$
F_a(\la)\,=\!\prod_{1\le c<a}\}f(\la_{ca})
\prod_{a<c\le N}\}\bigl(f(\la_{ac}-\gm)\bigr)\vpb{-1},
\vvmm.3:-.1>
$$
is a \mcclo/. If $f(\xi)\)=\)g(\xi+\gm)\)/\)g(\xi)$,
\vv.1>
then $\Fb$ is a \mult/ \$1$-coboundary: $\Fb=d\>G$, where
$G(\la)\>=\!\prod_{1\le b<c\le N}g(\la_{bc})$.
\vsk.5>
If $(F_1\lc F_N)$ is a \mcclo/, then so is $(\)G_1\lc G_N)$:
$G_a(\la)\)=F_a(\)-\)\la+\gm\>\eps_a)$.

\Sect[norm]{Proof of Theorem \{\[PBW]}
In this section we introduce another ordering on generators of $\els$, called
ordering by rows, and prove the analogue of Theorem~\[PBW] for the monomials
\obr/, \cf. Theorem~\[pbw]\). Since the number or \obrm/s of degree $k$ equals
the number of \norm/s of the same degree, and each monomial can be transformed
to a linear combination of \norm/s, Theorem~\[pbw] implies Theorem~\[PBW]\).
We introduce \obrm/s for technical reason because this allows us to reduce
the number of cases to be examined at some stage of the proof.
\Par
Consider the \em{ordering by rows} of generators: ${t_{ab}\]\prec\)t_{cd}}$
if $a<c$, or $a=c$ and $b<d$. Say that the monomial $\tabk$ is \em{\obr/}
if $t_{a_ib_i}\!\prec\)t_{a_jb_j}$ for any $i<j$, or $k=0$.
\Th{pbw}
For any $k\in\Zp$ the \obrm/s of degree $k$ form a basis of \,$\eg_k\]$ over
$\Funt$.
\endpro
\vsk-.5>\vsk0>
\Pf.
To save space from now on we write \emph{ordered} \)instead of \emph{\obr/}.
We call an equality
\vvn-.2>
$$
\text{\it dis\orm/\/}\,=\,\text{\it linear combination of \orm/s\/}
\vvmm0:-.15>
$$
an \em{\orule/} for the monomial in \lhs/.
The commutation relations \(tbc)\,--\,\(tac) have the important property:
\vv.1>
\itemhalf{A.}
Any relation is a linear combination of \orule/s, and the complete list of
linear independent \orule/s contains precisely one rule for each disordered
product of generators.
\vsk.3>
For ${k=0}$ and ${k=1}$ the claim of Theorem~\[pbw] is immediate. Let ${k>1}$.
First we prove that \orm/s of degree $k$ span \)$\eg_k$ over $\Funt$. Indeed,
one can transform any monomial $\tabk$ to a linear combination of \orm/s by
the following procedure. Pick up any disordered product of adjacent factors
and replace it by a sum of ordered products using the \orule/, then do the same
for each of the obtained monomials. To see that the procedure terminates and,
hence, produces a linear combination of \orm/s, introduce auxilary gradings
on monomials by the rule
\vvn-.1>
$$
r\)(\tabk)\>=\>\tsum_{i=1}^k\>i\>a_i\,,\kern4em
r'(\tabk)\>=\>\tsum_{i=1}^k\>i\>b_i\,,
\vv.2>
$$
and observe that at each nontrivial step of the procedure we replace a monomial
by a sum of monomials of either less degree $r$, or the same degree $r$ and
less degree $r'\]$. We call the described procedure a \em{\regtr/} of
the monomial $\tabk$ to a linear combination of \orm/s.
\vsk.2>
If an \orule/ is applied to a product $t_{a_ib_i}t_{a_{i+1}b_{i+1}}$ in
$\tabk$, we say that the \orule/ is applied \em{at \$\}i$-th place}.
\vsk.3>
By the standard reasoning the property A of the commutation relations
\(tbc)\,--\,\(tac) implies that Theorem~\[pbw] follows from Proposition~\[uni]\).
\epf
\Prop{uni}
Any \regtr/ of the monomial $\tabk$ produces the same linear combination of
\orm/s.
\endpro
\Pf.
For $k=2$ the claim follows from the property A. For $k=3$ the claim can be
verified in a straightforward way. We discuss more details of $k=3$ case at
the end of the proof.
\vsk.3>
Let $k>3$. For the proof we use induction \wrt/ the \lex/ ordering on monomials
defined by a pair of degrees $(r\),r')$. The claim of the proposition is
obvious for \orm/s, which provides the base of induction.
\vsk.2>
Let $\I$ and $\II$ be two \regtr/s of the monomial $\tabk$ to linear
\mmgood
combinations of \orm/s. If for $\I$ and $\II$ the first steps coincide, then
they produce the same results by the induction assumption. Otherwise, let us
construct two additional \regtr/s $\III$ and $\IV$ \st/ the first steps of $\I$
and $\III$ coincide, the same holds for $\II$ and $\IV$, and $\III$ and $\IV$
produce the same results. By the previous remark this proves the proposition.
\vsk.2>
Assume that for the \trans/ $\I$ an \orule/ at the first step is applied
at \$\]i$-th place, and for the \trans/ $\II$ at \$\]j$-th place.
If $|\>i-j\)|>1$ then we define the \trans/ $\III$ as follows: first apply an
\orule/ at \$\]i$-th place, next apply an \orule/ at \$\]j$-th place
for all monomials obtained at the first step, then continue in any possible
way. The \trans/ $\IV$ is defined similarly with $i$ and $j$ interchanged.
By the induction assumption the \trans/s $\III$ and $\IV$ produce the same
results because after the first two steps of both $\III$ and $\IV$ one has
identical linear combinations of monomials, each of them being \lex/ly smaller
than the initial monomial $\tabk$.
\vsk.2>
The cases $j=i\pm 1$ are similar. Assume for example that $j=i+1$.
Define the \trans/ $\III$ as follows: apply a \regtr/ of the product $\tabii$
to a linear combination of ordered triple products, making the first step at
\$\]i$-th place, then continue in any possible way. Define the \trans/
$\IV$ similarly, but making the first step at \$\](i+1)$-th place. Then
the claim of the proposition for $k=3$ shows that after the first stages of
both $\III$ and $\IV$ one has identical linear combinations of monomials,
each of them being \lex/ly smaller than the initial monomial $\tabk$.
Therefore, by the induction assumption the \trans/s $\III$ and $\IV$
produce the same results.
\vsk.2>
It remains to complete the proof for $k=3$. For any monomial $\tabt$
there are at most two \regtr/s, and the \regtr/ is unique unless
$t_{a_3b_3}\!\prec\)t_{a_2b_2}\!\prec\)t_{a_1b_1}$. The rest of the proof
is given by the straightforward calculations. The simplest cases occur
if $a_1\]=a_2\]=a_3$ or $b_1\]=b_2\]=b_3$. We present below the most
bulky example when $a_1\]>a_2\]>a_3$ and $b_1\),\)b_2\),\)b_3$ are \pd/.
\epf
\Ex
Regular \trans/s of the monomial $t_{34}\>t_{25}\>t_{16}$.
Let \fn/s $\al(u\),\xi)$ and $\bt\)(u\),\xi)$ be given by \(albt).
Making the first step at the first place we have
\ifMag
$$
\align
t_{34}\>t_{25}\> & t_{16}\,=\,
\al(\la\+1_{32}\},\la\+2_{54})\,t_{25}\>t_{34}\>t_{16}\>-\)
\bt\)(\la\+1_{32}\},\la\+2_{45})\,t_{24}\>t_{35}\>t_{16}\,={}
\\
\ald
\nn9>
{}=\,{} & \al(\la\+1_{32}\},\la\+2_{54})\>
\bigl(\)\al(\la\+1_{31}\},\la\+2_{64})\,t_{25}\>t_{16}\>t_{34}\>-\)
\bt\)(\la\+1_{31}\},\la\+2_{46})\,t_{25}\>t_{14}\>t_{36}\bigr)\>+{}
\\
\nn6>
{}-\,{} & \bt\)(\la\+1_{32}\},\la\+2_{45})\>
\bigl(\)\al(\la\+1_{31}\},\la\+2_{65})\,t_{24}\>t_{16}\>t_{35}\>-\)
\bt\)(\la\+1_{31}\},\la\+2_{56})\,t_{24}\>t_{15}\>t_{36}\bigr)\,={}
\\
\ald
\nn9>
{}=\,{} & \al(\la\+1_{32}\},\la\+2_{54})\,\al(\la\+1_{31}\},\la\+2_{64})\>
\bigl(\)\al(\la\+1_{21}\},\la\+2_{65})\,t_{16}\>t_{25}\>t_{34}\>-\)
\bt\)(\la\+1_{21}\},\la\+2_{56})\,t_{15}\>t_{26}\>t_{34}\bigr)\>+{}
\\
\nn6>
{}-\,{} & \al(\la\+1_{32}\},\la\+2_{54})\,\bt\)(\la\+1_{31}\},\la\+2_{46})\>
\bigl(\)\al(\la\+1_{21}\},\la\+2_{45})\,t_{14}\>t_{25}\>t_{36}\>-\)
\bt\)(\la\+1_{21}\},\la\+2_{54})\,t_{15}\>t_{24}\>t_{36}\bigr)\>+{}
\\
\nn6>
{}-\,{} & \bt\)(\la\+1_{32}\},\la\+2_{45})\,\al(\la\+1_{31}\},\la\+2_{65})\>
\bigl(\)\al(\la\+1_{21}\},\la\+2_{64})\,t_{16}\>t_{24}\>t_{35}\>-\)
\bt\)(\la\+1_{21}\},\la\+2_{46})\,t_{14}\>t_{26}\>t_{35}\bigr)\>+{}
\\
\nn6>
{}+\,{} & \bt\)(\la\+1_{32}\},\la\+2_{45})\,\bt\)(\la\+1_{31}\},\la\+2_{56})\>
\bigl(\)\al(\la\+1_{21}\},\la\+2_{54})\,t_{15}\>t_{24}\>t_{36}\>-\)
\bt\)(\la\+1_{21}\},\la\+2_{45})\,t_{14}\>t_{25}\>t_{36}\bigr)\,,
\\
\cnn.1>
\endalign
$$
\else
$$
\align
t_{34}\>t_{25}\>t_{16}\, &{}=\,
\al(\la\+1_{32}\},\la\+2_{54})\,t_{25}\>t_{34}\>t_{16}\>-\)
\bt\)(\la\+1_{32}\},\la\+2_{45})\,t_{24}\>t_{35}\>t_{16}\,={}
\\
\ald
\nn9>
& {}=\,\al(\la\+1_{32}\},\la\+2_{54})\>
\bigl(\)\al(\la\+1_{31}\},\la\+2_{64})\,t_{25}\>t_{16}\>t_{34}\>-\)
\bt\)(\la\+1_{31}\},\la\+2_{46})\,t_{25}\>t_{14}\>t_{36}\bigr)\>+{}
\\
\nn6>
&\){}-\,\bt\)(\la\+1_{32}\},\la\+2_{45})\>
\bigl(\)\al(\la\+1_{31}\},\la\+2_{65})\,t_{24}\>t_{16}\>t_{35}\>-\)
\bt\)(\la\+1_{31}\},\la\+2_{56})\,t_{24}\>t_{15}\>t_{36}\bigr)\,={}
\\
\ald
\nn9>
& {}=\,\al(\la\+1_{32}\},\la\+2_{54})\,\al(\la\+1_{31}\},\la\+2_{64})\>
\bigl(\)\al(\la\+1_{21}\},\la\+2_{65})\,t_{16}\>t_{25}\>t_{34}\>-\)
\bt\)(\la\+1_{21}\},\la\+2_{56})\,t_{15}\>t_{26}\>t_{34}\bigr)\>+{}
\\
\nn6>
&\){}-\,\al(\la\+1_{32}\},\la\+2_{54})\,\bt\)(\la\+1_{31}\},\la\+2_{46})\>
\bigl(\)\al(\la\+1_{21}\},\la\+2_{45})\,t_{14}\>t_{25}\>t_{36}\>-\)
\bt\)(\la\+1_{21}\},\la\+2_{54})\,t_{15}\>t_{24}\>t_{36}\bigr)\>+{}
\\
\nn6>
&\){}-\,\bt\)(\la\+1_{32}\},\la\+2_{45})\,\al(\la\+1_{31}\},\la\+2_{65})\>
\bigl(\)\al(\la\+1_{21}\},\la\+2_{64})\,t_{16}\>t_{24}\>t_{35}\>-\)
\bt\)(\la\+1_{21}\},\la\+2_{46})\,t_{14}\>t_{26}\>t_{35}\bigr)\>+{}
\\
\nn6>
&\){}+\,\bt\)(\la\+1_{32}\},\la\+2_{45})\,\bt\)(\la\+1_{31}\},\la\+2_{56})\>
\bigl(\)\al(\la\+1_{21}\},\la\+2_{54})\,t_{15}\>t_{24}\>t_{36}\>-\)
\bt\)(\la\+1_{21}\},\la\+2_{45})\,t_{14}\>t_{25}\>t_{36}\bigr)\,,
\\
\cnn.1>
\endalign
$$
\fi
while making the first step at the second place we have
\vvn.1>
\ifMag
$$
\align
t_{34}\>t_{25}\> & t_{16}\,=\,
\al(\la\+1_{21}\},\la\+2_{65})\,t_{34}\>t_{16}\>t_{25}\>-\)
\bt\)(\la\+1_{21}\},\la\+2_{56})\,t_{34}\>t_{15}\>t_{26}\,={}
\\
\ald
\nn9>
\ald
{}=\,{} & \al(\la\+1_{21}\},\la\+2_{65})\>
\bigl(\)\al(\la\+1_{31}\},\la\+2_{64})\,t_{16}\>t_{34}\>t_{25}\>-\)
\bt\)(\la\+1_{31}\},\la\+2_{46})\,t_{14}\>t_{36}\>t_{25}\bigr)\>+{}
\\
\nn6>
{}-\,{} & \bt\)(\la\+1_{21}\},\la\+2_{56})\>
\bigl(\)\al(\la\+1_{31}\},\la\+2_{54})\,t_{15}\>t_{34}\>t_{26}\>-\)
\bt\)(\la\+1_{31}\},\la\+2_{45})\,t_{14}\>t_{35}\>t_{26}\bigr)\,={}
\\
\ald
\nn9>
\ald
{}=\,{} & \al(\la\+1_{21}\},\la\+2_{65})\,\al(\la\+1_{31}\},\la\+2_{64})\>
\bigl(\)\al(\la\+1_{32}\},\la\+2_{54})\,t_{16}\>t_{25}\>t_{34}\>-\)
\bt\)(\la\+1_{32}\},\la\+2_{45})\,t_{16}\>t_{24}\>t_{35}\bigr)\>+{}
\\
\nn6>
{}-\,{} & \al(\la\+1_{21}\},\la\+2_{65})\,\bt\)(\la\+1_{31}\},\la\+2_{46})\>
\bigl(\)\al(\la\+1_{32}\},\la\+2_{56})\,t_{14}\>t_{25}\>t_{36}\>-\)
\bt\)(\la\+1_{32}\},\la\+2_{65})\,t_{14}\>t_{26}\>t_{35}\bigr)\>+{}
\\
\nn6>
{}-\,{} & \bt\)(\la\+1_{21}\},\la\+2_{56})\,\al(\la\+1_{31}\},\la\+2_{54})\>
\bigl(\)\al(\la\+1_{32}\},\la\+2_{64})\,t_{15}\>t_{26}\>t_{34}\>-\)
\bt\)(\la\+1_{32}\},\la\+2_{46})\,t_{15}\>t_{24}\>t_{36}\bigr)\>+{}
\\
\nn6>
{}+\,{} & \bt\)(\la\+1_{21}\},\la\+2_{56})\,\bt\)(\la\+1_{31}\},\la\+2_{45})\>
\bigl(\)\al(\la\+1_{32}\},\la\+2_{65})\,t_{14}\>t_{26}\>t_{35}\>-\)
\bt\)(\la\+1_{32}\},\la\+2_{56})\,t_{14}\>t_{25}\>t_{36}\bigr)\,.
\\
\cnn.1>
\endalign
$$
\else
$$
\align
t_{34}\>t_{25}\>t_{16}\, &{}=\,
\al(\la\+1_{21}\},\la\+2_{65})\,t_{34}\>t_{16}\>t_{25}\>-\)
\bt\)(\la\+1_{21}\},\la\+2_{56})\,t_{34}\>t_{15}\>t_{26}\,={}
\\
\ald
\nn9>
\ald
& {}=\,\al(\la\+1_{21}\},\la\+2_{65})\>
\bigl(\)\al(\la\+1_{31}\},\la\+2_{64})\,t_{16}\>t_{34}\>t_{25}\>-\)
\bt\)(\la\+1_{31}\},\la\+2_{46})\,t_{14}\>t_{36}\>t_{25}\bigr)\>+{}
\\
\nn6>
&\){}-\,\bt\)(\la\+1_{21}\},\la\+2_{56})\>
\bigl(\)\al(\la\+1_{31}\},\la\+2_{54})\,t_{15}\>t_{34}\>t_{26}\>-\)
\bt\)(\la\+1_{31}\},\la\+2_{45})\,t_{14}\>t_{35}\>t_{26}\bigr)\,={}
\\
\ald
\nn9>
\ald
& {}=\,\al(\la\+1_{21}\},\la\+2_{65})\,\al(\la\+1_{31}\},\la\+2_{64})\>
\bigl(\)\al(\la\+1_{32}\},\la\+2_{54})\,t_{16}\>t_{25}\>t_{34}\>-\)
\bt\)(\la\+1_{32}\},\la\+2_{45})\,t_{16}\>t_{24}\>t_{35}\bigr)\>+{}
\\
\nn6>
&\){}-\,\al(\la\+1_{21}\},\la\+2_{65})\,\bt\)(\la\+1_{31}\},\la\+2_{46})\>
\bigl(\)\al(\la\+1_{32}\},\la\+2_{56})\,t_{14}\>t_{25}\>t_{36}\>-\)
\bt\)(\la\+1_{32}\},\la\+2_{65})\,t_{14}\>t_{26}\>t_{35}\bigr)\>+{}
\\
\nn6>
&\){}-\,\bt\)(\la\+1_{21}\},\la\+2_{56})\,\al(\la\+1_{31}\},\la\+2_{54})\>
\bigl(\)\al(\la\+1_{32}\},\la\+2_{64})\,t_{15}\>t_{26}\>t_{34}\>-\)
\bt\)(\la\+1_{32}\},\la\+2_{46})\,t_{15}\>t_{24}\>t_{36}\bigr)\>+{}
\\
\nn6>
&\){}+\,\bt\)(\la\+1_{21}\},\la\+2_{56})\,\bt\)(\la\+1_{31}\},\la\+2_{45})\>
\bigl(\)\al(\la\+1_{32}\},\la\+2_{65})\,t_{14}\>t_{26}\>t_{35}\>-\)
\bt\)(\la\+1_{32}\},\la\+2_{56})\,t_{14}\>t_{25}\>t_{36}\bigr)\,.
\\
\cnn.1>
\endalign
$$
\fi
The coefficients of the monomials $t_{16}\>t_{25}\>t_{34}$,
$t_{15}\>t_{26}\>t_{34}$ and $t_{16}\>t_{24}\>t_{35}$ coincide identically.
To check that the coefficients for other monomials are the same one should take
into account the explicit form of $\al(u\),\xi)$ and $\bt\)(u\),\xi)$ and use
summation formulae for the theta \fn/.
\enddemo

\Sect[sl2]{Elliptic \qg/ $\{\eslt$}
This section is an illustration of general constructions in the simplest case
${N=2}$. In this case $\om_1=\eps_1=-\eps_2=\rho=\al_1/2$ and
$\hga\}=\C\>\om_1$. For any $\la\in\hga\}$ we have $\la_1=-\)\la_2=\la_{12}/2$
and $\la=\la_{12}\>\om_1$. We identify \fn/s on $\hga\}$ with \fn/s of one
\var/ $\la_{12}$.
\vsk.2>
The \oalg/ $\elst$ is generated over $\C$ by elements
$t_{11}\),\)t_{12}\),\)t_{21}\),\)t_{22}$ and \fn/s $f\}\in\Funt$. According
to Theorem~\[PBW] monomials
$t_{21}^{\>k_{21}}\)t_{11}^{\>k_{11}}\)t_{22}^{\>k_{22}}\)t_{12}^{\>k_{12}}\}$,
\>$k_{11}\),\) k_{12}\),\) k_{21}\),\){k_{22}\in\Zp}$, form a basis of $\elst$
as a \$\Funt$-module.
\vsk.2>
Let $\Qh\)=(\)Q_1, Q_2)$ be a \mccl/, which in this case means that
\vvnm.1>
$$
Q_1(\la_{12})\,Q_2(\la_{12}\]-\gm)\,=\,Q_1(\la_{12}\]+\gm)\,Q_2(\la_{12})\,.
\vvm.1>
$$
\vv-.1>
A \Vmod/ $\MmQ$ of \hw/ ${\mu\in\hga\}}$ and \dhw/ $\Qh$ over $\eslt$ is
constructed as follows. As an \hmod/ ${\MmQ\,=\)\Plus_{k=0}^\8\>\C\>\vmk k}$,
\vv.1>
\>the vector $\vmk k$ being of \wt/ $\mu-k\)\al_1$. Notice that $\vmQ=\vmk 0$.
The generators $\tth_{11}\),\)\tth_{12}\),\)\tth_{21}\),\)\tth_{22}$ act on
$\Fun(\MmQ)$ by the rule
\vvnn-.1:-.5>
$$
\alds
\gather
\tth_{21}\,\vmk k\,=\,\vmk{k+1}\,,
\\
\nnm8:10>
\tth_{11}\,\vmk k\,=\,Q_1(\la_{12}\]+k\)\gm)
\ {\tht(\la_{12}\]+k\)\gm)\over\tht(\la_{12})}\ \vmk k\,,
\\
\nnm6:8>
\tth_{22}\,\vmk k\,=\,Q_2(\la_{12}\]+k\)\gm)
\ {\tht\bigl(\la_{12}\]-(\)\mu_{12}\]-k)\>\gm\bigr)\over
\tht\bigl(\la_{12}\]-(\)\mu_{12}\]-2\)k)\>\gm\bigr)}\ \vmk k\,,
\\
\nnm5:8>
\ifMag
{\align
\tth_{12}\, & \vmk k\,=\,
\\
{}=\,{} & Q_1(\la_{12}\]+k\)\gm)\,Q_2\bigl(\la_{12}\]+(k-1)\>\gm\bigr)
\ {\tht\bigl((\)\mu_{12}\]-k+1)\>\gm\bigr)\,\tht(k\)\gm)\over\tht(\la_{12})\,
\tht\bigl(\la_{12}\]-(\)\mu_{12}\]-2\)k+2)\>\gm\bigr)}\ \vmk{k-1}\,,
\endalign}
\else
\tth_{12}\,\vmk k\,=\,
Q_1(\la_{12}\]+k\)\gm)\,Q_2\bigl(\la_{12}\]+(k-1)\>\gm\bigr)
\ {\tht\bigl((\)\mu_{12}\]-k+1)\>\gm\bigr)\,\tht(k\)\gm)\over\tht(\la_{12})\,
\tht\bigl(\la_{12}\]-(\)\mu_{12}\]-2\)k+2)\>\gm\bigr)}\ \vmk{k-1}\,,
\fi
\endgather
$$
see \(aa1)\), \(a1a) and Lemma~\[abba]\). Taking into account relations \(tth)
and \(Tcu)\), and Corollary \[ehom]\), one reproduces from these formulae
the construction, given in \Cite{FV1}, of the \eVmod/ over $\Eqgt$ tensored
with a \onedim/ \rep/ of $\Eqgt$.
\vsk.3>
In general, to compute the \dSform/ on $\elst$ one has to find
$$
S\bigl(\)f(\la\+1_{12},\la\+2_{12}))\,
\tth_{21}^{\>k_{21}}\)\tth_{11}^{\>k_{11}}\)
\tth_{22}^{\>k_{22}}\)\tth_{12}^{\>k_{12}}\],
\,g(\la\+1_{12}\},\la\+2_{12})\,\tth_{21}^{\>m_{21}}\)\tth_{11}^{\>m_{11}}\)
\tth_{22}^{\>m_{22}}\)\tth_{12}^{\>m_{12}}\bigr)
\Tag{Sfg}
$$
for any $f\),\)g\in\Funt$ and any monomials $\tth_{21}^{\>k_{21}}\)
\tth_{11}^{\>k_{11}}\)\tth_{22}^{\>k_{22}}\)\tth_{12}^{\>k_{12}}\}$,
\>$\tth_{21}^{\>m_{21}}\)\tth_{11}^{\>m_{11}}\)\tth_{22}^{\>m_{22}}\)
\tth_{12}^{\>m_{12}}\}$. The answer is zero unless $k_{12}=\)m_{12}=\)0$ and
$k_{21}=\)m_{21}$. Moreover, contributions of the \fn/s $f\),\)g$ and factors
$t_{aa}$ are easy to calculate, and we find that in the nonzero case expression
\(Sfg) equals
\ifMag
$$
\align
f\bigl(\>-\>\la\+2_{12}+(k_{11}\]-k_{21}\]-k_{22})\>\gm\>,
\>-\>\la\+1_{12}+(k_{11}\]+k_{21}\]-k_{22})\>\gm\bigr)\,\x{} &
\\
\nn6>
{}\x\,g\bigl(\)\la\+1_{12}-(k_{11}\]+k_{21}\]-k_{22})\>\gm\),
\>\la\+2_{12}-(k_{11}\]-k_{21}\]-k_{22})\>\gm\bigr)\,\x{}&
\\
\nn6>
{}\x\,q_1^{\>k_{11}}\>q_2^{\>k_{22}}\,
S(\)\tth_{21}^{\>k_{21}}\>,\)\tth_{21}^{\>k_{21}})\,
q_1^{\>m_{11}}\>q_2^{\>m_{22}}\,\) &
\\
\cnn-.5>
\endalign
$$
\else
$$
\align
& f\bigl(\>-\>\la\+2_{12}+(k_{11}\]-k_{21}\]-k_{22})\>\gm\>,
\>-\>\la\+1_{12}+(k_{11}\]+k_{21}\]-k_{22})\>\gm\bigr)\,\x{}
\\
\nn6>
& \!{}\x\,g\bigl(\)\la\+1_{12}-(k_{11}\]+k_{21}\]-k_{22})\>\gm\),
\>\la\+2_{12}-(k_{11}\]-k_{21}\]-k_{22})\>\gm\bigr)\,
q_1^{\>k_{11}}\>q_2^{\>k_{22}}\,
S(\)\tth_{21}^{\>k_{21}}\>,\)\tth_{21}^{\>k_{21}})\,
q_1^{\>m_{11}}\>q_2^{\>m_{22}}
\endalign
$$
\fi
where
\vvnm-.5>
$$
S(\)\tth_{21}^{\>k}\>,\)\tth_{21}^{\>k})\,=\,(-1)^k\,\prod_{j=0}^{k-1}\,
{\tht(\la\+1_{12}\!-\la\+2_{12}\!-j\)\gm))\,\tht\bigl((j+1)\>\gm\bigr)
\over\tht(\la\+1_{12}\!-j\)\gm)\,\tht(\la\+2_{12}\!+j\)\gm)}
\ q_1^{\>k}\>q_2^{\>k}\,,\qquad k\in\Zp\,.\kern-2em
\vv.4>
$$
For the corresponding part of the \dSpair/ for $\MmQ$ we have
\ifMag
\vvnm-.2>
$$
\align
\SmQ\)\bigl(\)\tth_{21}^{\>k}\>,\)\vmk k\)\bigr)\,=\,(-1)^k\,\prod_{j=0}^{k-1}
\, & \Bigl(Q_1\bigl(\la_{12}\]+(j+1)\>\gm\bigr)\,Q_2(\la_{12}\]+j\)\gm)\,\x{}
\\
\nn1>
{}\x\,\){}&\){\tht\bigl((\)\mu_{12}\]-j)\>\gm\bigr)\,\tht\bigl((j+1)\>\gm\bigr)
\over\tht(\la_{12}\]-j\)\gm)\,
\tht\bigl(\la_{12}\]-(\)\mu_{12}\]-j)\>\gm\bigr)}\,\Bigr)\,.
\\
\cnn.2>
\endalign
$$
\else
\vvn-.1>
$$
\SmQ\)\bigl(\)\tth_{21}^{\>k}\>,\)\vmk k\)\bigr)\,=\,(-1)^k\,\prod_{j=0}^{k-1}
\,\Bigl(Q_1\bigl(\la_{12}\]+(j+1)\>\gm\bigr)\,Q_2(\la_{12}\]+j\)\gm)
\ {\tht\bigl((\)\mu_{12}\]-j)\>\gm\bigr)\,\tht\bigl((j+1)\>\gm\bigr)
\over\tht(\la_{12}\]-j\)\gm)\,
\tht\bigl(\la_{12}\]-(\)\mu_{12}\]-j)\>\gm\bigr)}\,\Bigr)\,.
$$
\fi
The \cform/ \,$\CmQ:\Fun(\MmQ)\oxC\Fun(\MmQt)\>\to\)\FunC$ \,is given by
\vvnn0:-.2>
$$
\align
\CmQ\)\bigl(\vmk k\>,\)v_{\mu,\Qti\)}[\>l\>]\)\bigr)\,=\,
\dl_{kl}\>(-1)^k\,\prod_{j=0}^{k-1}\,
& \Bigl(Q_1\bigl(\la_{12}\]+(j+1)\>\gm\bigr)\,Q_2(\la_{12}\]+j\)\gm)\,\x{}
\\
\nn1>
{}\x\,\){}&\){\tht\bigl((\)\mu_{12}\]-j)\>\gm\bigr)\,\tht\bigl((j+1)\>\gm\bigr)
\over\tht\bigl(\la_{12}\]+(j+1)\>\gm\bigr)\,
\tht\bigl(\la_{12}\]-(\)\mu_{12}\]-j-k)\>\gm\bigr)}\,\Bigr)\,.
\endalign
$$
\par
Let ${\gm\)\nin\)\QtQ}$ and assume that $\Qh$ is \ndeg/. Then the module
$\MmQ$ is \irr/ provided that ${\mu_{12}\nin\Zp}$. If ${\mu_{12}\in\Zp}$,
\vvm.08>
then $\vmk{\mu_{12}\]+1}$ is a \regsv/ generating an \irr/ submodule $\NmQ$
\vvm-.05>
isomorphic to $M_{-\mu-2\rho,\)\Qti}$ where
\vv.1>
$\Qti(\la_{12})=\)\Qh(\la_{12}\]-\mu_{12}\]-2)$. The quotient module
\vvmm.05:-.1>
${\VmQ\)=\MmQ\)/\NmQ}$ is the \irr/ \hw/ \$\eslt$-module
with \hw/ $\mu$ and \dhw/ $\Qh$, and it has dimension $\mu_{12}\]+1$, the same
as the \irr/ \$\sltwo$-module of \hw/ $\mu$.

\Sect[B]{Proof of Theorem \{\[Ect]}
In \Cite{\)TV2\), Section~2.6\)} for any semisimple Lie algebra $\g$ we have
defined \raf/s $B_w(\la)$ of $\la$ labeled by elements of the Weyl group.
The \fn/s takes values in a certain completion of $U(\g)$. Here we need
the \fn/ $B_w(\la-\rho\))$ for the particular case of $\g=\sln$ and $w=w_X\:$,
the longest element of the Weyl group. For brevity we denote this \fn/ by
$B_X\:(\la)$. We list required properties of $B_X\:(\la)$ below using
the notation of the present paper. All of them easily follow from
the properties of $B_w(\la)$ given in \Cite{TV2}\).
\vsk.3>
Consider the following series
\vvnn-.1:-.4>
$$
G_{ab}(\la)\,=\,\sum_{k=0}^\8\,e_{ba}^{\>k}\)e_{ab}^{\>k}\,
\prod_{j=1}^k\,{1\over j\>(\la_{ab}\]-e_{aa}\]+e_{bb}-j\))}\;.
$$
Then $B_X\:(\la)$ equals the ordered product
\vv.04>
$\dsize\prod_{\abn}\!\!G_{ab}(\la)$ \,where the factor $G_{ab}$ is to the
left from the factor $G_{cd}$ if $a<c$, or $a=c$ and $b>d$. For instance,
\vvm.1>
if $N\]=3$, then $B_X\:(\la)\)=\)G_{12}(\la)\>G_{13}(\la)\>G_{23}(\la)$.
\goodbreak
\vsk.2>
The series $B_X\:(\la)$ acts on any \fd/ \smod/ $V\}$, commuting with
the \$\hg$-action. The action of $B_X\:(\la)$ gives an element of $\Rend V$,
which tends to $1$ as $\la$ goes to infinity in a generic direction and,
therefore, is invertible for generic $\la$.
\Prop{BVW}
For any \fd/ \smod/s $V,\)W\}$ we have
\vvnn.3:.5>
$$
B_X\:(\la)\vst{V\ox\)W}\>
(X\]\ox\]X)\>J_{VW}\:\bigl(w_X\:(\la)\bigr)\)(X\]\ox\]X)\1\)=\,
J_{VW}\:(\la)\)\bigl(B_X\:(\la-h\"2)\)\ox\)B_X\:(\la)\bigr)\,.
\vvm-.1>
$$
\endpro
\Cr{RBB}
\ifMag
\vvn-.1>
$$
\Rti_{VW}(\la-\rho\))\,=\,\bigl(B_X\:(\la-h\"2)\)\ox\)B_X\:(\la)\bigr)
\vpb{-1}R_{VW}\:(\la)\)\bigl(B_X\:(\la)\)\ox\)B_X\:(\la-h\"1)\bigr)\,.
$$
\else
\>$\Rti_{VW}(\la-\rho\))\,=\,\bigl(B_X\:(\la-h\"2)\)\ox\)B_X\:(\la)\bigr)
\vpb{-1}R_{VW}\:(\la)\)\bigl(B_X\:(\la)\)\ox\)B_X\:(\la-h\"1)\bigr)$.
\fi
\endpro
\vsk.2>
\nt
The statement follows from Corollary~\[RRt] and cocommutativity of
the coproduct $\Dl$.
\Lm{BU}
Let $U\]$ be the vector \rep/ of $\sln$. Then
\ifMag
$$
B_X(\la)\vst U\>=\,\tsuan E_{aa}\!\tprod_{1\le b<a}\!(1+\la_{ba}\1)\,.
\vv-.2>
$$
\else
\,$B_X(\la)\vst U\)=\>\suan E_{aa}\!\tprod_{1\le b<a}\!(1+\la_{ba}\1)$.
\fi
\endpro
\Pf of Theorem~\[Ect]\).
By the definition of isomorpic functors the assertion of the theorem
\vv.05>
means that for any \fd/ \smod/ $V\}$ there exists an \iso/
$\psi_V\:\]\in\Mor\bigl(\Ect(V)\>,\)\Ec(V)\bigr)$,
and for any morphism $\phi:V\]\to\)W\}$ of \smod/s one has
\,$\Ec(\phi)\o\psi_V\:=\)\psi_W\:\}\o\Ect(\phi)$.
\vsk.4>
It follows from formula \(Gc)\), Proposition~\[L2t]\), Corollary~\[RBB] and
\vvm.06>
Lemma~\[BU] that $B_X\:(\la)\vst V$ is a morphism from $\Ect(V)$ to $\Ec(V)$.
\ifMag\else\vv.1>\fi
It is an \iso/, since $B_X\:(\la)\vst V$ is invertible for generic $\la$.
\vvmm0:-.08>
Moreover, for any morphism $\phi:V\]\to\)W\}$ of \smod/s we have
\,$\phi\,B_X\:(\la)\vst V\]=B_X\:(\la)\vst W\,\phi$.
Since both $\Ec$ and $\Ect$ send $\phi$ to the constant \fn/
$\phi\in\Rat\bigl(\Hom(V,W)\bigr)$, the theorem is proved.
\wwgood.3:0>
\epf

\myRefs
\widest{FV1}
\parskip.1\bls

\ref\Key ABBR
\by D\&Arnaudon, E\&Buffenoir, E\&Ragoucy and Ph\&Roche
\paper Universal \sol/s of quantum \DYB/s
\jour \LMP/ \vol 44 \issue 3 \yr 1998 \pages 201\)--\)214
\endref

\ref\Key C
\by \Cher/
\paper ``Quantum'' deformations of \irr/ \fd/ \rep/s of $\gln$
\jour Soviet Math\. Dokl. \vol 33 \issue 2 \yr 1986 \pages 507\)--\)510
\endref

\ref\Key ESS
\by \Etingof/, T\]\&Schedler and O\&Schiffmann
\paper Explicit quantization of dynamical \$r\]$-matrices for \fd/
semisimple Lie algebras
\jour J\. \AMS/ \vol 13 \yr 2000 \issue 3 \pages 595\>--\>609
\endref

\ref\Key EV1
\by \Etingof/ and \Varch/
\paper Solutions of the quantum \DYB/ and \dqg/s
\jour \CMP/ \vol 196 \yr 1998 \pages 591\)--\>640
\endref

\ref\Key EV2
\by \Etingof/ and \Varch/
\paper Exchange \dqg/s \jour \CMP/ \vol 205 \issue 1 \yr 1999
\pages 19\>--\)52
\endref

\ref\Key F
\by \Feld/
\paper Elliptic \qg/s
\inbook in Proceedings of the ICMP, Paris 1994 \ed D\&Iagolnitzer 
\yr 1995 \publ Intern.\ Press \publaddr Cambridge, MA \pages 211\)--\)218
\endref

\ref\Key FV1
\by \Feld/ and \Varch/
\paper On \rep/s of the \qg/ $E_{\tau,\eta}(sl_2)$
\jour \CMP/ \vol 181 \issue 3 \yr 1996 \pages 741\)--\)761
\endref

\ref\Key FV2
\by \Feld/ and \Varch/
\paper Elliptic \qg/s and Ruijsenaars models
\jour \JSP/ \vol 89 \issue 5\)--\)6 \yr 1997 \pages 963\>--\>980
\endref

\ref\Key FTV
\by \Feld/, \VT/ and \Varch/
\paper Solutions of elliptic \qKZBe/s and \Ba/ I
\jour \AMS/ Transl.,\ Ser\&\)2 \vol 180 \yr 1997 \pages 45\>--\)75
\endref

\ref
\by \Feld/, \VT/ and \Varch/
\paper Monodromy of \sol/s of the elliptic quantum \KZvB/ \deq/s
\jour \IJM/ \vol 10 \issue 8 \yr 1999 \pages 943\>--\>975
\endref

\ref\Key HW
\by B\&\]Y\]\&Hou and H\&\]Wei
\paper Algebras connected with the $Z\sb n$ elliptic \sol/ of the \YB/
\jour J\. Math\. Phys. \vol 30 \issue 12 \yr 1989 \pages 2750\>--\)2755
\endref

\ref\Key N
\by \MN/
\paper Yangians and Capelli identities
\jour \AMS/ Transl.,\ Ser\&\)2 \vol 181 \yr 1998 \pages 139\>--\)163
\endref

\ref\Key S
\by \Skl/
\paper On some algebraic structures related to the \YB/
\jour \FAA/ \vol 16 \issue 4 \yr 1982 \pages 263\>--\)270 
\endref

\ref
\by \Skl/
\paper On some algebraic structures related to the \YB/.
Representations of the quantum algebra
\jour \FAA/ \vol 17 \issue 4 \yr 1983 \pages 273\>--\)284 
\endref

\ref\Key TV1
\by \VT/ and \Varch/
\paper Geometry of \$q$-\hgeom/ \fn/s, \qaff/s and \eqg/s
\jour \Astq/ \vol 246 \yr 1997 \pages 1\)--\)135
\endref

\ref\Key TV2
\by \VT/ and \Varch/
\paper Difference \eq/s compatible with \tri/ \KZ/ \difl/ \eq/s
\jour \IMRN/ \yr 2000 \issue 15 \pages 801\)--\>829
\endref

\endRefs

\newpage
\nopagenumbers
\contents
\ifMag
\Entcd{}{Introduction}{1}
\Entcd{1}{Basic notation}{2}
\Entcd{2}{Elliptic quantum group $\Eqg $}{5}
\Entcd{3}{Small elliptic quantum group $\esl $}{6}
\Entcd{4}{Highest weight modules over $\esl $}{10}
\Entcd{5}{Dynamical Shapovalov form}{13}
\Entcd{6}{Contragradient modules over $\esl $ and contravariant form}{16}
\Entcd{7}{Rational dynamical quantum group $\eslr $}{19}
\Entcd{8}{Finite-dimensional highest weight modules over $\esl $}{23}
\Entcd{9}{Definition of functor $\Ec $}{27}
\Entcd{10}{Exchange quantum group $\Fsl $}{31}
\Appencd
\Entcd{\char 65}{Commutation relations in $\els $}{34}
\Entcd{\char 66}{Quantum determinant}{36}
\Entcd{\char 67}{Multiplicative forms}{39}
\Entcd{\char 68}{Proof of Theorem \[PBW]}{39}
\Entcd{\char 69}{Elliptic quantum group $\eslt $}{42}
\Entcd{\char 70}{Proof of Theorem \[Ect]}{44}
\Refcd{References}{45}
\else
\Entcd{}{Introduction}{1}
\Entcd{1}{Basic notation}{2}
\Entcd{2}{Elliptic quantum group $\Eqg $}{4}
\Entcd{3}{Small elliptic quantum group $\esl $}{5}
\Entcd{4}{Highest weight modules over $\esl $}{8}
\Entcd{5}{Dynamical Shapovalov form}{10}
\Entcd{6}{Contragradient modules over $\esl $ and contravariant form}{12}
\Entcd{7}{Rational dynamical quantum group $\eslr $}{15}
\Entcd{8}{Finite-dimensional highest weight modules over $\esl $}{18}
\Entcd{9}{Definition of functor $\Ec $}{21}
\Entcd{10}{Exchange quantum group $\Fsl $}{24}
\Appencd
\Entcd{\char 65}{Commutation relations in $\els $}{26}
\Entcd{\char 66}{Quantum determinant}{28}
\Entcd{\char 67}{Multiplicative forms}{30}
\Entcd{\char 68}{Proof of Theorem \[PBW]}{30}
\Entcd{\char 69}{Elliptic quantum group $\eslt $}{32}
\Entcd{\char 70}{Proof of Theorem \[Ect]}{33}
\Refcd{References}{34}
\fi
\endco

\bye